\documentclass[12pt]{amsart}
\usepackage{amssymb,amsfonts,amsmath,amsopn,amstext,amscd,latexsym,xy}
\usepackage{epsfig}

\input xy
\xyoption{all}

\theoremstyle{remark}

\newtheorem{para}{\bf}[subsection]

\theoremstyle{definition}

\theoremstyle{plain}

\newtheorem{thm}[para]{Theorem}
\newtheorem{lemma}[para]{Lemma}

\newtheorem{cor}[para]{Corollary}
\newtheorem{prop}[para]{Proposition}

\newenvironment{numequation}
{\addtocounter{para}{1}\begin{equation}}{\end{equation}}

\newcommand{\al}{{\alpha}}

\newcommand{\vep}{\varepsilon}
\newcommand{\vpi}{\varpi}
\newcommand{\vphi}{\varphi}

\newcommand{\cA}{{\mathcal A}}

\newcommand{\cC}{{\mathcal C}}

\newcommand{\cE}{{\mathcal E}}
\newcommand{\cF}{{\mathcal F}}

\newcommand{\cI}{{\mathcal I}}
\newcommand{\cJ}{{\mathcal J}}
\newcommand{\cJL}{{\mathcal J \mathcal L}}

\newcommand{\cL}{{\mathcal L}}
\newcommand{\cM}{{\mathcal M}}

\newcommand{\cO}{{\mathcal O}}

\newcommand{\cR}{{\mathcal R}}
\newcommand{\cS}{{\mathcal S}}

\newcommand{\cV}{{\mathcal V}}

\newcommand{\cX}{{\mathcal X}}
\newcommand{\cY}{{\mathcal Y}}

\newcommand{\bbC}{{\mathbb C}}

\newcommand{\bbF}{{\mathbb F}}
\newcommand{\bbG}{{\mathbb G}}

\newcommand{\bbP}{{\mathbb P}}
\newcommand{\bbQ}{{\mathbb Q}}

\newcommand{\bbX}{{\mathbb X}}

\newcommand{\bbZ}{{\mathbb Z}}

\newcommand{\Bx}{B^{\times}}

\newcommand{\fronr}{{\hat{\fro}^{nr}}}
\newcommand{\hFnr}{{\hat{F}^{nr}}}

\newcommand{\Fbh}{{{\bar F}^{\wedge}}}

\newcommand{\Qlb}{{\overline{\mathbb{Q}_\ell}}}

\newcommand{\froxB}{{{\mathfrak o}_B^{\times}}}

\newcommand{\froB}{{{\mathfrak o}_B}}

\newcommand{\fra}{{\mathfrak a}}

\renewcommand{\frm}{{\mathfrak m}}
\newcommand{\frn}{{\mathfrak n}}
\newcommand{\fro}{{\mathfrak o}}
\newcommand{\frp}{{\mathfrak p}}
\newcommand{\frq}{{\mathfrak q}}

\newcommand{\frx}{{\mathfrak x}}
\newcommand{\fry}{{\mathfrak y}}

\newcommand{\frE}{{\mathfrak E}}

\newcommand{\frM}{{\mathfrak M}}

\newcommand{\frP}{{\mathfrak P}}
\newcommand{\frQ}{{\mathfrak Q}}
\newcommand{\frR}{{\mathfrak R}}
\newcommand{\frS}{{\mathfrak S}}
\newcommand{\frT}{{\mathfrak T}}

\newcommand{\frX}{{\mathfrak X}}

\newcommand{\frZ}{{\mathfrak Z}}

\newcommand{\End}{{\rm End}}
\newcommand{\Ext}{{\rm Ext}}
\newcommand{\Hom}{{\rm Hom}}

\newcommand{\bksl}{\backslash}

\newcommand{\diag}{{\rm diag}}
\newcommand{\hra}{\hookrightarrow}

\newcommand{\lra}{\longrightarrow}
\newcommand{\midc}{\;|\;}

\newcommand{\ra}{\rightarrow}
\newcommand{\Spa}{{\rm Spa}}
\newcommand{\Spec}{{\rm Spec}}
\newcommand{\Spf}{{\rm Spf}}
\newcommand{\sub}{\subset}
\newcommand{\supp}{{\rm supp}}

\begin{document}

\begin{center}\Large

{\bf Deformation spaces
of one-dimensional formal modules and their cohomology}\\

\bigskip
\normalsize

{\bf Matthias Strauch}\\

\medskip

{\it Department of Pure Mathematics and Mathematical Statistics\\
Centre for Mathematical Sciences, University of Cambridge\\
Wilberforce Road, Cambridge, CB3 0WB, United Kingdom\\
e-mail: M.Strauch@dpmms.cam.ac.uk}\\

\end{center}

\bigskip
\bigskip

{\small {\bf Abstract.} Let $\cM_m$ be the formal scheme which
represents the functor of deformations of a one-dimensional formal
module over $\bar{\bbF}_p$ equipped with a level-$m$-structure. By
work of Boyer (in mixed characterististic) and Harris and Taylor,
the $\ell$-adic \'etale cohomology of the generic fibre $M_m$ of
$\cM_m$ realizes simultaneously the local Langlands and
Jacquet-Langlands correspondences. The proofs given so far use
Drinfeld modular varieties or Shimura varieties to derive this
local result. In this paper we show without the use of global
moduli spaces that the Jacquet-Langlands correspondence is
realized by the Euler-Poincar\'e characteristic of the cohomology.
Under a certain finiteness assumption on the cohomology groups, it
is shown that the correspondence is realized in only one degree.
One main ingredient of the proof consists in analyzing the
boundary of the deformation spaces and in studying larger spaces
which can be considered as compactifications of the spaces $M_m$.}\\

\normalsize

\tableofcontents

\section{Introduction}

Let $F$ be a local non-Archimedean field and denote by $\fro$ its
ring of integers. In \cite{Dr} V. G. Drinfeld introduced the
notion of a formal $\fro$-module which generalizes the concept of
a formal group. A formal $\fro$-module is a smooth formal group
over an $\fro$-algebra $R$ which is equipped with an action of
$\fro$ such that the induced $\fro$-module structure on the
tangent space comes from the morphism $\fro \ra R$. In this paper
we consider the deformations of a fixed formal $\fro$-module
$\bbX$ over the algebraic closure $\bbF$ of the residue field of
$\fro$, which has the property that multiplication by a
uniformizer $\vpi$ of $\fro$ on $\bbX$ is an isogeny of degree
$q^n$ where $q$ is the number of elements of the residue field of
$\fro$. Extending work of Lubin and Tate who considered the case
$\fro = \bbZ_p$, Drinfeld had shown that the functor of
deformations is representable by an affine formal scheme.
Moreover, he introduced the important concept of a level structure
for finite height one-dimensional formal modules and showed that
the functor which associates to some $\fronr$-algebra $R$ the set
of pairs $(X,\phi)$ of deformations $X$ of $\bbX$ to $R$, together
with a level-$m$-structure $\phi$, is representable by an affine
formal scheme $\cM_m$. To such a formal scheme one can attach its
generic fibre $M_m$, which is an analytic space over $\hFnr$ (the
completion of the maximal unramified extension of $F$). These
spaces can be considered as rigid analytic spaces, or as
non-Archimedean analytic spaces in the sense of Berkovich, or as
adic spaces in the sense of Huber. We are interested in the
$\ell$-adic \'etale cohomology groups of the spaces $M_m$, where
$\ell$ is a prime number different from the residue
characteristic of $F$. \\

{\it Historical overview.} We will give a brief historical
overview of the development of the study of these spaces and their
cohomology. Their significance, in the context of Shimura curves,
was first pointed out by H. Carayol, who studied the bad reduction
of these curves in \cite{Ca1}. In \cite{Ca2} he constructed a
representation of $GL_2(F) \times \Bx \times W_F$ in terms of
vanishing cycles at singular points in the special fibre at primes
of bad reduction of these curves. Here $B$ is the quaternion
division algebra over $F$, and $W_F$ is the Weil group. This
representation, which he calls {\it repr\'esentation locale
fondamentale}, is constructed purely locally and its existence
rests solely on the representability of the deformation spaces
with level structures, and can hence be carried out in the equal
characteristic case too. The importance of this representation is
that it realizes simultaneously the Jacquet-Langlands and the
Langlands correspondence. It is, {\it very roughly speaking}, of
the form

$$\bigoplus \pi \otimes \rho \otimes \sigma$$

\medskip

where $\pi$ is a representation of $GL_2(F)$, $\rho$ a
representation of $\Bx$ and $\sigma$ a representation of $W_F$,
and a tensor product $\pi \otimes \rho \otimes \sigma$ occurs if
and only if $\pi$ and $\rho$ correspond under the
Jacquet-Langlands correspondence {\it and} $\pi$ and $\sigma$ are
related by the Langlands correspondence (cf. \cite{Ca2} for a
precise statement). Later, in the seminal article {\it Non-abelian
Lubin-Tate theory}, cf. \cite{Ca3}, Carayol generalized the
construction of this representation to $GL_n$ and explained how to
find this representation in the cohomology of Shimura varieties.
He exhibited in fact another local representation whose
construction rests on the existence of certain non-trivial
coverings of $p$-adic symmetric spaces, which is another discovery
of Drinfeld, \cite{Dr2}. Following Carayol we call this latter
setting the {\it rigid} setting, and the former, relying on
deformation spaces of formal modules, the {\it vanishing cycle}
setting. In each case he conjectured that the construction
realizes simultaneously the Jacquet-Langlands and Langlands
correspondence\footnote{Carayol attributes these conjectures
partly to P. Deligne and V. G. Drinfeld. Furthermore, I was told
by Y. Varshavsky that Drinfeld gave the construction of the
representation and stated the conjectures in an unpublished
manuscript. For the rigid setting Drinfeld made already a less
precise conjecture in \cite{Dr2}.}, and he outlined a program how
to prove these conjectures using the cohomology of Shimura
varieties. The first to take up these conjectures was G. Faltings,
who proved in \cite{F1}, in the context of the rigid setting, that
the Euler-Poincar\'e characteristic of the cohomology of the
coverings realizes the Jacquet-Langlands correspondence. His work
is important because it does not use global arguments. On the
other hand his method does not prove much about the representation
of the Weil group (except, e.g., its dimension and compatibility
with twisting). But it explains nicely the occurrence of the
Jacquet-Langlands correspondence as being the consequence of a
Lefschetz trace formula for rigid analytic spaces. In \cite{H1} M.
Harris considered the rigid setting too, but used the theory of
$p$-adic uniformization of Shimura varieties which was shortly
before generalized by Rapoport and Zink, cf. \cite{RZ}. In this
context and as consequences of that paper he was able to prove the
existence of the local Langlands conjecture for $n < p$, cf.
\cite{H2}. Around that time the equal-characteristic case began to
be studied and it was already mentioned by Carayol that the
modular varieties of Drinfeld and Stuhler furnish the global
context in which to place the local construction. Then P. Boyer
proved Carayol's conjecture in the vanishing cycle setting for
local fields of equal characteristic, cf. \cite{Bo}, thereby
reproving the local Langlands correspondence which had been shown
earlier, without Carayol's construction but with the use of the
Drinfeld-Stuhler varieties, by Laumon, Rapoport and Stuhler.
Somewhat later Th. Hausberger showed that such modular varieties
possess at certain places rigid-analytic uniformization, and could
subsequently prove Carayol's conjecture in the rigid setting in
positive characteristic (\cite{Ha}). Finally, M. Harris and R.
Taylor accomplished Carayol's program in the vanishing cycle
setting in characteristic zero and proved for the first time the
local Langlands correspondence for a local field of characteristic
zero, cf. \cite{HT}. Among the more recent papers pursuing
investigations in this field we would like to mention T. Yoshida's
local study of the cohomology of the tame-level space $M_1$, cf.
\cite{Yo}, and Laurent Fargue's thesis \cite{Fa2}, where he
considers the cohomology of other Rapoport-Zink spaces (the spaces
$M_m$ as well as Drinfeld's coverings of $p$-adic symmetric
domains are examples of such spaces, cf. \cite{RZ}). Most recently
S. Wewers analyzed the 'stable reduction' of the deformation
spaces $M_m$ in the height-2-case, gives a precise description of
the action of the Weil group and proves Carayol's conjecture in
this case (cf. \cite{W1}, \cite{W2}).
\\

{\it The present approach.} In this paper we will study the spaces
$M_m$ and their cohomology groups by local means\footnote{For
technical reasons we need to use schemes in two instances,
algebraizing the formal schemes under consideration. The usage of
algebraizations can possibly be replaced by results on adic spaces
which we expect to hold but have not yet been proven.}. We prove
first, without the use of global modular varieties, that, for
supercuspidal representations, the Jacquet-Langlands
correspondence is realized by the Euler-Poincar\'e characteristic
of vanishing cycles in all degrees, for $F$ of arbitrary
characteristic. Then we prove that there cannot be any
supercuspidal representations except in one degree. At this step
we have to assume that the cohomology groups of certain boundary
strata are finite-dimensional. This should certainly be the case
but has not been proved yet. So we will
work under the assumption of the finiteness of these cohomology groups.\\

Following Faltings' approach, we already considered in \cite{St}
the Euler-Poincar\'e characteristic and showed that a suitable
Lefschetz trace formula for the group action of elements $(g,b)
\in GL_n(F) \times \Bx$ on the $\ell$-adic cohomology of the
spaces $M_m$ would show that the Jacquet-Langlands correspondence
is realized by that representation. However, at that time we could
only prove that such a formula holds in the case $n = 2$, using a
trace formula proved by R. Huber, cf. \cite{Hu5} (who's work also
rests on Faltings' paper \cite{F1}). Whereas it is clear that it
suffices to consider pairs $(g,b)$ with $b$ being regular
elliptic, our work was for a long time blocked by the fact that,
even for regular elliptic $b$, the pair $(g,b)$ may have fixed
points 'at the boundary' of the spaces $M_m$. To make this precise
we found out that one can in fact speak of the boundary of these
spaces, and it is a consequence of the 'quasi-compactifications'
we construct, that for pairs $(g,b)$ with {\it both} elements
being regular elliptic, a trace formula with the desired
properties is indeed provable. Due to a remark Laurent Fargues
made to me, it is actually possible to show that one can restrict
oneself to pairs $(g,b)$ with $g$ {\it and} $b$ being regular
elliptic. Fargues studied the deformation spaces too, cf.
\cite{Fa}, where he introduced generalized canonical subgroups for
$\fro$-modules of any height. Thereby he is able to define a
stratification in a neighborhood of what we call the 'boundary',
and he found out independently that regular elliptic $g$ act
without fixed points in a neighborhood of the boundary.\\

{\it Summary of content.} We are going to give a brief description
of the content of the paper. In the second section we recall the
definitions of the deformation functors, following \cite{Dr}, and
define the group actions. We show that the formal schemes $\cM_m$,
which represent the deformation functors with level structures,
are algebraizable. This uses only Artin's theorem on the
algebraization of formal moduli and no Shimura or Drinfeld-Stuhler
varieties. Then we introduce the generic fibres of the formal
schemes $\cM_m$, namely the spaces spaces $M_m$. Because of later
considerations, when we study the boundary of these spaces, we
consider $M_m$ as an adic space. The first essential ingredient of
our proof is the fixed point formula which gives a
group-theoretical expression for the number of fixed points on a
space $M_m$. This result is taken from \cite{St} and uses in a
crucial way the period map of Hopkins and Gross, which was
generalized by Rapoport and Zink for other moduli spaces of
$p$-divisible groups. For the sake of completeness we include it
here again. In section three we introduce the main tool of the
whole paper, namely the spaces $\overline{M}_m$ which we call
'quasi-compactifications', for lack of a better term. We refrain
from calling them compactifications, because they are no longer
spaces over the non-archimedean field $\hFnr$. A compactification
in the usual sense should be a space which is proper over $\hFnr$
and contains $M_m$ as a dense open subspace. However, the boundary
$\partial M_m = \overline{M}_m - M_m$ of the spaces $M_m$ is no
longer defined over $\hFnr$, it is a space where $\vpi = 0$ on the
structure sheaf. We therefore call it the {\it
characteristic-$\vpi$-boundary}. But $\overline{M}_m$ as well as
$\partial M_m$ are quasi-compact analytic adic spaces, which have
a natural interpretation in terms of formal models. Namely,
$\overline{M}_m$ can be identified with Fujiwara's Zariski-Riemann
space associated to $\cM_m$. This is the projective limit of all
admissible blow-ups of $\cM_m$. The special fibre of each such
blow-up is a projective scheme over $\bbF$, and $\overline{M}_m$
is therefore, as a topological space, a projective limit of
schemes which are proper over $\bbF$. On $\overline{M}_m$ there is
still the universal level structure $\phi: \vpi^{-m}\fro^n/\fro^n
\ra X^{univ}[\vpi^m]$ ($X^{univ}$ is the universal deformation of
$\bbX$), and we can consider the subsets $\partial_A M_m$ of
$\overline{M}_m$ where the kernel of this level structure is equal
to some direct summand $A \sub \vpi^{-m}\fro^n/\fro^n$ which is
free over $\fro/\vpi^m$. These sets define a nice stratification
of $\overline{M}_m$ by adic subspaces, and this stratification is
respected by the group action. Then it is easy to see that regular
elliptic $g$ will permute the boundary strata, mapping none of
them to itself, for $m$ sufficiently large (provided $g$ maps
$M_m$ to itself, which we assume here, for simplicity). This in
turn has the consequence that there are suitable formal models
such that the fixed point locus of $g$ on the special fibre is
contained in the 'interior', i.e. in the complement of the image
of the boundary under the specialization map. Then we can use
Fujiwara's techniques and his result on the specialization of
local terms in the trace formula, cf. \cite{Fu2}, to conclude that
the trace of $(g,b)$ on the cohomology of $M_m$\footnote{in this
introduction, when we speak of the cohomology of $M_m$ we rather
mean the Euler-Poincar\'e characteristic of the cohomology of $M_m
\times_\hFnr \Fbh$} is equal to the number of rigid-analytic fixed
points. In the forth section we use this result, together with a
d\'evissage to the regular elliptic locus, to prove that the
Jacquet-Langlands correspondence is realized on the
Euler-Poincar\'e characteristic of the cohomology, as far as
supercuspidal representations are concerned. The last step
consists then in proving that there are no supercuspidal
representations except in the middle degree. To this end we
consider another kind of compactification, which is this time a
space which lives over $\hFnr$. One can describe it as follows.
Let $\frM_m$ be an algebraization of $\cM_m$, i.e. a scheme of
finite type over $\fronr$ such that its completion at a closed
point $\frx_m$ is isomorphic to $\cM_m$. Then one associates an
analytic adic space $\frM_m^{ad}$ over $\Fbh$ to it and considers
the preimage $\overline{M}^{\vpi}_m \sub \frM_m^{ad}$ of $\frx_m$
under the specialization map. This is a pseudo-adic space which is
proper over $\Fbh$, and it contains $M_m \times_\hFnr \Fbh$ as an
open subset (which can be shown to be dense in it). We call

$$\partial^{\vpi} M_m = \overline{M}^{\vpi}_m - M_m \times_\hFnr \Fbh$$

\medskip

the {\it $\vpi$-adic boundary}. There is a canonical map from
$\overline{M}^{\vpi}_m$ to $\overline{M}_m$ and we can pull back
the previously defined strata. Combining theorems of Berkovich and
Huber, it follows that the cohomology of $\overline{M}^{\vpi}_m$
is canonically isomorphic to the cohomology of $M_m$. Because the
cohomology of $M_m$ vanishes in degree $\ge n$, in the long exact
sequence relating the cohomology with compact support of $M_m$ to
the cohomology of $\overline{M}^{\vpi}_m$ and the cohomology of
the $\vpi$-adic boundary $\partial^{\vpi} M_m$, we find
surjections

$$H^i((\partial^{\vpi} M_m) \times_\hFnr \Fbh, \Qlb)
\twoheadrightarrow H^i_c(M_m \times_\hFnr \Fbh, \Qlb) \,,$$

\medskip

if $i \ge n$. From the stratification of $\partial^{\vpi} M_m$ we
deduce that its cohomology is a successive extension of
parabolically induced representations, provided that the
cohomology groups of the strata of $\partial^{\vpi} M_m$ are
finite-dimensional. Although this should definitively be the case,
we can unfortunately not prove it yet. Hence we will work under
the assumption of finiteness of the cohomology. At the end of
section \ref{vpi-adic boundary} we will make some remarks that
should make it plausible why this assumption should be true. Then
one can conclude that there are no supercuspidal representations
in the cohomology except in degree $n-1$. In the appendix we prove
in a rather elementary way that the deformation rings are
algebraizable in the equal-characteristic case. It is based on the
fact that in this case the multiplication of $\vpi$ on the
universal deformation can be described by a polynomial. Lastly, we
show that two endomorphisms of an affinoid rigid space have the
same number of fixed points, if one of them has only finitely many
simple fixed points and if the endomorphisms
are sufficiently close to another.\\

At some instances it would have been possible to shorten the
exposition a little. This would have been the effect of a
systematic use of correspondences throughout the paper. As general
group elements $g$ will {\it not} map a space $M_m$ to itself, we
were led to consider 'intermediate' spaces $M_K$, for compact
subgroups $K \sub GL_n(F)$. If one works with correspondences,
that would have been superfluous. Another point to improve on
would be the use of {\it proper} algebraizations of $\cM_m$, i.e.
schemes $\frM_m$ which are proper over $\fronr$ and such that
$\cM_m$ is the completion of $\frM_m$ at a closed point. Then we
could have applied immediately Fujiwara's results on the
topological trace formula, instead of proving the existence of
suitable formal models explicitly. However, where it seemed to
possible for us to do so, we preferred a more direct method. \\

{\it Acknowledgements.} From what has been said above it is plain
how much this paper rests on the work of Berkovich, Drinfeld,
Carayol, Faltings, Fujiwara and Huber, to name just a few. I am
very thankful to Laurent Fargues for many discussions about the
topic of this paper. He made the crucial remark that it suffices
to consider only regular elliptic group elements, and he explained
to me how to algebraize the formal schemes. I am very grateful to
Roland Huber for equally many discussions about adic spaces, which
were always very helpful. V.G. Berkovich helped me to find a proof
of some technical result (cf. \ref{continuity of fixed points}),
and I am very thankful to him for his quick response to my
question. Part of this paper was written during a stay at IH\'ES,
Bures-sur-Yvette, and I would like to thank heartily C. Breuil and
L. Lafforgue for their invitation, and IH\'ES for its hospitality
and support. I studied this problem and learned about the material
connected with this paper while being at the department of
Mathematics of the University of M\"unster, and I am thankful to
S. Bosch and P. Schneider and their study groups for the framework
they provided. The SFB 478 "Geometrische Strukturen in der
Mathematik" at the University of M\"unster financially supported
travels to various conferences and seminars where I had the
opportunity to discuss this and related topics or give talks about
it.
\\

{\bf Notation.} In this paper, $F$ will be a non-Archimedean local
field, with ring of integers $\fro$, and $\vpi$ will be a
uniformizer of $F$. The number of elements of the residue field
will be denoted by $q$, and the residue field itself by $\bbF_q$.
The valuation $v$ on $F$ will be normalized by $v(\vpi) = 1$. This
is the only valuation which we write additively; all valuations
which appear in the context of adic spaces will be written
multiplicatively. We denote by $\bbF$ an algebraic closure of
$\bbF_q$. $\hFnr$ is the completion of the maximal unramified
extension of $F$, $\fronr$ its ring of integers, and $\Fbh \supset
\hFnr$ is the completion of an algebraic closure of $F$ whose ring
of integers we denote by $\fro_\Fbh$. If $\cX$ is a scheme over
$\fronr$ we often write $\cX \times_\fronr \fro_\Fbh$ instead of
$\cX \times_{\Spec(\fronr)} \Spec(\fro_\Fbh)$, and we use similar
abuse of notation in other instances of base change. If $R$ is a
local ring we denote by $\frm_R$ its maximal ideal. $G$ always
stands for the group $GL_n(F)$ (the number $n$ will be fixed
throughout the paper), $K_0$ is its maximal compact subgroup
$GL_n(\fro)$, and $K_m = 1+\vpi^mM_n(\fro)$ is the
$m^{\mbox{\scriptsize{th}}}$ principal congruence subgroup. $B$
will denote a central division algebra over $F$ with invariant
$\frac{1}{n}$, and $N: B \ra F$ the reduced norm.
$\ell$ denotes a prime number not dividing $q$. \\

\section{Deformation spaces with level structures}

\subsection{Deformation functors}

\begin{para}\label{deformations}
Let $\bbX$ be a one-dimensional formal group over $\bbF$ that is equipped
with an action of $\fro$, i.e. we assume given a homomorphism
$\fro \ra End_\bbF(\bbX)$ such that the action of $\fro$ on the
tangent space is given by the reduction map $\fro \ra \bbF_q \sub \bbF$.
Such an object is  called a {\it formal $\fro$-module} over $\bbF$.
Moreover, we assume that $\bbX$ is of $F$-height $n$, which means that
the kernel of multiplication by $\vpi$ is a finite group scheme of
rank $q^n$ over $\bbF$.\\

It is known that for each $n \in \bbZ_{>0}$ there exists a formal
$\fro$-module of $F$-height $n$ over $\bbF$, and that it is unique
up to isomorphism \cite{Dr}, Prop. 1.6, 1.7.\\

Let $\cC$ be the category of complete local noetherian
$\fronr$-algebras with residue field $\bbF$.
A {\it deformation} of $\bbX$ over an object $R$ of $\cC$ is
a pair $(X,\iota)$, consisting of a formal $\fro$-module $X$
over $R$ which is equipped with an isomorphism
$\iota: \bbX \ra X_{\bbF}$ of formal $\fro$-modules over $\bbF$,
where $X_{\bbF}$ denotes the reduction of $X$ modulo the maximal
ideal $\frm_R$ of $R$. Sometimes we will omit $\iota$ from
the notation.\\

Following Drinfeld \cite{Dr}, sec. 4B, we define
a {\it structure of level m} or {\it level-m-structure}
on a deformation $X$ over $R \in \cC$ ($m \ge 0$)
as an $\fro$-module homomorphism

$$ \phi: (\vpi^{-m}\fro/\fro)^n \lra \frm_R \,,$$

\medskip

such that, after having fixed a coordinate $T$ on the formal group
$X$, the power series $[\vpi]_X(T) \in R[[T]]$, which describes the
multiplication by $\vpi$ on $X$, is divisible by

$$ \prod_{a \in (\vpi^{-1}\fro/\fro)^n} (T-\phi(a)) \,.$$

\medskip

Here, $\frm_R$ is given the structure of an $\fro$-module via $X$.\\

For each $m \ge 1$ let $K_m = 1+\vpi^mM_n(\fro)$ be the
$m^{\mbox{\scriptsize{th}}}$ principal congruence subgroup inside
$K_0=GL_n(\fro)$. Define the following set-valued functor
$\cM^{(0)}_{K_m}$ on the category $\cC$. For an object $R$ of
$\cC$ put

$$\cM^{(0)}_{K_m}(R)=\{(X,\iota,\phi) \,|\, (X,\iota)
\mbox{ is a def. over $R$, $\phi$ is a
level-$m$-structure on $X$} \}/ \simeq \,,$$

\medskip

where $(X,\iota,\phi) \simeq (X',\iota',\phi')$ if and only if
there is an isomorphism $(X,\iota) \ra (X',\iota')$ of formal
$\fro$-modules over $R$, which is compatible with the level structures.
For $0 \le m' \le m$ one gets by restricting a level-$m$-structure

$$\phi:  (\vpi^{-m}\fro/\fro)^n \lra \frm_R$$

\medskip

to $(\vpi^{-m'}\fro/\fro)^n \sub (\vpi^{-m}\fro/\fro)^n$ a
level-$m'$-structure

$$\phi' = \phi|_{(\vpi^{-m'}\fro/\fro)^n}:
(\vpi^{-m'}\fro/\fro)^n \lra \frm_R \,,$$

\medskip

and hence a natural transformation

$$\cM^{(0)}_{K_m} \lra \cM^{(0)}_{K_{m'}} \,.$$

\end{para}

\begin{thm}\label{representability of def functors}
(i) The functor $\cM^{(0)}_{K_m}$ is representable by a regular
local ring $R_m$ of dimension $n$. Hence there is a universal
formal $\fro$-module $X^{univ}$ over $R_0$ which defines on the
maximal ideal $\frm_{R_m}$ of $R_m$ the structure of an
$\fro$-module. There is a universal level-m-structure

$$\phi^{univ}_m: (\vpi^{-m}\fro/\fro)^n \lra \frm_{R_m}$$

\medskip

such that, if $a_1,...,a_n$ is a basis of $(\vpi^{-m}\fro/\fro)^n$
over $\fro/(\vpi^m)$, then

$$\phi^{univ}_m(a_1),...,\phi^{univ}_m(a_n)$$

\medskip

is a regular system of parameters for $R_m$.

(ii) Let $0 \le m' \le m$. The ring homomorphism $R_{m'} \ra R_m$
which corresponds to the natural transformation
$\cM^{(0)}_{K_m} \lra \cM^{(0)}_{K_{m'}}$ described above makes
$R_m$ a finite and flat $R_{m'}$-algebra. Moreover, $R_m[\frac{1}{\vpi}]$
is \'etale and galois over $R_{m'}[\frac{1}{\vpi}]$ with Galois group
isomorphic to $K_{m'}/K_m$.

(iii) $R_0$ is (non-canonically) isomorphic to
$\fronr[[u_1,\ldots,u_{n-1}]]$.

\end{thm}

{\it Proof.} (i) This result is \cite{Dr}, Prop. 4.3.\\

(ii) That $R_m$ is finite and flat over $R_{m'}$ is again
\cite{Dr}, Prop. 4.3. For the second assertion it suffices to
treat the case $m' = 0$ (by \cite{FK}, Ch.
1, \S 1, Prop. 1.7. (3)). \\

We show first that the universal level-$m$-structure is injective.
Suppose $\phi^{univ}_m(a) = 0$ for some non-zero
$a \in (\vpi^{-m}\fro/\fro)^n$. Then there is also a non-zero $a'$ in
the subgroup

$$(\vpi^{-1}\fro/\fro)^n \sub (\vpi^{-m}\fro/\fro)^n$$

\medskip

which is mapped to zero by $\phi^{univ}_m$. But the restriction of
$\phi^{univ}_m$ to the subgroup $(\vpi^{-1}\fro/\fro)^n$ is equal
to the composition of $\phi^{univ}_1$ with the map $R_1 \ra R_m$.
If $a' \in (\vpi^{-1}\fro/\fro)^n$ is non-zero, $a'$ is part of a
basis, and hence $\phi^{univ}_1(a')$ is part of a regular system
of parameters,
in particular $\phi^{univ}_1(a') \neq 0$.\\

Denote by $X^{univ}$ the universal deformation of $\bbX$
over $R_0$. We fix a coordinate $T$ and denote for any $\al \in \fro$
by $[\al]_{X^{univ}}(T)$ the power series with coefficients in $R_0$
which describes the multiplication of $\al$ on $X^{univ}$. For a given
$m \ge 0$ write

$$[\vpi^m]_{X^{univ}}(T) = P_m(T)e_m(T) \,,$$

\medskip

with a polynomial $P_m(T) \in R_0[T]$ and a unit $e_m(T) \in R_0[[T]]$,
by the Weierstrass Preparation Theorem. Because $\bbX$ was assumed to
be of $F$-height $n$, $P_m(T)$ is a polynomial of degree $q^{nm}$
whose zeros in $R_m$ are exactly the elements $\phi^{univ}_m(a)$,
$a \in (\vpi^{-m}\fro/\fro)^n$. Therefore, as the universal
level-$m$-structure is injective, $P_m(T)$ is separable over
$Frac(R_0)$. Moreover, for any prime ideal $\frq \in R_0$
the zeros of the image of $P_m(T)$ in $Frac(R_0/\frq)[T]$
are the images of the elements $\phi^{univ}_m(a)$ in $Frac(R_0/\frq)$.
We will show
that the image of $P_m(T)$ in $Frac(R_0/\frq)[T]$ is separable if
$\frq$ does not contain $\vpi$. This means that for any such prime
ideal $\frq$ and any zero $\xi = \phi^{univ}_m(a)$ one has
$P_m'(\xi) \notin \frq$, which is equivalent to
$[\vpi^m]_{X^{univ}}'(\xi) \notin \frq$. Let us calculate this
derivative. For a given $r \in \bbZ_{>0}$ we can write

$$X^{univ}(\xi,\vpi^r) = \xi + \vpi^r\tau(\xi,r)$$

\medskip

with $\tau(\xi,r) \in 1 + \frm_{R_m}$, and the limit
$\tau(\xi) = \lim_{r \ra \infty}\tau(\xi,r)$ exists in $R_m$
and is again an
element of $1 + \frm_{R_m}$. Then we compute:

$$\begin{array}{rl}
[\vpi^m]_{X^{univ}}'(\xi) &
= \lim_{r \ra \infty}\frac{1}{\vpi^r}[\vpi^m]_{X^{univ}}(\xi + \vpi^r) \\
 & \\
 & = \lim_{r \ra \infty}\frac{1}{\vpi^r\tau(\xi,r)}
[\vpi^m]_{X^{univ}}(\xi + \vpi^r\tau(\xi,r)) \\
 & \\
 & = \lim_{r \ra \infty}\frac{1}{\vpi^r\tau(\xi,r)}
[\vpi^m]_{X^{univ}}(X^{univ}(\xi,\vpi^r)) \\
 & \\
 & = \lim_{r \ra \infty}\frac{1}{\vpi^r\tau(\xi,r)}
X^{univ}([\vpi^m](\xi),[\vpi^m](\vpi^r)) \\
 & \\
 & = \lim_{r \ra \infty}\frac{1}{\vpi^r\tau(\xi,r)}
[\vpi^m](\vpi^r) = \lim_{r \ra \infty}\frac{1}{\tau(\xi,r)}
\vpi^m = \frac{1}{\tau(\xi)} \vpi^m \,. \\
\end{array}$$

\medskip

It follows that the image of $P_m(T)$ in $Frac(R_0/\frq)[T]$ is
separable if $\frq$ does not contain $\vpi$. We are going to use
this to show
that $R_m[\frac{1}{\vpi}]$ is unramified over $R_0[\frac{1}{\vpi}]$.\\

Let us recall from the proof of \cite{Dr}, Prop. 4.3,
how the ring $R_m$ can be build up explicitely. Suppose first
that $m = 1$, and put $L_0 = R_0$. Then, for $0 \le i < n$,
there is a map

$$\varphi_i: (\vpi^{-1}\fro/\fro)^i \ra \frm_{L_i}$$

\medskip

such that, in $L_i[[T]]$, $[\vpi]_{X^{univ}}(T)$ is divisible by

$$\prod_{a \in (\vpi^{-1}\fro/\fro)^i}(T-\varphi_i(a)) \,.$$

\medskip

Put

$$f_i(T) = \frac{[\vpi]_{X^{univ}}(T)}
{\prod_{a \in (\vpi^{-1}\fro/\fro)^i}(T-\varphi_i(a))} \,, \,\, \mbox{and}
\,\, L_{i+1} = L_i[[\theta]]/(f_i(\theta)) \,.$$

\medskip

Then we have $L_n = R_1$. The power series $f_i$ all divide
the power series $[\vpi]_{X^{univ}}(T)$, and because the polynomial
$P_1(T)$ from above is separable over $R_0[\frac{1}{\vpi}]$
the successive extensions $L_i|L_{i-1}$ are \'etale
after inverting $\vpi$ (cf. \cite{Mi}, Ch. I, Ex. 3.4).
Hence $R_1[\frac{1}{\vpi}]$ is \'etale over $R_0[\frac{1}{\vpi}]$.\\

For $m > 1$ Drinfeld describes $R_m$ as follows: let
$\phi^{univ}_{m-1}$ be the universal level-($m$-$1$)-structure,
and let $a_1,...,a_n$ a basis of $(\vpi^{-(m-1)}\fro/\fro)^n$ over
$\fro/(\vpi^{m-1})$. Then

$$R_m = R_{m-1}[[y_1,...,y_n]]/
([\vpi]_{X^{univ}}(y_1) - \phi^{univ}_{m-1}(a_1),...,
[\vpi]_{X^{univ}}(y_n) - \phi^{univ}_{m-1}(a_n)) \,.$$

Write for $j = 1,...,n$

$$[\vpi]_{X^{univ}}(T) - \phi^{univ}_{m-1}(a_j) = P_{m,j}(T)e_{m,j}(T)$$

\medskip

with a polynomial $P_{m,j}(T) \in R_{m-1}[T]$ and a unit
$e_{m,j}(T) \in R_{m-1}[[T]]$. We will show that $P_{m,j}$ has
only simple zeros ouside the vanishing locus of $\vpi$. If
$P_{m,j}(\xi) = 0$ then $[\vpi^m]_{X^{univ}}(\xi) = 0$, and if
furthermore $P_{m,j}'(\xi) = 0$ then $[\vpi]_{X^{univ}}'(\xi) = 0$.
As

$$[\vpi^m]_{X^{univ}}'(T)
= ([\vpi^{m-1}]_{X^{univ}}([\vpi]_{X^{univ}}(T)))'
= [\vpi^{m-1}]_{X^{univ}}'([\vpi]_{X^{univ}}(T)) \cdot
[\vpi]_{X^{univ}}'(T) \,,$$

\medskip

$[\vpi]_{X^{univ}}'(\xi) = 0$ would imply
$[\vpi^m]_{X^{univ}}'(\xi) = 0$, and so any multiple zero of
$P_{m,j}(T)$ would be a multiple zero of $[\vpi^m]_{X^{univ}}(T)$,
which do not exist outside the vanishing locus of $\vpi$. This means
that $R_m[\frac{1}{\vpi}]$ is \'etale over $R_0[\frac{1}{\vpi}]$
for any $m$.\\

From the description of $R_m$ just recalled it follows that the
degree of $R_m$ over $R_0$ is equal to the cardinality of
$GL_n(\fro/(\vpi^m))$. On the covering $R_m[\frac{1}{\vpi}]$ over
$R_0[\frac{1}{\vpi}]$ there is an obvious action of
$GL_n(\fro/(\vpi^m))$. To prove that the covering is galois with
this group we only have to show that no non-trivial elements acts
trivially. Suppose $g \in GL_n(\fro/(\vpi^m))$ would act trivially
on $R_m[\frac{1}{\vpi}]$. Then it would follow that
$\phi^{univ}_m(g(a_i)) = \phi^{univ}_m(a_i)$ for $i= 1,...,n$.
But the universal level-$m$-structure is injective, hence $g = 1$.\\

(iii) This is \cite{Dr}, Prop. 4.2. \hfill $\Box$ \\

{\it Remark.} One can also show the statement about
\'etaleness by showing that
the finite flat group scheme of $\vpi^m$-torsion points is
unramified over the base outside the vanishing locus of $\vpi$.
This is the case because the $\fro$-action on the module of
relative differentials $\Omega$, which is induced by the structure
of an formal $\fro$-module, coincides with the usual action via
the inclusion $\fro \hra R_0$. As, on the one hand, $\vpi^m$ acts
trivial on $X^{univ}[\vpi^m]$, it acts as $0$ on $\Omega$, but, on
the other hand, $\vpi$ is invertible outside the vanishing locus
of $\vpi$,
hence the module $\Omega$ has to be zero.\\

The fact that $\fronr[[u_1,\ldots,u_{n-1}]]$ represents
$\cM^{(0)}_{K_0}$ is due to Lubin and Tate (for $F=\bbQ_p$), cf.
\cite{LT}. For this reason $\cM^{(0)}_{K_0}$, the deformation
space without level structures, is sometines called the
{\it Lubin-Tate moduli space}, cf. \cite{HG}, \cite{Ch}.\\

\subsection{Group actions}

\begin{para}\label{defs with quasi-isogenies}
Let $X$ be a formal $\fro$-module over $R \in \cC$ such that
$X_\bbF$ has $F$-height $n$, in which case we say that the
formal $\fro$-module $X$ has height $n$. As pointed out above,
$X_\bbF$ is then isomorphic to $\bbX$. Denote $End_\fro(\bbX)$
by $\froB$; this $\fro$-algebra is the maximal compact
subring of $B:= \froB \otimes_\fro F$, which is a central
division algebra over $F$ with invariant $\frac{1}{n}$,
cf. \cite{Dr}, Prop. 1.7.\\

Any non-zero element of $Hom_\fro(\bbX,X_\bbF) \otimes_\fro F$
is called an {\it $\fro$-quasi-isogeny} from $\bbX$ to $X_\bbF$.
For such an element $\iota$ we define its {\it F-height} by

$$ F\mbox{-height}(\iota) = F\mbox{-height}(\vpi^r\iota)-nr \,,$$

\medskip

where we choose some $r \in \bbZ$ such that $\vpi^r\iota$ lies in
$Hom_\fro(\bbX,X_\bbF)$, and for an element $\iota'$ of this
latter set, its $F$-height is $h$ if ker($\iota'$) is a group
scheme of rank $q^h$ over $\bbF$.\\

Define for $j \in \bbZ$ a set-valued functor $\cM^{(j)}_{K_m}$ on
$\cC$ as follows: for $R \in \cC$ the set $\cM^{(j)}_{K_m}(R)$
consists of equivalence classes of triples $(X,\iota,\phi)$, where
$X$ is a formal $\fro$-module of height $n$ over $R$, $\iota$ is
an $\fro$-quasi-isogeny from $\bbX$ to $X_\bbF$ of $F$-height $j$,
and $\phi$ is a level-$m$-structure on $X$. Now put

$$\cM_{K_m} = \coprod_{j \in \bbZ}\cM^{(j)}_{K_m} \,.$$

By the uniqueness of $\bbX$ (up to isomorphism), we have
$\cM^{(j)}_{K_m} \simeq \cM^{(0)}_{K_m}$, but there is no
distinguished isomorphism.
\end{para}

\begin{para}\label{group actions}
There is an action of $\Bx$ {\it from the right} on the functors
$\cM_{K_m}$ given by

$$[X,\iota,\phi].b = [X,\iota \circ b,\phi] \,,$$

\medskip

where we denote by $[X,\iota,\phi]$ the equivalence class of
$(X,\iota,\phi)$, and where $b \in \Bx$. If $[X,\iota,\phi]$
belongs to $\cM^{(j)}_{K_m}(R)$, then $[X,\iota,\phi].b$ is an
element of $\cM^{(j+v(N(b)))}_{K_m}(R)$,
where $N:B \ra F$ denotes the reduced norm.\\

Next we will describe the 'action' of the group $G=GL_n(F)$
on the tower $(\cM_{K_m})_m$. Let $g \in G$ and suppose first
that $g^{-1} \in M_n(\fro)$. For integers $m \ge m' \ge 0$
such that

$$ g\fro^n \sub \vpi^{-(m-m')}\fro^n$$

\medskip

(this inclusion is meant to be inside $F^n$) we will
define a natural transformation

$$g_{m,m'}:\cM_{K_m} \ra \cM_{K_{m'}} \,.$$

\medskip

Let $[X,\iota,\phi] \in \cM_{K_m}(R)$, $R \in \cC$. The following
construction gives an element $[X',\iota',\phi']$ of
$\cM_{K_{m'}}(R)$ that is the image under the corresponding map

$$(g_{m,m'})_R:\cM_{K_m}(R) \ra \cM_{K_{m'}}(R)$$

\medskip

on $R$-valued points and it will be denoted by $[X,\iota,\phi].g$.\\

The conditions imposed on $g$ show that $g\fro^n$
contains $\fro^n$ and that $g\fro^n/\fro^n$ can
naturally be regarded as a subgroup of $\vpi^{-m}\fro^n/\fro^n$,
so we may define a formal $\fro$-module $X'$ over $R$ by
taking the quotient of $X$ by the finite subgroup
$\phi(g\fro^n/\fro^n)$ (cf. \cite{Dr}, Prop. 4.4):

$$X' = X/\phi(g\fro^n/\fro^n) \,.$$

\medskip

Moreover, left multiplication with $g$ induces an injective homomorphism

$$\vpi^{-m'}\fro^n/\fro^n \stackrel{g}{\lra}
\vpi^{-m}\fro^n/g\fro^n = (\vpi^{-m}\fro^n/\fro^n)/(g\fro^n/\fro^n)$$

\medskip

and the composition with the morphism induced by $\phi$,

$$(\vpi^{-m}\fro^n/\fro^n)/(g\fro^n/\fro^n)
\lra X/\phi(g\fro^n/\fro^n) = X' \,,$$

\medskip

gives by \cite{Dr}, Prop. 4.4, a level-$m'$-structure

$$\phi': \vpi^{-m'}\fro^n/\fro^n \ra X'[\vpi^{m'}] \,.$$

\medskip

Finally define $\iota'$ to be the composition of $\iota$ with the
projection

$$X_\bbF \ra (X')_\bbF \,.$$

\medskip

This construction is independent of the
representative $(X,\iota,\phi)$ and gives a morphism of
functors. If $[X,\iota,\phi]$ lies in $\cM^{(j)}_{K_m}(R)$ then
$[X,\iota,\phi].g$ is an element of $\cM^{(j-v(detg))}_{K_{m'}}(R)$.\\

For an arbitrary element $g \in G$, choose $r \in \bbZ$ such
that $(\vpi^{-r}g)^{-1} \in M_n(\fro)$. Then, for $m \ge m' \ge 0$
with

$$\vpi^{-r}g\fro^n \sub \vpi^{-(m-m')}\fro^n$$

\medskip

and $[X,\iota,\phi] \in \cM_{K_m}(R)$, define
$[X',\iota',\phi']=[X,\iota,\phi].(\vpi^{-r}g)$ as above and put

$$[X,\iota,\phi].g = [X',\iota' \circ \vpi^{-r},\phi'] \,.$$

\medskip

This construction gives a natural transformation

$$g_{m,m'}:\cM_{K_m} \ra \cM_{K_{m'}} \,,$$

\medskip

which depends neither on $\vpi$ nor on the integer $r$
(among all $r$'s with $\fro^n \sub \vpi^{-r}g\fro^n \sub
\vpi^{-(m-m')}\fro^n$). In particular, one gets for each $m$
an action {\it from the right} of $K_0 = GL_n(\fro)$
on $\cM_{K_m}$ which commutes with
the action of $\Bx$. If $g_i \in G$, $i = 1,2$, satisfy

$$\fro^n \sub \vpi^{-r_1}g_1\fro^n \sub \vpi^{-(m_0-m_1)}\fro^n \, ,
\,\, \fro^n \sub \vpi^{-r_2}g_2\fro^n \sub \vpi^{-(m_1-m_2)}\fro^n $$

\medskip

for integers $r_1,r_2$ and $0 \le m_2 \le m_1 \le m_0$ then one has clearly

$$\fro^n \sub \vpi^{-(r_1 + r_2)}g_1g_2\fro^n \sub
\vpi^{-(m_0-m_2)}\fro^n \,.$$

\medskip

Therefore, we have morphisms

$$\begin{array}{l}
(g_1)_{m_0,m_1}: \cM_{K_{m_0}} \ra \cM_{K_{m_1}} \\
 \\
(g_2)_{m_1,m_2}: \cM_{K_{m_1}} \ra \cM_{K_{m_2}} \\
 \\
(g_1g_2)_{m_0,m_2}: \cM_{K_{m_0}} \ra \cM_{K_{m_2}} \,,
\end{array}$$

\medskip

and it is easy to check that the composition of the two former morphisms
is identical to the latter one:

\begin{numequation}\label{composition}
(g_2)_{m_1,m_2} \circ (g_1)_{m_0,m_1} = (g_1g_2)_{m_0,m_2} \,.
\end{numequation}

Note that we defined the morphisms $g_{m_i,m_j}$ so as to obtain
a right action.

\medskip
\end{para}

\begin{para}\label{formal quotients}
We need to consider also quotients of the formal schemes $\cM_{K_m}$.
For an open subgroup $K \sub K_0$ we choose $m \ge 0$ so that
$K_m \sub K$. Then let $\cM_{K_m}^{(j)} = \Spf(R_m^{(j)})$ and put

$$R^{(j)}_K = (R_m^{(j)})^K \, , \;\;
\cM_K^{(j)} = \Spf(R^{(j)}_K) \, , \;\; \cM_K = \coprod_{j \in
\bbZ} \cM_K^{(j)} \, ,$$

\medskip

where $(R_m^{(j)})^K$ denotes the subring of $K$-invariant elements in
$R_m^{(j)}$. \\
\end{para}

\begin{prop}\label{ring of invariants} Fix $j \in \bbZ$ and
$0 \le m' \le m$ such that $K_m \sub K \sub K_{m'}$. Whenever
$0 \le m_1 \le m_2$ we identify $R_{m_1}^{(j)}$ with a subring
of $R_{m_2}^{(j)}$ (cf. Theorem \ref{representability of def functors}).

(i) $(R_{m''}^{(j)})^K = (R_m^{(j)})^K$ for all $m''$ with
$K_{m''} \sub K$, and hence $R^{(j)}_K$ is well-defined. In particular,
if $K = K_{m'}$ then $R^{(j)}_K = R_{m'}^{(j)}$.

(ii) $R^{(j)}_K$ is a local noetherian ring,
which is finite over $R_{m'}^{(j)}$. It is complete with respect to
the topology defined by the maximal ideal.

(iii) $R^{(j)}_K$ is the integral closure of $R_{m'}^{(j)}$ in the
field $Frac(R^{(j)}_K)$. It is thus integrally closed. The set of
prime ideals of $R^{(j)}_K$ is in bijection with the orbits of
$K/K_m$ on the set of prime ideals of $R_m^{(j)}$:

$$\Spec(R_m^{(j)})/K \stackrel{\simeq}{\lra} \Spec(R^{(j)}_K) \,.$$

\medskip

(iv) $R^{(j)}_K[\frac{1}{\vpi}]$ is \'etale over
$R_{m'}^{(j)}[\frac{1}{\vpi}]$, and it is galois with Galois group
$K_{m'}/K$ if $K$ is normal in $K_{m'}$. $R^{(j)}_m[\frac{1}{\vpi}]$
is \'etale and galois over $R^{(j)}_K[\frac{1}{\vpi}]$
with Galois group $K/K_m$.
\end{prop}

{\it Proof.} (i) Suppose $m \ge m''$. Then $K_m \sub K_{m''} \sub K$
and

$$(R_m^{(j)})^K = ((R_m^{(j)})^{K_{m''}})^{(K/K_{m''})}
= (R_{m''}^{(j)})^K \,.$$

\medskip

(ii) Denote by $\frm$ the maximal ideal of $R_m^{(j)}$. The set of
$K$-invariant
elements $\frm^K$ in $\frm$ is a proper ideal of $R^{(j)}_K$. If now
$f \in  R^{(j)}_K$ is not in $\frm^K$, it is invertible in $R_m^{(j)}$,
$gf = 1$, say, with $g \in R_m^{(j)}$. Then, for any $k \in K$ we have
$g - k(g) = gf(g - k(g)) = g(fg - k(fg)) = 0$, hence $g \in R^{(j)}_K$,
and it follows that $R^{(j)}_K$ is local. $R^{(j)}_K$ is submodule of the
finite $R_{m'}^{(j)}$-module $R_m^{(j)}$, and as $R_{m'}^{(j)}$ is noetherian,
$R^{(j)}_K$ is finite over $R_{m'}^{(j)}$ too. This implies that
$R^{(j)}_K$ is a noetherian ring. The completeness of $R^{(j)}_K$
with respect to the topology defined by the maximal ideal follows from
the fact that the automorphisms $k \in K$ act continuously on $R_m^{(j)}$. \\

(iii) An element of $Frac(R^{(j)}_K) = Frac(R_m^{(j)})^K$ which is
integral over $R_{m'}^{(j)}$ belongs to $R_m^{(j)}$ because the latter ring
is integrally closed (as it is regular, hence factorial). So it is an
element of $R_m^{(j)} \cap Frac(R_m^{(j)})^K = R^{(j)}_K$. The assertion about
the orbits of $K$ on $\Spec(R_m^{(j)})$ is a general fact, cf. \cite{Bou},
Ch. V, \S 2.1, Thm. 1, \S 2.2, Thm. 2. \\

(iv) $R_m^{(j)}[\frac{1}{\vpi}]$ is the integral closure of
$R_{m'}^{(j)}[\frac{1}{\vpi}]$ in $Frac(R_{m'}^{(j)})$,
and by Thm. \ref{representability of def functors} it is flat
and unramified over $R_{m'}^{(j)}[\frac{1}{\vpi}]$.
Further, $R_K^{(j)}[\frac{1}{\vpi}]$ is the integral closure
of $R_{m'}^{(j)}[\frac{1}{\vpi}]$ in $Frac(R_{m'}^{(j)})$.
By \cite{Bou}, Ch. V, Ex. 19 (c) to \S 2, $R_K^{(j)}[\frac{1}{\vpi}]$
is unramified over $R_{m'}^{(j)}[\frac{1}{\vpi}]$ and
$R_m^{(j)}[\frac{1}{\vpi}]$ is unramified over $R_K^{(j)}[\frac{1}{\vpi}]$.
By \cite{FK}, Ch. 1, \S 1, Lemma 1.5., $R_K^{(j)}[\frac{1}{\vpi}]$
is \'etale over $R_{m'}^{(j)}[\frac{1}{\vpi}]$,
and by \cite{FK}, Ch. 1, \S 1, Prop. 1.7.,
$R_m^{(j)}[\frac{1}{\vpi}]$ is \'etale over
$R_K^{(j)}[\frac{1}{\vpi}]$. \hfill $\Box$

\medskip

We equip $R^{(j)}_K$ with the adic topology defined by the maximal ideal. \\

\begin{para}\label{action on formal quotients}
Consider $g \in G$ and suppose that $K' = g^{-1}Kg$ is contained in $K_0$.
Choose $0 \le m' \le m$ and $r \in \bbZ$ such that

\medskip

\begin{itemize}
\item $\fro^n \sub \vpi^{-r}g\fro^n \sub \vpi^{-(m-m')}\fro^n \,,$
\item $K_m \sub K, \, K_{m'} \sub g^{-1}Kg \,.$
\end{itemize}

\medskip

Put $d = v(\det(g))$. Then we have for any $k \in K$ and
$j \in \bbZ$ a commutative diagram (by \ref{composition})

$$\xymatrixcolsep{3pc}\xymatrix{
\cM_{K_m}^{(j)} \ar[r]^{g_{m,m'}} \ar[d]^k & \cM_{K_{m'}}^{(j-d)}
\ar[d]^{g^{-1}kg}\\
\cM_{K_m}^{(j)} \ar[r]^{g_{m,m'}} & \cM_{K_{m'}}^{(j-d)}}$$


\medskip

This shows that the ring homomorphism

$$g_{m,m'}^\sharp: R_{m'}^{(j-d)} \lra R_m^{(j)}$$

\medskip

maps the $g^{-1}Kg$-invariants in $R_{m'}^{(j-d)}$ to the
$K$-invariants in $R_m^{(j)}$, and defines hence a morphism

$$g_K^\sharp: R^{(j-d)}_{g^{-1}Kg} \lra R^{(j)}_K \,,$$

\medskip

which in turn induces a morphism of formal schemes

$$g_K: \cM_K \lra \cM_{g^{-1}Kg} \,.$$

\bigskip

\end{para}

\subsection{Algebraization of the formal schemes and the group action}

\hfill{\space} \newline

We will show that the formal schemes $\cM^{(j)}_K$ are completions
of schemes of finite type over $\fronr$ at a closed point. This
implies that the cohomology groups of the rigid-analytic spaces
attached to $\cM^{(j)}_K$ are finite-dimensional. Moreover, we
will need to know that the action of endomorphisms of these formal
schemes can be approximated. Later we will need such an
algebraicity statement when we are going to apply a result of K.
Fujiwara on the specialization of local terms.

\medskip

\begin{thm}\label{algebraicity}
Fix an open subgroup $K \sub K_0$ and an integer $j \in \bbZ$.

There is an affine scheme of finite type $\frM^{(j)}_K$ over $\fronr$
and a closed point $\frx_K$ on $\frM^{(j)}_K$ such that the following
assertions hold:

(i) there is an isomorphism of formal schemes over $\fronr$

$$\widehat{\frM^{(j)}_{K,\frx_K}} \stackrel{\simeq}{\lra} \cM^{(j)}_K$$

between the completion $\widehat{\frM^{(j)}_{K,\frx}}$ of
$\frM^{(j)}_K$ at $\frx_K$ and $\cM^{(j)}_K$;

(ii) the scheme $\frM^{(j)}_K$ carries an action of $K_0$, and
the subgroup $K$ acts trivially on $\frM^{(j)}_K$;

(iii) the point $\frx_K$ is a fixed point of the action of $K_0$,
and the isomorphism between $\widehat{\frM^{(j)}_{K,\frx_K}}$ and
$\cM^{(j)}_K$ is $K_0$-equivariant.
\end{thm}

\medskip

{\it Proof.} We may assume without loss of generality that $j = 0$,
and we identify $R^{(0)}_0$ with a power series ring in the variables
$u_1,...,u_{n-1}$ over $\fronr$. For the proof we drop the superscripts.

(a) We show first that the group scheme $X^{univ}[\vpi^m]$
is defined over some subring $\frR \sub R_0$ which is of
finite type over $\fronr$, and whose completion at a maximal ideal
is isomorphic to $R_0$. To do this we will reason along
the same lines as explained after Thm. 1.7. of \cite{A2}. In \cite{F3}
Faltings introduces what one may call
'truncated Barsotti-Tate $\fro$-modules' (of level $m$).
This is the analogue of the notion of a
truncated Barsotti-Tate group (in the sense of \cite{Me}, and \cite{Il})
in the context of group schemes with strict $\fro$-action.
Such an object is in particular of finite presentation over the base
scheme. Put $A = \fronr[u_1,...,u_{n-1}]$,
$S = \Spec(A)$, and consider
the functor $BT_m^\fro / S$ on the category of
$S$-schemes which
associates to $T / S$ the set of isomorphism classes of
truncated Barsotti-Tate $\fro$-modules of level $m$ over $T$.
It is a functor which
is locally of finite presentation over $S$. We let $\widehat A$ denote
the completion of $A$ with respect to the ideal
$\frm = (\vpi,u_1,...,u_{n-1})$, i.e. $\widehat{A} = R_0$,
and put $\widehat{S} = \Spec(\widehat{A})$.
We have $X^{univ}[\vpi^m] \in BT_m^\fro(\widehat{S})$.
By Artin's approximation theorem, \cite{A1}, Cor. 2.2, there
is an \'etale neighborhood (which we may assume to be affine)

$$\xymatrix{
 & S' = \Spec(A') \ar[d]\\
\Spec(A/\frm) \ar[r] \ar[ru] & S}$$


\medskip

of the maximal ideal $\frm$ and an element $\frX_m \in
BT_m^\fro(S')$ such that $\frX_m$ and $X^{univ}[\vpi^m]$ have the
same images in $\widehat{A}/\frm^2 \hat{A} = A'/\frm^2 A'$. Now we
use the fact that the pair $(\widehat{A},X^{univ}[\vpi^m])$ is an
effective formal versal deformation of $\bbX[\vpi^m]$ over $\bbF$
in the sense of \cite{A1}, sec. 1. In the case $\fro = \bbZ_p$
this statement is Cor. 4.8 (ii) in \cite{Il}. For arbitrary $\fro$
the versality of $\bbX[\vpi^m]$ follows from the results of
Faltings, cf. \cite{F3}, p. 276, with the same arguments as in
\cite{Il}. We conclude as explained after Thm. 1.7. of \cite{A2}:
there is an $\fronr$-linear automorphism of $\widehat A$ such that
$\frX_m \times_{S'} \widehat{S}$ is isomorphic to the pull back of
$X^{univ}[\vpi^m]$ via this isomorphism. We put $\frR = A'$ and
$\frM = S'$. This is the algebraization we are looking for. We
denote the point corresponding to the maximal ideal $\frm \sub
\frR$ by $\frx$.

(b) We define the notion of a Drinfeld level structure for
truncated Barsotti-Tate $\fro$-modules as in \cite{HT}, sec. II.2.
By \cite{HT}, Lemma II.2.1 (6.), there exists a scheme $\frM_m$
which classifies level-$m$-structures on $\frX_m \times_\frM T$
over schemes $T$ over $\frM$. As $\frM_m$ is finite over $\frM$,
it is affine, and we put $\frR_m = \cO_{\frM_m}(\frM_m)$. $\frM_m$
carries a natural action of $GL_n(\fro/\vpi^m)$ (action on the
level structure). As there is only one level-$m$-structure on
$\bbX[\vpi^m]$, there is only one point $\frx_m$ of $\frM_m$ over
$\frx$. Passing to the completion at $\frx_m$, we get that
$\widehat{\frM_{m,\frx_m}}$ classifies level-$m$-structures on
$X^{univ}[\vpi^m]$, hence is isomorphic to $\cM_m$. Finally, if $K
\sub K_0$ is an open subgroup containing $K_m$, we put
$\frR_K = (\frR_m)^K$ and $\frM_K = \Spec(\frR_K)$.
The assertions (i)-(iii) are now clear. \hfill $\Box$ \\

{\it Remark.} This proof is due to L. Fargues, cf. \cite{Fa}, sec.
9.2.1, Prop. 5, who proved an even stronger statement which gives
a better control on the action of $GL_n(\fro/\vpi^m)$ on the
'boundary' (in the sense of loc. cit., Prop. 5).\\

Whereas the action of $GL_n(\fro/\vpi^m)$ on $\cM^{(j)}_K$
extends, by our construction, to the algebraizations
$\frM^{(j)}_K$, this is not the case for the action of elements in
$\froxB$. For our purposes, it is enough to know that any given
endomorphism $\gamma$ of one of our formal schemes $\cM^{(j)}_K$
can be approximated by a suitable correspondence. More precisely,
we can find an \'etale neighborhood $\widetilde{\frM}^{(j)}_K$ of
the point $\frx_K$ (cf. Prop. \ref{algebraicity}) in
$\frM^{(j)}_K$ and a morphism from $\widetilde{\frM}^{(j)}_K$ to
$\frM^{(j)}_K$ such that such that the induced endomorphism on the
completions, both being isomorphic to $\cM^{(j)}_K$, approximate
the given endomorphism $\gamma$ up to prescribed order.

\medskip

\begin{thm}\label{approximation}
Suppose $\gamma: \cM^{(j)}_K \ra \cM^{(j)}_K$ is a morphism
of formal schemes over $\fronr$. Then, for any $c \in \bbZ_{\ge 0}$
there is an \'etale neighborhood $\widetilde{\frM}^{(j)}_K$ of the point
$\frx_K \in \frM^{(j)}_K$

$$\xymatrix{
\Spec(\kappa(\tilde{\frx}_K)) \ar[r] \ar[d]^{id} &
\widetilde{\frM}^{(j)}_K \ar[d]\\
\Spec(\kappa(\frx_K)) \ar[r] & \frM^{(j)}_K }$$


\medskip

and a morphism of schemes over $\fronr$

$$\gamma_c: \widetilde{\frM}^{(j)}_K \lra \frM^{(j)}_K $$

such that $\gamma_c(\tilde{\frx}_K) = \frx_K$, and $\gamma_c$
induces a morphism

$$\hat{\gamma_c}: \cM^{(j)}_K \lra \cM^{(j)}_K $$

such that $\hat{\gamma_c} \equiv \gamma \,\, (\frm_K^c)$ where $\frm_K$
is the maximal ideal of $R^{(j)}_K$, i.e. if

$$\gamma^\sharp, \hat{\gamma_c}^\sharp: R^{(j)}_K \lra R^{(j)}_K$$

are the corresponding ring homomorphisms, then, for all
$f \in R^{(j)}_K$

$$\gamma^\sharp(f) \equiv \hat{\gamma_c}^\sharp(f) \,\, (\frm_K^c) \, .$$
\end{thm}

{\it Proof.} This is \cite{A1}, Cor. 2.5. \hfill $\Box$

\bigskip

\subsection{Associated analytic spaces}

\hfill{\space} \newline

The next step is to introduce the analytic spaces whose
$\ell$-adic \'etale cohomology groups we are going to study. There
are different possible methods how to construct such spaces: as
rigid-analytic spaces, as non-Archimedean analytic spaces as
defined and studied by V. G. Berkovich \cite{Be1}, as
Zariski-Riemann spaces (\cite{Fu1}, \cite{Bom}) or finally as adic
spaces in the sense of R. Huber \cite{Hu1}. For each of these
kinds of spaces there has been defined an \'etale cohomology
theory (\cite{dJvdP}, \cite{Be2}, \cite{Fu2}, \cite{Hu3}) and
there are comparison theorems assuring that the resulting
cohomology groups for the spaces considered by us are the same
(\cite{Hu3}, sec. 8.3). In section \ref{The boundary of the
deformation spaces} we make use of the theory of adic or
Zariski-Riemann spaces. Therefore, we define the spaces we are
going to work with directly as adic spaces. Nevertheless, we will
now give brief references concerning the other constructions.\\

There is a construction of P. Berthelot, generalizing Raynaud's
construction for $\vpi$-adic formal schemes, which associates to
any formal scheme which is locally formally of finite type over a
discrete valuation ring $V$ a rigid analytic space (\cite{Ber},
0.2.6, \cite{RZ}, sec. 5.5) over the field of fractions. Secondly,
Berkovich has defined for so-called special formal schemes over
$V$ (these possess affine coverings of the form $\Spf(R)$ where
$R$ is a quotient of some algebra $V\langle x_1,...,x_m \rangle
[[y_1,...,y_n]]$) an associated non-Archimedean analytic space
\cite{Be3}, sec.1. Finally, R. Huber associates to a locally
noetherian formal scheme $\cX$ an adic space $t(\cX)$ (cf.
\cite{Hu2}, sec. 4, and, more generally, \cite{Hu3}, sec. 1.9). If
$\cX = \Spf(R)$ is affine, the set of points of the underlying
topological space of $t(\cX) = \Spa(R,R)$ consists of all
equivalence classes of continuous valuations $|\cdot|_v$ on $R$
such that $|f|_v \le 1$ for all $f \in R$. (We recall that all
valuations occuring in the context of adic spaces will be written
multiplicatively.) If $R = R^{(j)}_K$, the set of valuations
$|\cdot|_v$ with $|\vpi|_v=0$ is a closed subset which we denote
by $V(\vpi)$. The open complement inherits the structure of an
adic space and we put

$$M^{(j)}_K = t(\cM^{(j)}_K) - V(\vpi)\,,\,\,
\mbox{and}\,\, M_K = \coprod_{j \in \bbZ}M^{(j)}_K\,.$$

If $K \sub K'$ are open subgroups of $K_0$ there is a morphism of
adic spaces

$$M_K \lra M_{K'} \,,$$

induced from the corresponding morphism of formal schemes
$\cM_K \ra \cM_{K'}$ (cf. \ref{ring of invariants}). This morphism
is always \'etale, and it is galois with Galois group
$K'/K$ if $K$ is normal in $K'$, cf. \ref{ring of invariants} (iv).
In particular the Galois group of $M_{K_m}$ over
$M_{K_0}$ is $GL_n(\fro/(\vpi^m))$. By \ref{representability of
def functors} each space $M^{(j)}_{K_0}$ is isomorphic to an
open polydisc of dimension $n-1$; in particular:

$$M^{(j)}_{K_0}(\Fbh) \simeq \{(z_1,\ldots,z_{n-1}) \in
(\Fbh)^{n-1} \,|\, \mbox{for all } i: |z_i| < 1 \,\} \,.$$

\medskip

Moreover, we get induced group actions on the analytic spaces.
For any $g \in G$ and an open subgroup $K \sub K_0$
such that $g^{-1}Kg \sub K_0$ there is a morphism of analytic spaces

$$g: M_K \ra M_{g^{-1}Kg} \,,$$

\medskip

and these morphisms, for varying $g$ and $K$, are compatible with
each other whenever composition is defined. All the spaces $M_K$
also come with an induced action of $\Bx$ which commutes with the
morphisms induced by elements of $G$.\\

\bigskip

\subsection{$\ell$-adic cohomology and statement of the main result}

\hfill{\space} \newline

We are going to introduce the cohomology groups. We use the
\'etale cohomology theory as developed by Huber (\cite{Hu3}),
respectively Berkovich (\cite{Be2}). Because of the comparison
theorems in \cite{Hu3}, sec. 8.3, we can and will use results of
Berkovich for the \'etale cohomology of non-Archimedean analytic
spaces. From now on, we fix a prime number $\ell$ which is different
from the residue characteristic of $F$. So far, the cohomology
theories and the results concern mostly the cohomology of torsion
sheaves, and a general theory of $\ell$-adic cohomology has not
been developed yet. But for the spaces considered by us we can
show the following

\medskip

\begin{lemma}\label{finiteness of cohomology}
For any open subgroup $K \sub K_0$ and any $j \in \bbZ$, the
$\bbQ_\ell$-vector spaces

$$H^i(M^{(j)}_K \times_\hFnr \Fbh, \bbQ_\ell)
:= \left(\lim_{\stackrel{\longleftarrow}{r}} H^i(M^{(j)}_K
\times_\hFnr \Fbh, \bbZ/\ell^r\bbZ)\right) \otimes_{\bbZ_\ell} \bbQ_\ell$$

and

$$H^i_c(M^{(j)}_K \times_{\hFnr} \Fbh, \bbQ_\ell)
:= \left(\lim_{\stackrel{\longleftarrow}{r}} H^i_c(M^{(j)}_K
\times_\hFnr \Fbh,\bbZ/\ell^r\bbZ)\right)\otimes_{\bbZ_\ell} \bbQ_\ell$$

are finite-dimensional and the induced action of $\froxB$ on these
spaces is smooth. The cohomology groups $H^i(M^{(j)}_K
\otimes_{\hFnr} \Fbh, \bbQ_\ell)$ vanish for $i > n-1$, and the
cohomology groups $H^i_c(M^{(j)}_K \otimes_{\hFnr} \Fbh, \bbQ_\ell)$
vanish for $i < n-1$ and $i > 2(n-1)$.
\end{lemma}

{\it Proof.} (Cf. \cite{HT}, Lemma I.5.1) By \cite{Be3}, Prop.
2.4, the cohomology group $H^i_c(M^{(j)}_K \otimes_{\hFnr}
\Fbh,\bbZ/\ell^r\bbZ)$ is isomorphic to the stalk of
$R^i\Theta_{F^s}(\bbZ/\ell^r\bbZ)$ at the single point in the special
fibre of $\cM^{(j)}_K$, where $F^s$ denotes the separable closure
of $\hFnr$ and $\Theta_{F^s}$ is Berkovich's vanishing cycle
functor, cf. \cite{Be3}, sec. 2. By \ref{algebraicity}, the formal
scheme $\cM^{(j)}_K$ is algebraizable. By Berkovich's comparison
theorem, \cite{Be3}, Thm. 3.1, the stalk of
$R^i\Theta_{F^s}(\bbZ/\ell^r\bbZ)$ at the closed point is canonically
isomorphic to the stalk of the corresponding vanishing cycle sheaf
for the scheme $\frM^{(j)}_K$ at the closed point $\frx_K$, cf.
\ref{algebraicity}. By Nagata's theorem on compactifications,
together with Thm. 3.2 in \cite{De} and Prop. 5.3.1 in Exp. V in
\cite{SGA5}, the sheaf of vanishing cycles $R^i\Phi(\bbQ_\ell)$ on
the special fibre of $\frM^{(j)}_K$ is constructible, and hence
its stalk at $\frx_K$ is finite-dimensional. Similarly, for the
cohomology groups with compact support, we use \cite{Be3}, Cor.
2.5 and Thm. 3.1, and the
finiteness theorem Thm. 3.2 of \cite{De}.\\

The spaces $\cM^{(j)}_K$ are quasi-affine in the sense of
\cite{Be3}, Def. 5.5., and hence the statement about the vanishing
of cohomology in degree $> n-1$ (resp. in degree $< n-1$) follows
from \cite{Be3}, Thm. 6.1, Cor. 6.2. The smoothness of the action
of $\froxB$ follows from \cite{Be3}, Cor. 4.5. \hfill $\Box$

\bigskip

Next we put

$$H^i_c(M_K) = H^i_c(M_K \times_\hFnr \Fbh, \Qlb)
= \bigoplus_{j \in \bbZ} H^i_c(M^{(j)}_K \times_\hFnr \Fbh,
\bbQ_\ell)\otimes_{\bbQ_\ell} \Qlb \,.$$

On each $\Qlb$-vector space $H^i_c(M_K)$ there is an induced
action of $K_0 \times \Bx$ and for each $g \in G$ there is an
isomorphism

$$H^i_c(M_{g^{-1}Kg}) \ra H^i_c(M_K) \,,$$

as soon as $g^{-1}Kg \sub K_0$. These give rise to a
representation of $G \times \Bx$ on

$$H^i_c := \lim_{\stackrel{\lra}{K}}H^i_c(M_K) \,,$$

where the limit is taken over all compact-open subgroups $K$ of
$G$ contained in $K_0$.\\

{\it Remark.} The inertia group $Gal(F^{sep}/F^{nr})$ acts also on
$H^i_c(M_K)$, and this action can be extended to an action of the
Weil group $W_F$, cf. \cite{Bo}, Prop. 2.3.2, \cite{RZ}, sec.
3.48. Then one gets a smooth/continuous action of $G \times \Bx
\times W_F$ on $H^i_c$. In this paper however we pay only
attention to the representations of $G$ and $\Bx$.\\

The main result of this paper is the following theorem. Its
content is implied by Boyer's Theorem, \cite{Bo}, Thm. 3.2.4, in
the equal characteristic case, and it follows from the work of
Harris and Taylor \cite{HT} in the mixed characteristic case.
Whereas in these works Shimura or Drinfeld-Stuhler varieties play
a decisive role, the proof given here does not use modular
varieties. The second assertion of the following theorem is Thm.
\ref{JL on Euler characteristic}. For the last assertion, which
follows immediately from Thm. \ref{coh of boundary} (iv), we have
to assume that certain cohomology groups are finite-dimensional,
cf. hypothesis (H) in \ref{finiteness hypothesis}.

\begin{thm}\label{main theorem}
Let $\pi$ be an irreducible supercuspidal representation of $G$,
supposed to be realized over $\Qlb$.\\

(i) For each $i$ the representation $\Hom_G(H^i_c,\pi)$ of $\Bx$
is finite-dimensional and smooth.\\

(ii) In the Grothendieck group of admissible representations of
$\Bx$ one has

$$\sum_i (-1)^{i+n-1}\Hom_G(H^i_c,\pi) =  n \cdot
\cJL(\pi) \,,$$

\medskip

where $\cJL(\pi)$ is the representation of $\Bx$ associated to
$\pi$ by the Jacquet-Langlands correspondence.\\

(iii) Suppose hypothesis (H) in \ref{finiteness hypothesis} does
hold. Then, if $i \neq n-1$ one has $\Hom_G(H^i_c,\pi) = 0$, and

$$\Hom_G(H^{n-1}_c,\pi) \simeq \cJL(\pi)^{\oplus n} \,.$$

\medskip

($\cJL(\pi)^{\oplus n}$ is the direct sum of $n$ copies of
$\cJL(\pi)$.)
\end{thm}

{\it Proof of (i).} The element $\vpi$ of the center of $G$ acts
as a scalar on $\pi$, and this scalar we can write as $c^n$ for
some $c \in \Qlb$. Put $\zeta(g) = c^{-v(\det(g))}$. Then:

$$\begin{array}{rl}
\Hom_G(H^i_c,\pi) &
= \Hom_G(H^i_c \otimes \zeta , \pi \otimes \zeta)\\
 & \\
 & = \Hom_G((H^i_c \otimes \zeta) / \langle v - c^{-n} \cdot
 \vpi.v \,|\, v \in H^i_c \rangle , \pi \otimes \zeta) \,,
\end{array}$$

\medskip

where $\vpi.v$ denotes the action of $\vpi$, considered as an
element of $G$, on $v$, considered as an element of $H^i_c$.\\

Next, $(H^i_c\otimes \zeta ) / \langle v - c^{-n} \cdot \vpi.v
\,|\, v \in H^i_c \rangle$ is isomorphic, as a representation of
$G \times B^{\times}$, to the natural representation of $G \times
B^{\times}$ on

$$\left(\lim_{\stackrel{\lra}{K}}H^i_c(M_K / \vpi^{\bbZ})\right)
\otimes \xi \,,$$

\medskip

where

$$H^i_c(M_K / \vpi^{\bbZ})
= H^i_c((M_K / \vpi^{\bbZ}) \times_\hFnr \Fbh, \Qlb)$$

\medskip

and the limit is taken over all compact-open subgroups $K \sub
K_0$ of $G$, and $\xi$ is the character of $B^{\times}$ given by
$\xi (b) = c^{-v(N(b))}$. The map is defined as follows: an
element $v \in H^i_c(M^{(j)}_K, \Qlb)$ is mapped to
$c^{nk}\vpi^{-k}.v \in H^i_c(M^{(j_0)}_K, \Qlb)$, where $j = j_0 +
nk$ with $0 \le j_0 < n$. It is not difficult to check that this
map gives a $G \times B^{\times}$-equivariant isomorphism

$$(H^i_c\otimes \zeta ) / \langle v - c^{-n} \cdot \vpi.v
\,|\, v \in H^i_c \rangle \stackrel{\simeq}{\lra}
\left(\lim_{\stackrel{\lra}{K}}H^i_c(M_K / \vpi^{\bbZ})\right)
\otimes \xi \,.$$

\medskip

Hence we get the following identity of representations of
$B^\times$:

$$\Hom_G(H^i_c,\pi) = \Hom_G\left( H^i_c(M_\infty / \vpi^\bbZ),
\pi \otimes \zeta \right) \otimes \xi^{-1} \,,$$

\medskip

where

$$H^i_c(M_\infty / \vpi^\bbZ) := \lim_{\stackrel{\lra}{K}}
H^i_c((M_K / \vpi^\bbZ) \times_\hFnr \Fbh, \Qlb) \,.$$

\medskip

The representation $H^i_c(M_\infty / \vpi^\bbZ)$ is admissible if
we regard it as a representation of $G$, because if $K \sub K_0$
is an open subgroup, its subspace of $K$-invariant vectors is just
$H^i_c(M_K / \vpi^\bbZ)$, which is finite-dimensional, by
\ref{finiteness of cohomology}. Therefore, if we denote by
$H^i_c(M_\infty / \vpi^\bbZ)_{cusp}$ its cuspidal part, we have:

$$H^i_c(M_\infty / \vpi^\bbZ)_{cusp} = \bigoplus_{i \in I} \pi_i
\otimes \rho_i \,,$$

\medskip

where $I$ is some countable set (because $G/\vpi^\bbZ$ has compact
center), the representations $\pi_i$ are supercuspidal and
pairwise non-isomorphic, and $\rho_i$ is a finite-dimensional
smooth representation of $\Bx$, because $\pi_i$ has finite
multiplicity in $H^i_c(M_\infty / \vpi^\bbZ)$. This proves the
first assertion. \hfill $\Box$

\bigskip

\subsection{The period morphism and fixed points}

\begin{para}
To count fixed points we will use the period map from the moduli
spaces $M_K$ to a projective space of dimension $n-1$. This map
was first studied by M. Hopkins and B. Gross \cite{HG}. Later on,
M. Rapoport and Th. Zink introduced these morphisms for moduli
spaces for $p$-divisible groups \cite{RZ}, thereby giving a
unified account of $p$-adic period maps that have been studied
before. The set-up of Gross and Hopkins is insofar closer to our
situation as they work with formal $\fro$-modules (hence treat the
mixed and equal characteristic case simultaneously), herein
following Drinfeld. On the other hand, Gross and Hopkins only work
with one component of the moduli space $\cM_{K_0}$, namely the
component $\cM_{K_0}^{(0)} = \Spf(R_0^{(0)})$ where the
quasi-isogeny on the special fibre has height zero. After
recalling the main results of \cite{HG} in the next section, we
will explain how to define the period map on the whole space
$M_{K_0}$.
\end{para}

\medskip

\begin{para}
Let $X^{univ}$ be the universal formal $\fro$-module over the
formal scheme $\cM_{K_0}^{(0)}$, and denote by $\cE$ the universal
extension of $X^{univ}$ with additive kernel. This is a formal
$\fro$-module of dimension $n$ which sits in an exact sequence

$$0 \ra \cV \ra \cE \ra X^{univ} \ra 0 \,,$$

\medskip

where $\cV = \bbG_a \otimes
\Hom_{R_0^{(0)}}(\Ext(\cX,\bbG_a),R_0^{(0)})$. This exact sequence
furnishes an exact sequence

$$0 \ra {\cL}ie(\cV) \ra {\cL}ie(\cE) \ra {\cL}ie(X^{univ}) \ra 0 \,,$$

\medskip

of vector bundles on the formal scheme $\cM_{K_0}^{(0)}$, and an
analogous sequence

$$0 \ra {\cL}ie(\cV)^{ad} \ra {\cL}ie(\cE)^{ad} \ra {\cL}ie(X^{univ})^{ad}
\ra 0 \,,$$

\medskip

on the generic fibre of this formal scheme, i.e. on the space
$M_{K_0}^{(0)}$.
\end{para}

\medskip

\begin{prop}
([HG], Prop. 22.4, 23.2, 23.4) (i) There is a basis
$c_0,\ldots,c_{n-1}$ of ${\cL}ie(\cE)^{ad}$ such that the
$\hFnr$-subspace generated by these global sections is stable by
the action of $\froxB$. More precisely, the canonical map of
vector bundles on $M_{K_0}^{(0)}$

$$\langle c_0,\ldots,c_{n-1}\rangle_{\hFnr} \otimes \cO_{M_{K_0}^{(0)}}
\lra {\cL}ie(\cE)^{ad}$$

\medskip

is an $\froxB$-equivariant isomorphism, where $\froxB$ acts
diagonally on the left hand side. The representation of $\froxB$
on $\langle c_0,\ldots,c_{n-1} \rangle_{\hFnr}$ is equivalent to
the representation of $\froxB$ on $B \otimes_{F_n}\hFnr$ given by
left multiplication (where $F_n / F$ is the unramified
extension of degree $n$ in $\hFnr$).\\

(ii) Let $w_i$ be the image of $c_i$ in ${\cL}ie(X^{univ})^{ad}$,
$i=0,\ldots,n-1$, and denote by $W$ the space generated by these
global sections over $\hFnr$. Then the sections $w_i$ have no
common zeroes, and they are linearly independent
over $\hFnr$.\\

(iii) Denote by $\bbP(W)$ the projective space of hyperplanes in
$W$, and by $\bbP(W)^{ad}$ the associated analytic adic space.
Define

$$\pi_{K_0}^{(0)}: M_{K_0}^{(0)} \ra \bbP(W)^{ad} $$

\medskip

by sending $x \in M_{K_0}^{(0)}$ to the hyperplane

$$\{w=\al_0w_0+...+\al_{n-1}w_{n-1} \in W\otimes \hFnr(x) \,|\,
\al_0w_0(x)+...+\al_{n-1}w_{n-1}(x) = 0 \,\} \,.$$

\medskip

This map is an \'etale morphism. It is $\froxB$-equivariant and
surjective on $\Fbh$-valued points.
\end{prop}

\medskip

\begin{para}
Choose an element $\vpi_B \in \fro_B$ whose reduced norm is a
uniformizer of $F$. The action of $\Bx$ on $\cM_{K_0}$ furnishes
for each $j \in \bbZ$ an isomorphism

$$ \vpi_B^j : \cM_{K_0}^{(0)} \ra \cM_{K_0}^{(j)} \,.$$

\medskip

Define $\pi_{K_0}^{(j)} : M_{K_0}^{(j)} \ra \bbP(W)$ by
$\pi_{K_0}^{(j)} = \vpi_B^j \circ \pi_{K_0}^{(0)} \circ
\vpi_B^{-j}$. Because of the $\froxB$-equivariance of
$\pi_{K_0}^{(0)}$, this map is does not depend on the choice of
$\vpi_B$. Finally we get the {\it period map} on the whole space
$M_{K_0}$ by putting

$$\pi_{K_0} = \coprod_{j \in \bbZ}\pi_{K_0}^{(j)}: M_{K_0}
= \coprod_{j \in \bbZ} M_{K_0}^{(j)} \lra \bbP(W)^{ad} \,.$$

\medskip

More generally, for any open subgroup $K \subset K_0$ we let
$\pi_K$ be the composition of the projection $M_K \ra M_{K_0}$
with $\pi_{K_0}$, and refer to $\pi_{K}$ as a period morphism. The
proposition above gives immediately the following assertion about
the morphisms $\pi_K$.
\end{para}

\medskip

\begin{prop}\label{period_morphism}
For any open subgroup $K \subset K_0$

$$\pi_K : M_K \rightarrow \bbP(W)^{ad} $$

\medskip

is an \'etale morphism of analytic adic spaces over $\hFnr$.
Moreover, $\pi_K$ is equivariant with respect to the action of
$N_G (K) \times B^{\times}$, where the normalizer $N_G (K)$ of $K$
in $G$ acts trivially on $\bbP(W)^{ad}$ and the action of $\Bx$ on
$\bbP(W)^{ad}$ is the one that is induced by the action of $\Bx$
on $W$.
\end{prop}

\medskip

\begin{para}
Now we are in a position to count fixed points. Let $b \in \Bx$ be
an element which is regular elliptic. Hence $b$ has $n$ distinct
simple fixed points on $\bbP(W)^{ad} \times \Fbh$. Let $K$ be a
compact-open subgroup of $G$ that is contained in $K_0$, and let
$g$ be an element of the normalizer of $K$ in $G$. By Proposition
\ref{period_morphism}, the action of the pair $(g, b^{-1})$ on
$M_K \otimes \Fbh $ stabilizes the fibre of $\pi_K$ over a fixed
point of $b^{-1}$ on $\bbP(W)^{ad} \times \Fbh$. Hence we need a
description of the fibres of $\pi_K$ together with the action of
$(g, b^{-1})$. The next proposition gives such a description.
\end{para}

\medskip

\begin{prop}
(i) Let $x \in M_{K_0}(\Fbh)$, and let $[X,\iota]$ be the
deformation of $\bbX$ corresponding to $x$. Then, the fibre of
$\pi_{K_0}$ through $x$ consists of all deformations which are
quasi-isogenous to $X$. More precisely, it consists of those pairs
$[X',\iota']$ such that there exists a quasi-isogeny $f:X' \ra X$
with the property that $f_\bbF \circ \iota' = \iota$, where
$f_\bbF$
is the the reduction of $f$.\\
(ii) The fibre of $\pi_{K_0}$ through $x$ can be identified with
the set of lattices in the rational Tate module $V(X) = T(X)
\otimes_\fro F$, where

$$T(X) = \lim_\leftarrow X[\vpi^m](\Fbh) \,.$$

\medskip

By fixing an isomorphism $\phi: F^n \ra V(X)$, this set gets
identified with $G/K_0$. More generally, let $K \sub K_0$ be an
open subgroup, and let $[X,\iota,\phi]$ be a point of $M_K(\Fbh)$.
Then, the fibre of $\pi_K$ through this point can be
identified with the coset $G/K$.\\
(iii) Consider an $\Fbh$-valued fixed point of $b$ on
$\bbP(W)^{ad}$, and choose a base point of the set of
$\Fbh$-valued points of the fibre of $\pi_K$ over this point.
Using this fixed point, identify this set with $G/K$, as in (ii).
Then there exists $g_b \in G$ with the same characteristic
polynomial as $b$ such that the action of $(g, b^{-1})$ on the
(set of $\Fbh$-valued points of the) fibre is given, in terms of
this identification, by

$$ hK \longmapsto g_bhg K \,.$$

\medskip

\end{prop}

{\it Proof.} The first assertion follows from Prop. 23.28 of
\cite{HG}. The relationship between lattices in the rational Tate
module and quasi-isogenies in the mixed characteristic case can be
found in Lubin's paper \cite{Lu}, Theorem 2.2. The same holds true
also in the equal characteristic case,
cf. \cite{Yu}, sec. 3. The second assertion of (ii) follows immediately.\\

Now we are going to prove part (iii). Fix an $\Fbh$-valued point
of $M_K$, given by a triple $[X,\iota,\phi]$. We can consider
$\phi$ as an isomorphism $\fro^n \ra T(X)$ which is determined up
to multiplication (from the right) by elements from $K$. Suppose
this point is mapped by $\pi_K$ onto a fixed point of $b$. Then it
follows from \cite{HG}, Prop. 23.28, that $b$ lifts to an
endomorphism $\tilde{b}: X \rightarrow X$ of the formal
$\fro$-module $X$ such that $\tilde{b}_\bbF \circ \iota = \iota
\circ b$, where $\tilde{b}_\bbF$ is the quasi-isogeny induced on
the special fibre. $\tilde{b}$ is mapped to $b$ under the
canonical map $\End_\fro(X) \otimes F \hookrightarrow
End_\fro(\bbX) \otimes F$. Therefore the caracteristic polynomial
of $\tilde{b}$ is the same as that of $b$. Let $g_b \in G$ be such
that the following diagram is commutative:

$$\begin{array}{rcl}
F^n & \stackrel{\phi}{\longrightarrow}& V(X) \\
\mbox{\scriptsize{$g_b$}}\downarrow & &
\downarrow \mbox{\scriptsize{$V(\tilde{b})$}} \\
F^n & \stackrel{\phi}{\longrightarrow} & V(X)
\end{array}$$

\medskip

Let $[X',\iota',\phi']$ be an element in the fibre of $\pi_K$.
Hence there is a quasi-isogeny $f:X' \ra X$ and an element $h \in
G$ such that the following diagram commutes:

$$\begin{array}{rcl}
F^n & \stackrel{\phi'}{\longrightarrow}& V(X') \\
\mbox{\scriptsize{$h$}}\downarrow & &
\downarrow \mbox{\scriptsize{$V(f)$}} \\
F^n & \stackrel{\phi}{\longrightarrow} & V(X)
\end{array}$$

\medskip

The class $hK \in G/K$ corresponds to the point
$[X',\iota',\phi']$. This point is mapped by $b^{-1}$ to
$[X',\iota' \circ b^{-1}, \phi']$. The map $\tilde{b} \circ f : X'
\ra X$ is then a quasi-isogeny, if we equip $X'$ with the map
$\iota' \circ b^{-1}: \bbX \ra (X')_\bbF$. Moreover, we see that
the following diagram commutes:

$$\begin{array}{rcl}
F^n & \stackrel{\phi'}{\longrightarrow}& V(X') \\
\mbox{\scriptsize{$g_bh$}}\downarrow & & \downarrow
\mbox{\scriptsize{$V(\tilde{b})V(f)$}} \\
F^n & \stackrel{\phi}{\longrightarrow} & V(X)
\end{array}$$

\medskip

The action of $b^{-1}$ on the fibre of $\pi_K$ is thus given by
sending $hK$ to $g_bhK$. It is straightforward to check that the
action of some $g \in N_G(K)$ on this fibre is given by sending
$hK$ to $hgK$.
This proves the third assertion. \hfill $\Box$\\

\begin{thm}\label{fixedpointtheorem}
Let $g_b$ be in the conjugacy class corresponding to $b$. Then the
number of (intersection theoretic) fixed points of the pair
$(g,b^{-1})$ on the space $(M_K/\vpi^\bbZ) \times_\hFnr \Fbh$ is
equal to

$$n \cdot \# \{h \in G / \vpi^\bbZ K \,|\, h^{-1} g_b h = g^{-1}\} \,,$$

\medskip

and each fixed point is simple. The identity $h^{-1} = g_b h g^{-1}$
means that for some (and hence any) representative
$\dot{h}$ of $h \in G / \vpi^\bbZ K$ the cosets $\dot{h}^{-1} g_b
\dot{h} \vpi^\bbZ K$ and $g^{-1} \vpi^\bbZ K$ are equal. By the
fact pointed out in \ref{characters of cuspidals}, the number of
such $h \in G/\vpi^\bbZ K$ is always finite.
\end{thm}

{\it Proof.} Let $b \in B^{\times}$ be regular elliptic and
consider the fibre of the induced map

$$(M_K / \varpi^{\bbZ})(\Fbh) \rightarrow \bbP(W)(\Fbh)$$

\medskip

over a fixed point of $b^{-1}$. By the preceding proposition, we may
identify this set with $G / \vpi^\bbZ K$ and the action of $(g,
b^{-1})$, $g$ in the normalizer of $K$ in $G$, is given by

$$h \vpi^\bbZ K \longmapsto g_b h g \vpi^\bbZ K \,,$$

\medskip

where $g_b \in G$ has the same characteristic polynomial as $b$.
Hence the number of fixed points on such a fibre is

$$\# \{h \in G / \vpi^\bbZ K \midc h^{-1} g_b h = g^{-1} \} \,.$$

\medskip

Because there are $n$ simple fixed points of $b$ on $\bbP(W)$ and
the morphism $\pi_K$ is \'etale, all fixed points are simple and
the total number of fixed points is

$$ n \cdot \# \{h \in G / \vpi^\bbZ K \midc h^{-1} g_b h = g^{-1} \} \,.$$

\hspace*{\fill} $\Box$

\medskip

\section{The characteristic-$\vpi$-boundary}
\label{The boundary of the deformation spaces}

\subsection{Quasi-Compactifications}

\begin{para}\label{compactifications}
For any open subgroup $K \sub K_0$ and integer $j \in \bbZ$ we
consider the adic spaces

$$\overline{M}^{(j)}_K = t(\cM^{(j)}_K)_a
=  \Spa(R^{(j)}_K, R^{(j)}_K) - V(\frm_{R^{(j)}_K}) \,,$$

\medskip

where $V(\frm_{R^{(j)}_K})$ is the one-point-set consisting of the
single valuation of $R^{(j)}_K$ which factorizes over the residue
field. This valuation is also the single non-analytic point of
$\Spa(R^{(j)}_K, R^{(j)}_K)$. We put

$$\overline{M}_K = \coprod_{j \in \bbZ}\overline{M}^{(j)}_K \,.$$

\medskip

The group action extends to these spaces. Namely, if $g \in G$ is
such that $g^{-1}Kg \sub K_0$, then the induced map $\cM_K \ra
\cM_{g^{-1}Kg}$ induces a morphism of adic spaces

$$\overline{M}_K \lra \overline{M}_{g^{-1}Kg} \,.$$

\medskip

We will show below that these spaces contain the previously
defined spaces $M^{(j)}_K$ and $M_K$ as open subspaces, and we
denote by

$$\partial M^{(j)}_K = \overline{M}^{(j)}_K - M^{(j)}_K \,, \,\,
\partial M_K = \overline{M}_K - M_K$$

\medskip

their complements. These complements have natural stratifications.
In order to introduce them, we denote for a given integer $h$, $0
\le h \le n-1$, and $m > 0$, by $\cS_{m,h}$ the set of direct
summands $A \sub (\vpi^{-m}\fro/\fro)^n$ which are free over
$\fro/(\vpi^m)$ of rank $h$. If $K \sub K_0$ is an open subgroup,
fix $m \in \bbZ_{> 0}$ such that $K$ contains $K_m$. Then $K$ acts
(from the left) on $\cS_{m,h}$, and we put

$$\cS_{K,h} = K \bksl \cS_{m,h} \,, \,\,
\cS_K = \coprod_{0 \le h \le n-1} \cS_{K,h} \,.$$

\medskip

The so defined sets do not depend on $m$. Note that even for $K = K_0$
we have chosen $m>0$ so that $\cS_{K_0}$ is the disjoint
union of $n$ one-element-sets. Consider the universal
level-$m$-structure

$$\phi^{univ}_m: (\vpi^{-m}\fro/\fro)^n \ra X^{univ}[\vpi^m] \,,$$

\medskip

where we consider the universal deformation $X^{univ}$ as a formal
$\fro$-module over $\overline{M}_{K_0}$, which is possible,
because the universal deformation is already an object over the
schemes $\Spec(R^{(j)}_{K_0})$ (cf. Prop. \ref{representability of
def functors}). By base change we consider it as an object over
any of the spaces $\overline{M}_K$, $K \supset K_m$. For any point
$x \in \overline{M}_{K_m}$, we can consider the formal
$\fro$-module

$$X^{univ} \otimes \kappa(x)$$

\medskip

over the residue field $\kappa(x)$ at that point. We can also
consider the associated $\vpi$-divisible module, i.e. the
ind-group scheme

$$(X^{univ} \otimes \kappa(x))[\vpi^\infty]
 = \lim_{\stackrel{\lra}{r}} \, (X^{univ} \otimes \kappa(x))[\vpi^r] \,
.$$

\medskip

The \'etale part of this ind-group scheme has $F$-rank $n-h$ for
some $h$ between $0$ and $n-1$. (The point which corresponds to
the maximal ideal is not in our spaces, and this is the only point
where the universal deformation is connected.) It follows that the
kernel of the induced level-$m$-structure

$$\phi^{univ}_m: (\vpi^{-m}\fro/\fro)^n \ra (X^{univ}\otimes
\kappa(x))[\vpi^m]$$

\medskip

is a free submodule $A_x \sub (\vpi^{-m}\fro/\fro)^n$ of rank $h$.
More generally, if $x \in \overline{M}_K$, there is always a point
$x' \in \overline{M}_{K_m}$ lying over $x$ (cf. \cite{En}, 13.2),
and if $x',x''$ both lie over $x$, then there is a $k \in K$ such
that $x'k = x''$ (cf. \cite{En}, 14.1). As the group action is via
the level structure, it is easily seen that for the kernels of the
induced level structures at $x'$ and $x''$ one has:

$$A_{x''} = k^{-1} A_{x'} \,.$$

\medskip

Hence we define the {\it kernel of the universal level structure
at $x$} to be the class of $A_{x'}$ in $\cS_{K,h}$ for any $x' \in
\overline{M}_{K_m}$ lying over $x$. Then, for a given element $A
\in \cS_{K,h}$ we denote by

$$\partial_A M_K \sub \overline{M}_K \,, \,\,
\partial_A M^{(j)}_K \sub \overline{M}^{(j)}_K$$

\medskip

the set of points $x$ of $\overline{M}_K$ (resp.
$\overline{M}^{(j)}_K$) such that the kernel of the universal
level-$m$-structure at $x$ is equal to $A$. (Again we point out
that $m$ was chosen to be positive.) Then we have a
(set-theoretical) decomposition of $\overline{M}_K$:

$$\overline{M}_K = \coprod_{A \in \cS_K}\partial_A M_K \,.$$

\medskip

The group action respects this decomposition. For $g \in G$
choose $r \in \bbZ$ and $m \ge m' \in \bbZ_{\ge 0}$ such that

$$\fro^n \sub \vpi^{-r}g\fro^n \sub \vpi^{-(m-m')}\fro^n \,,$$

\medskip

and $K_{m'} \sub g^{-1}Kg \sub K_0$. Put $g_1 = \vpi^{-r}g$. Then
we have the following commutative diagram:

\begin{multline}\label{action on boundary labels}
$$\xymatrixcolsep{4pc}\xymatrix{
 &  \vpi^{-m}\fro^n/\fro^n \ar[d] \ar[rd]^{\vpi^{(m-m')}} & \\
\vpi^{-m'}\fro^n/\fro^n \ar[r]^{g_1} & \vpi^{-m}\fro^n/g_1\fro^n
\ar[r]^{\vpi^{(m-m')}} & \vpi^{-m'}\fro^n/\fro^n}$$
\end{multline}

\medskip

The composition of the arrows in the bottom line is an isomorphism
which we denote by $g'$. Now consider an element $A \in
\cS_{K,h}$, and choose a representative $A_1 \sub
\vpi^{-m}\fro^n/\fro^n$. Then $\vpi^{(m-m')}A_1 \sub
\vpi^{-m'}\fro^n/\fro^n$ is a direct summand, and
$(g')^{-1}\vpi^{(m-m')}A_1$ is an element of $\cS_{m',h}$ whose
class in $\cS_{g^{-1}Kg,h}$, which we denote by $g^{-1}A$, does
not depend on the class of $A_1$ in $\cS_{K,h}$. Therefore, we
have defined in this way a well-defined map

$$\cS_{K,h} \lra \cS_{g^{-1}Kg,h} \,, \,\, A \mapsto g^{-1}A \,.$$

\medskip

It is easy to check that the morphism $g: \overline{M}_K \ra
\overline{M}_{g^{-1}Kg}$ induces maps

$$g: \partial_A M^{(j)}_K \lra
\partial_{g^{-1}A} M^{(j-v(det(g)))}_{g^{-1}Kg}
\,.$$

\medskip

If $K = K_m$, we identify elements of $\cS_K$ with submodules of
$(\vpi^{-m}\fro/\fro)^n$, and hence there is a partial ordering on
this set defined by $A_1 \prec A_2$ if and only if $A_1 \supset
A_2$. For $K \subset K_m$ as above we define analogously $A_1
\prec A_2$ for classes $A_1, A_2 \sub \cS_K$ if and only if there
are representatives $A_1', A_2' \in \cS_{K_m}$ such that $A_1'
\prec A_2'$. This defines a partial ordering on $\cS_K$.
Note that for this ordering the class of the zero module is the
{\it maximal} element.
\end{para}

\medskip

\begin{prop}\label{top prop of comp}
Let $K \sub K_0$ be an open subgroup.

(i) For any $j \in \bbZ$ and any $A \in \cS_K$ there is a prime
ideal $\frp_A \sub R^{(j)}_K$ with the following property:
$\partial_A M^{(j)}_K$ consists of exactly those valuations $v \in
\overline{M}^{(j)}_K$ such that the support of $v$ contains
$\frp_A$ but does not contain $\frp_{A'}$ for any $A' \prec A$,
$A' \neq A$. With the notation of adic spaces we have:

$$\partial_A M_K = \Spa(R^{(j)}_K/\frp_A,R^{(j)}_K/\frp_A)_a
- \bigcup_{A' \prec A, A' \neq A} V(\frp_{A'}) \,.$$

\medskip
If $A \neq A'$ then $\partial_A M_K \cap \partial_{A'} M_K = \emptyset$.

(ii) For any element $A \in \cS_K$ the subspace $\partial_A M_K
\sub \overline{M}_K$ is open in its closure, and this closure is
equal to the union of all $\partial_{A'} M_K$ with $A' \prec A$.
\end{prop}

{\it Proof.} It suffices to prove this for the case of a principal
congruence subgroup $K = K_m$, $m > 0$. The general case follows
by passage
to the quotient.\\

(i) Consider the universal level-$m$-structure

$$\phi^{univ}_m: (\vpi^{-m}\fro/\fro)^n \lra \frm_{R^{(j)}_K} \,.$$

\medskip

Let

$$\frp_A = (\phi^{univ}_m(a))_{a \in A}$$

\medskip

be the ideal generated by the images of the elements of $A$ under
$\phi^{univ}_m$. The points in $\partial_A M_K$ are clearly those
whose support contains $\frp_A$, but does contain $\frp_{A'}$ for
any $A' \prec A$, $A' \neq A$. To show that $\frp_A$ is a prime
ideal we may assume (as the situation is homogenous under the
group action) that $A$ is generated by the first $h$ standard
basis vectors of $(\vpi^{-m}\fro/\fro)^n$. Then $\frp_A$ is
generated by the images of the first $h$ standard basis vectors
under the universal level-$m$-structure. By \ref{representability
of def functors}, these form part of a system of parameters, and
hence generate a prime ideal.\\

(ii) We show more generally that if $R$ is a noetherian local
complete domain, and if $\{0\} \subsetneq \fra$ is any non-zero
ideal of $R$, then $\Spa(R,R)_a - V(\fra)$ is an open dense subset
of $\Spa(R,R)_a$. Firstly, $V(\fra)$ is the intersection of the
complements of the sets $\{v \midc |f|_v \le |f|_v \neq 0\}$ which
are open in $\Spa(R,R)_a$. Hence $V(\fra)$ is a closed subset.
Suppose that $U$ is an open subset of $\Spa(R,R)_a$ which is
contained in $V(\fra)$. We may assume that $U$ is a rational
subset of the form

$$\{v \in \Spa(R,R)_a \midc |f_i|_v \le |s|_v \neq 0\
\mbox{ for } i = 1,\ldots,n\}$$

\medskip

with elements $f_1,\ldots,f_n \in R$ generating an open ideal. We
consider the corresponding ring of functions on $U$ which is
$R\langle\frac{f_1}{s},\ldots,\frac{f_n}{s}\rangle$. As $U$ was
supposed to lie in $V(\fra)$, the images of the elements of $\fra$
in $R\langle\frac{f_1}{s},\ldots,\frac{f_n}{s}\rangle$ vanish. By
its very definition,
$R\langle\frac{f_1}{s},\ldots,\frac{f_n}{s}\rangle$ is the
completion of $R[\frac{1}{s}]$ with respect to the topology having
the sets $\frm^r_R R[\frac{f_1}{s},\ldots,\frac{f_n}{s}]$, $r \ge
0$, as a fundamental system of neighborhoods of zero. Suppose $U$
is non-empty, and let $|\cdot|$ be a continuous valuation of $R$
in $U$. Then $s \neq 0$ and we have $|f| < 1$ for all elements $f
\in \frn = \frm_R R[\frac{f_1}{s},\ldots,\frac{f_n}{s}]$. Hence
$\frn$ is a proper ideal in
$R[\frac{f_1}{s},\ldots,\frac{f_n}{s}]$. As this ring is a
noetherian domain, the intersection $\cap_r \frn^r$ is the zero
ideal. Therefore, the elements in $\fra$ already map to zero in
$R[\frac{1}{s}]$, which is impossible if $R$ is a domain and
$\fra$ is non-zero. \hfill $\Box$

\bigskip

\subsection{Consequences for formal models}

\begin{para}\label{set-up consequences}
In this subsection we fix a regular elliptic element $g \in G$ and
an element $b \in \Bx$. We assume $v(det(g)) + v(N(b)) = 0$, so
that the action of the pair $(g,b^{-1})$ does not change the
height of the quasi-isogeny. This implies that $(g,b^{-1})$ acts
on the space $M^{(j)}_K$ as soon as $g$ normalizes $K$. We are
going to study the fixed point locus of this pair, and begin with
the following
\end{para}

\medskip

\begin{lemma}\label{appropriate subgroups}
(i) There is a fundamental system of neighborhoods of $1$ in G
consisting of compact open subgroups $K \sub K_0$ which are
normalized by $g$.

(ii) There is an open subgroup $K' \sub K_0$, depending only on
$g$, with the property that if $K$ is an open subgroup of $K'$ and
normalized by $g$, then $(g,b^{-1})$ maps for any $A \in
\cS_{K,h}$ the stratum $\partial_A M^{(j)}_K$ to
$\partial_{g^{-1}A} M^{(j)}_K$ and the intersection of these sets
is empty if $h > 0$. In particular, $(g,b^{-1})$ does not have any
set-theoretical fixed points on $\partial M_K$.
\end{lemma}

{\it Proof.} (i) (Cf. \cite{He2}, sec. 6.1) The ring $E = F[g] \sub
M_n(F)$ is a separable field extension of $F$ of degree $n$.
Denote by $\fro_E$ its ring of integers, and by $\frp_E \sub
\fro_E$ its maximal ideal. We identify $G$ with
$\End_F(E)^\times$. For $r \in \bbZ_{> 0}$ put

$$K(r) = 1 + \bigcap_{r' \ge 0}
\Hom_\fro(\frp_E^{r'},\frp_E^{r+r'}) \,.$$

\medskip

These subgroups form a fundamental system of neighborhoods of $1$
in $G$, and a simple calculation shows that $g$ normalizes $K(r)$.\\

(ii) Fix $h$ with $1 \le h \le n-1$. We can identify the set
$\cS_{K,h}$ of labels of boundary components of
$\overline{M}^{(j)}_K$ with $K \bksl K_0/P(\fro)$, where $P(\fro)$
is the stabilizer of $\fro^h \sub \fro^n$. The quotient
$K_0/P(\fro)$ is the same as $G/P$, where $P$ is the stabilizer of
$F^h \sub F^n$, and the action by $g$ on $\cS_{K,h}$, cf. diagram
\ref{action on boundary labels}, is given on $K \bksl G/P$ by
multiplication by $g^{-1}$ from the left:

$$KxP \mapsto (g^{-1}Kg)g^{-1}xP = Kg^{-1}xP \,.$$

\medskip

An element of $G$ is regular elliptic if and only if its
characteristic polynomial is separable and irreducible over $F$.
It follows that the set of regular elliptic elements is open in
$G$. Therefore, if $K$ is sufficiently small and normalized by
$g$, cf. \ref{appropriate subgroups}, the set $Kg^{-1}$ will have
empty intersection with any conjugate of $P$. Now, if $KxP$ would
be a fixed point of the action of $g$, we would have $kg^{-1} \in
xPx^{-1}$ for some $k \in K$, contradicting what we just stated.
Hence there are no fixed points of $g$ on $\cS_{K,h}$. As
$(g,b^{-1})$ maps $\partial_A M^{(j)}_K$ to $\partial_{g^{-1}A}
M^{(j)}_K$ (the action of $b$ leaves the boundary components
unchanged, because $b$ does not act on the level structure), and
as different boundary components have empty intersection, we have
proved the second assertion. \hfill $\Box$

\medskip

Let $K \sub K_0$ be an open subgroup normalized by $g$ which has
the property that $\gamma = (g,b^{-1})$ has no fixed points on the
boundary of $M^{(j)}_K$. We want to show that there is a formal
model for the analytic space $\overline{M}^{(j)}_K$, i.e. an
admissible blow-up (in the sense of \cite{Fu1}, Def. 4.1.1) of
$\cM^{(j)}_K = \Spf(R^{(j)}_K)$, to which the action of $\gamma$
extends and such that the fixed point locus on the special fibre
does not meet the image of the boundary under the specialization
map. We recall the definition of the specialization map in the
affine case (cf. \cite{Hu2}, Prop. 4.1). Let $\cX = \Spf(A)$ be an
affine noetherian formal scheme and $v \in t(\cX) = \Spa(A,A)$ a
continuous valuation. Then the set

$$sp_\cX(v) = \{f \in A \midc |f|_v < 1 \}$$

\medskip

is an open prime ideal of $A$, hence a point in $\cX = \Spa(A)$.
The so defined map

$$sp_\cX: t(\cX) = \Spa(A,A) \lra \cX =  \Spf(A)$$

\medskip

is continuous, and even a morphism of topologically ringed spaces
(\cite{Hu2}, Prop. 4.1). We denote its restriction to the subspace
$t(\cX)_a \sub t(\cX)$ of analytic points by $sp_\cX$ as well. In
order to deduce the existence of formal models from topological
properties of the adic space, we are going to use the theory of
Fujiwara's Zariski-Riemann spaces. The sheaf $\cO^+_{t(\cX)_a}$
which appears below has been defined in \cite{Hu2}, sec. 1.

\medskip

\begin{prop}\label{adic and ZR spaces}
If $\cX$ is a noetherian formal scheme, then the locally ringed
space

$$(t(\cX)_a, \cO^+_{t(\cX)_a})$$

\medskip

having the open subspace $t(\cX)_a \sub t(\cX)$ as underlying
topological space, which is equipped with the sheaf
$\cO^+_{t(\cX)_a}$ is canonically isomorphic, via the
specialization maps, to the projective limit of all admissible
blow-ups of $\cX$:

$$\lim_{\longleftarrow} sp_{\cX'}: (t(\cX)_a, \cO^+_{t(\cX)_a})
\stackrel{\sim}{\lra} \lim_{\longleftarrow} \cX' \,,$$

\medskip

where $\cX'$ runs through the set of admissible blow-ups of $\cX$
and the projective limit on the right carries the projective limit
topology and structure of a locally ringed space (cf. \cite{Fu1},
4.1.3).
\end{prop}

{\it Proof.} If $\cX$ is of topologically finite type over a
discrete valuation ring, then the assertion is known by
\cite{vdPS}. For the general case we refer to \cite{Hu0}, 3.9.25.
\hfill $\Box$

\medskip

We come back to the assertion about formal models that we want to
prove.

\medskip

\begin{prop}\label{appropriate models}
Let $\gamma = (g,b^{-1}) \in G \times \Bx$ be as in \ref{set-up
consequences}. For any open subgroup $K'$ of $K_0$ there exists an
open subgroup $K \sub K'$ which is normalized by $g$, and such the
following assertions do hold:\\

(i) There exists an admissible blow-up $(\cM^{(j)}_K)'$ of
$\cM^{(j)}_K$ to which the action of $\gamma$ extends, and such
that the fixed point locus of the extended morphism

$$\gamma': (\cM^{(j)}_K)' \lra (\cM^{(j)}_K)'$$

\medskip

has no set-theoretical fixed points on the image of $\partial
M^{(j)}_K$ under the map $sp_{(\cM^{(j)}_K)'}$.

(ii) Denote by $\frM^{(j)}_K$ the scheme of finite type over
$\fronr$ from Thm. \ref{algebraicity}. For a given integer $c \in
\bbZ_{\ge 0}$ denote by $\widetilde{\frM}^{(j)}_K$ the \'etale
neighborhood of the point $\frx_K$ and by $\gamma_c$ the morphism

$$\gamma_c: \widetilde{\frM}^{(j)}_K \lra \frM^{(j)}_K$$

\medskip

from Thm. \ref{approximation}. Then the blow-up $(\cM^{(j)}_K)'$
of $\cM^{(j)}_K$ in (i) induces a blow-up $(\frM^{(j)}_K)'$ of
$\frM^{(j)}_K$, and $\gamma_c$ lifts to a blow-up
$(\widetilde{\frM}^{(j)}_K)'$ of $\widetilde{\frM}^{(j)}_K$:

$$\xymatrix{
(\widetilde{\frM}^{(j)}_K)' \ar[r]^{\gamma_c'} \ar[d]^{\tilde{pr}}
& (\frM^{(j)}_K)' \ar[d]^{pr}\\
\widetilde{\frM}^{(j)}_K \ar[r]^{\gamma_c} & \frM^{(j)}_K }$$

\medskip

Now, if $\cJ$ is any ideal of definition of $(\cM^{(j)}_K)'$, we
can choose $c$ sufficiently large, so that the completion of
$(\widetilde{\frM}^{(j)}_K)'$ along the closed subscheme
$\tilde{pr}^{-1}(\tilde{\frx}_K)$ will be equal to
$(\cM^{(j)}_K)'$ and the induced morphism

$$\hat{\gamma}_c': (\cM^{(j)}_K)' \lra (\cM^{(j)}_K)'$$

\medskip

is congruent to $\gamma'$ modulo $\cJ$. Moreover, if $c$ is large
enough, the correspondence $\gamma_c'$ will also have no
set-theoretical fixed points on $sp_{(\cM^{(j)}_K)'}(\partial
M^{(j)}_K)$.
\end{prop}

{\it Proof.} (i) We note first that the set of admissible blow-ups
to which the action of $\gamma$ lifts is cofinal in the set of all
admissible blow-ups. Namely, the maximal ideal $\frm_R$ of
$R^{(j)}_K$ is mapped by

$$\gamma^\sharp: R^{(j)}_K \lra R^{(j)}_K$$

\medskip

to itself, as $\gamma^\sharp$ is continuous. Furthermore, the
element $\gamma$ generates topologically a profinite subgroup of
$G \times \Bx / (\vpi,\vpi)^\bbZ$. To see this, we consider the
fields $E_g = F[g] \sub M_n(F)$ and $E_b = F[b] \sub B$. The
element $\gamma$ lies in the subgroup $H \sub  E_g^\times \times
E_b^\times$ which consists of all elements $(\alpha,\beta)$ with
$v(\alpha) - v(\beta) = 0$. The group $H/(\vpi,\vpi)^\bbZ$ is
compact and hence pro-finite. Therefore, there is for any $c>0$
some $N \in \bbZ_{>0}$ such that $\gamma^N$ is congruent to the
identity modulo $\frm_R^c$. Hence, if the admissible ideal $\cI$
of $R^{(j)}_K$ contains $\frm_R^c$, then the ideal

$$\cI' = \cI \cdot \gamma^\sharp(\cI) \cdot \ldots \cdot
(\gamma^\sharp)^N(\cI)$$

\medskip

is admissible again and invariant under $\gamma^\sharp$. The
action of $\gamma$ lifts to the blow-up of $\cI'$, and this
blow-up dominates the blow-up of $\cI$. \\

Write $\partial_A$ instead of $\partial_A M^{(j)}_K$, and let
$A_K$ be the class of $\{0\}$ in $\cS_{K,0}$ (which is the unique
{\it maximal} element of $\cS_K$). We show that there is a
covering $(U_A)_{A \in \cS_K}$ of $\overline{M}^{(j)}_K$ by
constructible subsets $U_A$ with the property that for all $A \in
\cS_K$:

$$\partial_A \sub \bigcup_{A' \prec A} U_A$$

\medskip

and $\gamma(U_A) \cap U_A = \emptyset$ for any $A \neq A_K$. Note
first that the complement of each of the closed subsets

$$\overline{\partial_A} = \bigcup_{A' \prec A} \partial_A$$

\medskip

is closed under specialization, because the support of a point in
an analytic adic space does not change under specialization. This
implies that for each of these subsets there is a descending
family

$$V^{(1)}_A \supset V^{(2)}_A \supset V^{(3)}_A \supset \ldots$$

\medskip

of open quasi-compact neighborhoods of $\overline{\partial_A}$
whose intersection is $\overline{\partial_A}$. Now let $A \neq
A_K$ be a minimal element in $\cS_K$. Then $\partial_A$ is just a
closed point. We have

$$\bigcap_i \left( \gamma(V^{(i)}_A) \cap V^{(i)}_A \right)
= \emptyset \,,$$

\medskip

because $\gamma$ is bijective and $\gamma(\partial_A) \cap
\partial_A = \emptyset$, and if the intersection of constructible
sets in a spectral space is empty, then the intersection of
already finitely many of them is empty. Namely, if we equip the
spectral space with its constructible topology, then it becomes
compact and the constructible subsets become closed. Hence there
is some $i$ such that $\gamma(V^{(i)}_A) \cap V^{(i)}_A = \emptyset$
and we put $U_A = V^{(i)}_A$. We proceed by induction
with respect to the partial ordering on $\cS_K$ and assume that
for a given $A \neq A_K$ we have already constructible subsets
$U_B$ for all $B \prec A$, $B \neq A$, such that

$$\partial_B \sub \bigcup_{B' \prec B} U_B$$

\medskip

and $\gamma(\partial_B) \cap \partial_B = \emptyset$. Next we note
that

$$\bigcap_i \Bigl( V^{(i)}_A - \bigcup_{B \prec A, B \neq A} U_B \Bigr)$$

\medskip

is contained in $\overline{\partial_A} - \bigcup_{B \prec A}
\partial_B = \partial_A$ and hence

$$\bigcap_i \Bigl(
\gamma \bigl( V^{(i)}_A - \bigcup_{B \prec A, B \neq A} U_B \bigr) \cap
\bigl( V^{(i)}_A - \bigcup_{B \prec A, B \neq A} U_B \bigr)
\Bigr) = \emptyset \,.$$

\medskip

By the same reasoning as above, there is some $i$ such that

$$\gamma \bigl( V^{(i)}_A - \bigcup_{B \prec A, B \neq A} U_B \bigr)
\cap \bigl( V^{(i)}_A- \bigcup_{B \prec A, B \neq A} U_B \bigr)
= \emptyset \,,$$

and we put $U_A = V^{(i)}_A - \bigcup_{B \prec A, B \neq A} U_B$.
Thus we have proved the assertion about the covering $(U_A)_A$.\\

We show next that there is a formal model $(\cM^{(j)}_K)'$ such
that the images of $\gamma(U_A)$ and $U_A$ under the
specialization map have empty intersection for $A \neq A_K$. This
is a very general fact. Put $\cX = \cM^{(j)}_K$ and
$X = \overline{M}^{(j)}_K$. Namely, if $Z \sub X$ is a constructible
subset one sees easily that

$$\bigcap_{\cX'} sp_{\cX'}^{-1} \bigl( sp_{\cX'}(Z) \bigr) = Z \,,$$

\medskip

where the intersection is over all models $\cX'$. If now $Z,Z'
\sub X$ are two constructible subsets having empty intersection,
one has

$$\bigcap_{\cX'} sp_{\cX'}^{-1} \bigl(
sp_{\cX'}(Z) \cap sp_{\cX'}(Z') \bigr) = \emptyset \,,$$

\medskip

and by the same reasoning as above, there is a model $\cX'$ such that
$sp_{\cX'}(Z) \cap sp_{\cX'}(Z')$ is empty. If now
$(\cM^{(j)}_K)'$ is a model such that the images of $\gamma(U_A)$
and $U_A$ under the specialization map have empty intersection
(for $A \neq A_K$), then this is the case for any model dominating
this one. Hence there is a model with this property for all $A
\neq A_K$, and we can even find a model $(\cM^{(j)}_K)'$ having
this property and such that the action of $\gamma$ lifts to

$$\gamma': (\cM^{(j)}_K)' \lra (\cM^{(j)}_K)' \,.$$

\medskip

If now $x \in (\cM^{(j)}_K)'$ lies in the image of $\partial
M^{(j)}_K$, it is equal to $sp_{(\cM^{(j)}_K)'}(z)$ for some $z$
in some $U_A$ with $A \neq A_K$. Then

$$\gamma'(x) = sp_{(\cM^{(j)}_K)'}(\gamma(z))
\in sp_{(\cM^{(j)}_K)'}(\gamma(U_A)) \,,$$

\medskip

and hence $x \neq \gamma'(x)$.\\

(ii) Let $\cI \sub R^{(j)}_K$ be the ideal that is blown up to
give the formal scheme $(\cM^{(j)}_K)'$. Let
$\frM^{(j)}_K = \Spec(\frR^{(j)}_K)$, and denote by
$\frm_\frR$ the maximal ideal
corresponding to $\frx_K$. Then $R^{(j)}_K$ is the completion of
$\frR^{(j)}_K$ at $\frm_\frR$. Put

$$\cI_\frR = \cI \cap \frR^{(j)}_K \,,$$

\medskip

and let

$$pr: (\frM^{(j)}_K)' \lra \frM^{(j)}_K$$

\medskip

be the blow-up of $\cI_\frR$ on $\frM^{(j)}_K$. As $\cI$ is
admissible, it contains a power of the maximal ideal $\frm_R$ of
$R^{(j)}_K$, and so $\cI_\frR$ contains a power of $\frm_\frR$. It
is hence invertible outside $\frx_K$, and can actually be
generated by finitely many elements of $\frR^{(j)}_K$ which
generate $\cI$ over $R^{(j)}_K$. It follows that the completion of
$(\frM^{(j)}_K)'$ along the preimage $pr^{-1}(\frx_K)$ is
isomorphic to $(\cM^{(j)}_K)'$. Let
$\widetilde{\frM}^{(j)}_K = \Spec(\widetilde{\frR}^{(j)}_K)$. Put

$$\cI_{\widetilde{\frR}} = \gamma_c^\sharp(\cI_\frR)
\cdot \widetilde{\frR}^{(j)}_K \,.$$

\medskip

Recall that $R^{(j)}_K$ is also the completion of
$\widetilde{\frR}^{(j)}_K$ at the unique closed point
$\tilde{\frx}_K$ lying over $\frx_K$. If we now choose $c$
sufficiently large so that $\cI_{\widetilde{\frR}}$ generates
$\cI$ over $R^{(j)}_K$, the blow-up

$$\tilde{pr}: (\widetilde{\frM}^{(j)}_K)' \lra \widetilde{\frM}^{(j)}_K$$

\medskip

of $\cI_{\widetilde{\frR}}$ on $\widetilde{\frM}^{(j)}_K$ will be
an \'etale neighborhood of $pr^{-1}(\frx_K)$. Hence the completion
of $(\widetilde{\frM}^{(j)}_K)'$ along
$\tilde{pr}^{-1}(\tilde{\frx}_K)$ will be isomorphic to
$(\cM^{(j)}_K)'$. The next assertion finally follows from a very
general fact. Namely, let $\vphi_1, \vphi_2$ be an endomorphism of
a noetherian formal scheme $\cX$, let $\cX' \ra \cX$ be an
admissible blow-up of $\cX$, such that $\vphi_i$ lifts to an
endomorphism $\vphi_i'$ of $\cX'$, $i = 1,2$. Then, for any ideal
$\cJ$ of definition of $\cX'$, there is an ideal of definition
$\cI$ of $\cX$, such that, if $\vphi_1 \equiv \vphi_2$ modulo
$\cI$ then $\vphi_1' \equiv \vphi_2'$ modulo $\cJ$. \\

Finally, if $\hat{\gamma}_c'$ approximates $\gamma'$ sufficiently
well, then the induced morphisms on the special fibre will be the
same, hence $\hat{\gamma}_c'$ has no fixed points on
$sp_{(\cM^{(j)}_K)'}(\partial M^{(j)}_K)$. \hfill $\Box$

\medskip

{\it Remark.} We think that it would also be possible to prove the
statements about the appropriate formal models by applying
\cite{Fu2}, Prop. 2.2.5. Fujiwara however is working in an
algebraic (not only formal) setting, and he supposes that the
ambient algebraic space is proper. Although our formal schemes are
also algebraizable, we preferred to prove the assertions directly.

\bigskip

\subsection{The trace of regular elliptic elements}

\hfill{\space} \newline

The following theorem gives an expression of the trace of pairs
$(g,b^{-1}) \in G \times \Bx$, both regular elliptic, on the cohomology
in terms of the number of fixed points.

\medskip

\begin{thm}\label{trace on cohomology}
Let $g \in G$, $b \in \Bx$ be both regular elliptic elements such
that $v(det(g)) + v(N(b)) = 0$. Then there is an open subgroup $K'
\sub K_0$ such that the following holds: if $K \sub K'$ is an open
subgroup normalized by $g$, then the alternating sum of traces of
the endomorphism induced by $(g,b^{-1})$ on the cohomology

$$tr((g,b^{-1})|H^*_c(M^{(j)}_K)) = \sum_i
(-1)^itr((g,b^{-1})|H^i_c(M^{(j)}_K \times_\hFnr \Fbh,\bbQ_\ell))$$

\medskip

is equal to the number of (intersection theoretic) fixed points of
$(g,b^{-1})$ on $M^{(j)}_K \times_\hFnr \Fbh$, which is finite.
\end{thm}

{\it Proof.} Put $\gamma = (g,b^{-1})$. We choose the following
objects such that the assertions of Prop. \ref{appropriate models}
are fulfilled: subgroups $K \sub K' \sub K_0$, the blow-ups

$$(\cM^{(j)}_K)' \stackrel{pr}{\lra} \cM^{(j)}_K \,, \,\,
(\frM^{(j)}_K)' \stackrel{pr}{\lra} \frM^{(j)}_K \,, \,\,
(\widetilde{\frM}^{(j)}_K)' \stackrel{pr}{\lra}
\widetilde{\frM}^{(j)}_K \,,$$

\medskip

the endomorphism

$$\gamma': (\cM^{(j)}_K)' \lra (\cM^{(j)}_K)' \,,$$

\medskip

and the morphism of schemes

$$\gamma_c': (\widetilde{\frM}^{(j)}_K)' \lra (\frM^{(j)}_K)'$$

which approximates $\gamma'$. Put

$$\frM' := (\frM^{(j)}_K)' \times_\fronr \fro_\Fbh \,,$$

\medskip

and

$$\frM'_\eta := (\frM^{(j)}_K)' \times_\fronr \Fbh \,.$$

\medskip

Let $(\frM')^{an}$ be the non-Archimedean analytic space
associated by Berkovich to the formal scheme $\widehat{\frM'}$,
$\widehat{\frM'}$ being the completion of $\frM'$ along its
special fibre (\cite{Be3}, sec. 1). Let $(\frM')^{ad}$ be the
analytic adic space associated by Huber to the formal scheme
$\widehat{\frM'}$ (\cite{Hu3}, sec. 1.9). It is known that
$(\frM')^{an}$ is the maximal Hausdorff quotient of
$(\frM')^{ad}$. Put $\overline{M} = \overline{M}^{(j)}_K$,
$\partial M = \partial M^{(j)}_K$. Let

$$sp^{ad}: (\frM')^{ad} \lra \widehat{\frM'}_s \,, \,\,
sp^{an}: (\frM')^{an} \lra \widehat{\frM'}_s \,, \,\, sp:
\overline{M} \lra (\cM^{(j)}_K)' $$

\medskip

be the specialization maps, and

$$j: \frM'_\eta \lra \frM' \,, \,\, i: \frM'_s \lra \frM'$$

\medskip

the open respectively closed embeddings. Put

$$M = M^{(j)}_K \times_\hFnr \Fbh \sub (\frM')^{ad}$$

\medskip

and let $M^{an}$ be its image in $(\frM')^{an}$. Let $Y pr^{-1}(\frx_K) \stackrel{i_Y}{\lra} \frM'_s$ be the inclusion.
Then we have the following commutative diagram:

\medskip

$$\xymatrixcolsep{3pc}\xymatrix{
(sp^{ad})^{-1}(\partial M) \ar[dd] \ar[r] &
(sp^{ad})^{-1}(Y) \ar[dd] \ar[rrrd] & & & \\
& & M \ar@{^{(}->}[lu] \ar@{^{(}->}[rr] \ar[rd] \ar[ld] & &
(\frM')^{ad} \ar[d] \\
\partial M \ar[r] & \overline{M} \ar[d]^{sp} & &
M^{an} \ar@{^{(}->}[r] \ar[ld] & (\frM')^{an} \ar[d]^{sp^{an}}\\
 & (\cM^{(j)}_K)' & Y \ar[rr] \ar[l] & & \frM'_s }$$

\medskip

Denote by $\cF = \bbZ/\ell^r\bbZ$ the constant torsion sheaf on the
scheme $\frM'_\eta$, on the analytic spaces $(\frM')^{ad}$ and
$(\frM')^{an}$, as well as on subspaces of these. $M^{an}$ is the
maximal Hausdorff quotient of $M$, and by \cite{Hu3}, Thm. 8.3.5,
the cohomology groups of $M$ and $M^{an}$ are canonically
isomorphic. As $M^{an}$ satisfies Poincar\'e duality, we can
compute the trace on the cohomology instead of the cohomology with
compact support and get the same result. By \cite{Be3}, Cor. 3.5,
there is a natural isomorphism

$$R\Gamma(M^{an},\cF) \simeq R\Gamma(Y,i_Y^*i^*R^+j_*\cF)
\,.$$

\medskip

And by \cite{Be3}, Thm. 4.1, we can choose $c$ large enough so
that the morphism induced by $\gamma_c'$ on
$R\Gamma(Y,i_Y^*i^*Rj_*\cF)$ is equal to the morphism induced by
$\gamma$. As $Y$ is proper and $i_Y^*i^*R^+j_*\cF$ is
constructible of finite tor-dimension, we can use the
Grothendieck-Verdier-Lefschetz trace formula, cf. \cite{SGA5},
Exp. III, Cor. 4.7, to compute the trace of $\gamma_c'$. This
trace is equal to the sum of local terms attached to the connected
components of the fixed point locus of $\gamma_c'$. By Prop.
\ref{appropriate models}, $\gamma'_c$ has no fixed points on the
closed subset $\partial Y$ which is defined to be the preimage of
$sp(\partial M) \sub (\cM^{(j)}_K)'$ under the map $Y \ra
(\cM^{(j)}_K)'$ ($Y$ is the special fibre of $(\cM^{(j)}_K)'$). We
show below in Lemma \ref{separation} that $Y -
\partial Y$ is open in $\frM'_s$. The fixed point locus of
$\gamma_c'$ on $Y$ is therefore open in the fixed point locus of
$\gamma_c'$ on $\frM'_s$. It follows that there is a proper
subscheme $D \sub \frM'$ which is contained and open in the fixed
point locus of the correspondence $\gamma_c'$ and whose special
fibre is equal to the fixed point locus of $\gamma_c'$ on $Y$.
Fujiwara's theorem on the specialization of local terms,
\cite{Fu2}, Prop. 1.7.1, then tells us that

$$\mbox{loc}_{D_s}((\gamma_c')_s, i_Y^*i^*Rj_* \cF) =
\sum_{D' \sub \pi_0(D_\eta)} \mbox{loc}_{D'} ((\gamma_c')_\eta,\cF) \,.$$

\medskip

We will show that the connected components $D'$ of $D_\eta$ are
just points with multiplicity one, and the number of these points
is equal to the number of fixed points of $\gamma$ on $M$,
provided $c$ is large enough. As we are now working only on $M$,
we can work with $\gamma_c$ instead of $\gamma_c'$ (cf. Prop.
\ref{appropriate models}), because the correspondences defined by
$\gamma_c$ and $\gamma_c'$ on the analytic spaces are the same.
$\gamma_c$ induces a morphism of formal schemes

$$\hat{\gamma}_c: \cM^{(j)}_K \lra \cM^{(j)}_K$$

\medskip

and the restriction of the correspondence $\gamma_c$ to $M$ is
just an endomorphism, which we denote by $\gamma_c$ as well. Let

$$\hat{\gamma}_c^\sharp: R^{(j)}_K \lra R^{(j)}_K$$

\medskip

be the induced morphism on complete local rings, let $\frm$ be the
maximal ideal of $R^{(j)}_K$, and let $f_1,\ldots,f_r$ be a system
of generators of $\frm$. Because $\hat{\gamma}_c^\sharp$ maps
$\frm$ to itself, the subsets

$$U_\nu = \{v \in M \midc |f_i|_v^\nu \le |\vpi|_v,
i = 1,\ldots,r \}$$

\medskip

of $M$ are stable under $\gamma_c$ and $\gamma$ as soon as $\nu$
is sufficiently large. $sp^{-1}(Y - \partial Y)$ is a
quasi-compact open subset of $M$, and is hence contained in
$U_\nu$, for $\nu \gg 0$. $\gamma_c$ and $\gamma$ do not have
fixed points on $sp^{-1}(\partial Y) \cap \overline{M}$. By Prop.
\ref{continuity of fixed points}, if $c$ is large enough,
$\gamma_c$ will have only finitely many fixed points on $sp^{-1}(Y
- \partial Y)$, these fixed points are all of multiplicity one,
and their number is equal to the number of fixed points of
$\gamma$. The connected components $D'$ of $D_\eta$ are hence just
points of multiplicity one. By \cite{SGA5}, Exp. III, 4.12, each
local term corresponding to a fixed point is then equal to $1$ (as
an element of $\bbZ/\ell^r\bbZ$). Now we pass to the limit $r \ra
\infty$. Hence we conclude that the trace of $\gamma$, which is
equal to the trace of $\gamma_c'$, is given by the number of fixed
points of $\gamma$ on $M$. \hfill $\Box$

\medskip

\begin{cor}\label{fixed point trace formula}
Let $g \in G$, $b \in \Bx$ be both regular elliptic elements. Then
there is an open subgroup $K' \sub K_0$ such that the following
holds: if $K \sub K'$ is an open subgroup normalized by $g$, then
the alternating sum of traces of the endomorphism induced by
$(g,b^{-1})$ on the cohomology

$$tr((g,b^{-1})|H^*_c(M_K/\vpi^\bbZ)) = \sum_i
(-1)^itr((g,b^{-1})|H^i_c((M_K/\vpi^\bbZ) \times_\hFnr
\Fbh,\bbQ_\ell))$$

is equal to

$$n \cdot \# \{h \in G / \vpi^\bbZ K \,|\, h^{-1} g_b h = g^{-1} \} \,.$$

\medskip
\end{cor}

{\it Proof.} We identify the cohomology of $(M_K/\vpi^\bbZ)
\times_\hFnr \Fbh$ with the direct sum of the cohomology of the
$M^{(j)}_K \times_\hFnr \Fbh$ for $j = 0,\ldots,n-1$. If
$v(det(g)) + v(N(b))$ is not a multiple of $n$ then the cohomology
groups

$$H^i_c(M^{(j)}_K \times_\hFnr \Fbh,\bbQ_\ell)$$

are permuted, none of them is mapped to itself, and  hence the
trace on the cohomology is zero. Similarly, for any $h$ the cosets

$$h^{-1} g_b h \vpi^\bbZ K \,\, \mbox{ and } \,\, g^{-1} \vpi^\bbZ K$$

\medskip

will be distinct. So suppose that $v(det(g)) + v(N(b))$ is $kn$
for some integer $k$. Then we can change $b$ to $b\vpi^{-k}$ and
assume that $v(det(g)) + v(N(b))$ is zero. Then the claimed
identity follows immediately from Thm. \ref{fixedpointtheorem} and
Thm. \ref{trace on cohomology}. \hfill $\Box$

\medskip

\begin{lemma}\label{separation}
Let $\frX$ be a scheme of finite type over $\fronr$, let $Y \sub
\frX_s$ be a closed subscheme and $\cY$ the completion of $\frX$
along $Y$. Let $t(\cY)_a$ be the analytic adic space associated to
$\cY$, and $V(\vpi)_\cY \sub t(\cY)_a$ the subspace of points
whose support contains $\vpi$. Put $Y' = sp_\cY(V(\vpi)_\cY) \sub
Y$, where $sp_\cY: \cY \ra Y$ is the specialization map. Then $Y -
Y'$ is open in $\frX_s$.
\end{lemma}

{\it Proof.} Denote by $\cX = \widehat{\frX}$ the $\vpi$-adic
completion of $\frX$. There is a canonical continuous map $\vphi:
t(\cY) \ra t(\cX)$. To describe it, assume $\frX$ is affine, $\frX
= \Spec(\frR)$ say. Denote by $\widehat{\frR}$ the $\vpi$-adic
completion of $\frR$, by $\fra$ the ideal corresponding to $Y$,
and by $\cR$ the $\fra$-adic completion of $\frR$. Then
$t(\cX) = \Spa(\widehat{\frR},\widehat{\frR})$ and
$t(\cY) = \Spa(\cR,\cR)$.
The canonical continuous morphism $\widehat{\frR} \ra \cR$ induces
the map $\vphi$, which is injective.
Without loss of generality we continue to assume
that $\frX$ is affine. $\vphi$ maps $V(\vpi)_\cY \sub t(\cY)$ to
the corresponding subset $V(\vpi)_\cX \sub t(\cX)$, hence induces
an injection

$$\vphi^{an}: t(\cY) - V(\vpi)_\cY \lra
t(\cX)_a = t(\cX) - V(\vpi)_\cX \,,$$

\medskip

which commutes with the specialization maps
$sp_\cY$ and $sp_\cX$ to give a commutative diagram

$$\xymatrixcolsep{5pc}\xymatrix{
t(\cY) - V(\vpi)_\cY \ar[d]^{sp_\cY} \ar[r]^{\vphi^{an}} &
t(\cX)_a \ar[d]^{sp_\cX} \\
Y \ar@{^{(}->}[r] & \frX_s }$$

\medskip

Let $f_1,\ldots,f_r$ be generators of $\fra$. Then $t(\cY) -
V(\vpi)_\cY$ is the union of the rational subsets

$$U_n := \{v \in t(\cY) \midc
\mbox{ for all }i: |f_i|^n_v \le |\vpi|_v \neq 0
\} \,, $$

\medskip

which are open rational subsets in $t(\cX)_a$ as well. We see in
particular, that $\vphi^{an}$ is an open embedding. As $t(\cX)_a$
is homeomorphic to the projective limit over all admissible
blow-ups, equipped with its projective limit topology, it is
enough to show that the subset $sp_\cX^{-1}(Y - Y')$ is open in
$t(\cX)_a$. Now $sp_\cY^{-1}(Y - Y')$ is open in $t(\cY) -
V(\vpi)_\cY$, so we have to show that $sp_\cX^{-1}(Y - Y')$ is
equal to $\vphi^{an}(sp_\cY^{-1}(Y - Y'))$. Let $v \in t(\cX)_a$
be a point such that $sp_\cX(v) \in Y$. Suppose $v$ does not lie in
$\vphi^{an}(t(\cY) - V(\vpi)_\cY)$. Then $v$ is a continuous
valuation of $\widehat{\frR}$ such that $|f|_v < 1$ for all $f \in
\fra$ but $v$ is not $\fra$-adically continuous. Let $\Gamma_v$ be
the value group of $v$ (supposed to be generated by the non-zero
images of $v$), and let $\Gamma_v' \sub \Gamma_v$ be the largest
convex subgroup such that for all $f \in \fra$ the value $|f|_v$
is cofinal for $\Gamma_v'$ (which means that for any $\delta \in
\Gamma_v'$ there is an $t > 0$ such that $|f|_v^t < \delta$), cf.
\cite{Hu1}, Lemma 2.4. By this Lemma, there is some $f \in \fra$
such that $|f|_v \in \Gamma_v'$. The valuation $w = v|\Gamma_v'$
(cf. \cite{Hu1}, sec. 2) is then $\fra$-adically continuous and
analytic for the $\fra$-adic topology (i.e., its support does not
contain $\fra$). But we have $|\vpi|_w = 0$, as otherwise
$|\vpi|_v$ would be contained in $\Gamma_v'$, contradicting our
assumption that $v$ is not $\fra$-adically continuous. Because
$sp_\cX(v) = sp_\cY(w) \in sp_\cY(V(\vpi)_\cY) = Y'$, we find that
all elements in $sp_\cX^{-1}(Y - Y')$ lie in the image of
$\vphi^{an}$, and hence
$sp_\cX^{-1}(Y - Y') = \vphi^{an}(sp_\cY^{-1}(Y - Y'))$
which is open in $t(\cX)_a$.
\hfill $\Box$

\bigskip
\bigskip

\section{The Jacquet-Langlands correspondence realized on the cohomology}
\label{The Jacquet-Langlands correspondence in the cohomology}

\subsection{The trace on the Euler-Poincar\'e characteristic}

\begin{para}\label{characters of cuspidals}
Let $\pi$ be an irreducible supercuspidal representation of $G$.
By the fundamental result of Bushnell-Kutzko \cite{BK} and Corwin
\cite{Co}, we know that $\pi$ is induced from a finite-dimensional
smooth irreducible representation $\lambda$ of some open subgroup
$K_{\pi} \subset G$ that contains and is compact modulo the centre
of $G$, cf. \cite{BK}, Thm. 8.4.1, for a more precise statement.
Hence we may write

$$\pi = c\mbox{-}Ind_{K_{\pi}}^G (\lambda) = Ind_{K_{\pi}}^G (\lambda)\,,$$

\medskip

where the second equality holds by \cite{Bu}, Thm.1. Here,
$Ind_{K_{\pi}}^G (\lambda)$ is by definition the subspace of smooth
vectors in the space

$$\{f:G \ra \lambda \,|\, {\rm for}\,{\rm all }\, k \in K_\pi, g \in G:
f(kg) = \lambda(k)f(g) \} \,.$$

\medskip

Moreover, the character of $\pi$ is a locally constant function
on the set of elliptic regular elements in $G$ (i.e. those whose
characteristic polynomial is separable and irreducible),
and for such an element $g \in G$ we have

$$ \chi_{\pi} (g) = \sum_{\scriptsize{
\begin{array}{c}
h \in G / K_{\pi} \\
h^{-1} g h \in K_{\pi}
\end{array}}}
\chi_\lambda (h^{-1} g h) \,. $$

\medskip

For regular elliptic $g$ the number of elements $h \in G/K_{\pi}$
such that $h^{-1} g h \in K_{\pi}$ is finite. This formula is
due to Harish-Chandra, proofs can be found in \cite{He1} and
\cite{Sa}. For the rest of this section we fix $\pi, K_{\pi}$, and
$\lambda$ with this property.
\end{para}

\begin{para}\label{Jacquet_Langlands}
We recall the Jacquet-Langlands correspondence. For $\pi$ as
above, the representation $\rho = \cJL(\pi)$ of $\Bx$ that
corresponds to $\pi$ via the Jacquet-Langlands correspondence is
characterized by the following identity. Let $g \in G$ and $b \in
B^{\times}$ be regular elliptic elements with the same
characteristic polynomial. Then the following character relation
holds

$$\chi_{\rho} (b) = (-1)^{n-1} \cdot \chi_{\pi} (g) \,,$$

\medskip

cf. \cite{DKV}, introduction, \cite{Ro}. Thm. 5.8., \cite{Ba}.
\end{para}

For an irreducible supercuspidal representation $\pi$ we know by
\ref{main theorem}, part (i), that the representation
$\Hom_G(H^i_c,\pi)$ is a finite-dimensional smooth representation
of $\Bx$. We consider

$$\Hom_G(H^*_c,\pi) = \sum_i (-1)^i\Hom_G(H^i_c,\pi) \,.$$

\medskip

as an element of the Grothendieck group of admissible
representations of $\Bx$.\\

\begin{thm}\label{JL on Euler characteristic}
Let $\pi$ be an irreducible supercuspidal representation of $G$.
Then, in the Grothendieck group of admissible representations of
$\Bx$ we have:

$$\Hom_G(H^*_c,\pi) = n \cdot (-1)^{n-1} \cJL(\pi)\,,$$

\medskip

\end{thm}

{\it Proof.} In the proof of \ref{main theorem}, part (i), we have
seen that as a representation of $\Bx$

$$\Hom_G(H^i_c,\pi)
= \Hom_G(H^i_c(M_\infty/\vpi^\bbZ),\pi \otimes \zeta) \otimes
\xi^{-1} \,,$$

\medskip

where the character $\zeta$ of $G$ is such that $\vpi \in G$ acts
as the identity on $\pi \otimes \zeta$ (notations as introduced in
\ref{main theorem}). As the Jacquet-Langlands correspondence is
compatible with twisting by characters, we can assume from now on
that $\vpi$ acts as the identity on $\pi$. We consider for any $i
\ge 0$ the admissible representation
$V^i = H^i_c(M_\infty/\vpi^\bbZ)$ of
$G/\vpi^\bbZ \times \Bx/\vpi^\bbZ$.
We put $V^i(\pi) = \Hom_G(V^i,\pi)$.\\

We want to compute the traces of Hecke operators on $\Bx$ on the
$V^i(\pi)$. With the notations as in \ref{characters of cuspidals}
we put

$$f_\pi = \frac{1}{vol(K_\pi/\vpi^\bbZ)}\chi_\lambda \cdot
{\bf 1}_{K_\pi/\vpi^\bbZ} \,,$$

\medskip

where ${\bf 1}_{K_\pi/\vpi^\bbZ}$ is the characteristic function
of $K_\pi/\vpi^\bbZ$. Define for any function $f$ on $G$ or $\Bx$
the function $f^*$ by $f^*(x) = f(x^{-1})$. Then the multiplicity
of $\pi$ in $V^i$ is $tr(f^*_\pi | V^i)$, because

$$tr(f^*_\pi | \pi) = 1 \,,$$

\medskip

and $tr(f^*_\pi | V) = 0$ for any admissible representation $V$
not containing $\pi$. It follows that for any compactly supported
function $f$ on $\Bx /\vpi^\bbZ$ one has

$$tr(f^*_\pi \cdot f^* | V^i) = tr(f | V^i(\pi)) \,.$$

\medskip

Let $K \sub K_0$ be an open subgroup normalized by $K_\pi$ and
such that $K$ lies in the kernel of $\lambda$. Then the trace of
$f_\pi$ on $V^i$ is the same as the trace of $f_\pi$ on $H^i_c(M_K
/\vpi^\bbZ)$. Let $W^i$ be the $(G/\vpi^\bbZ \times
\Bx/\vpi^\bbZ)$-subrepresentation of $V^i$ generated by $H^i_c(M_K
/\vpi^\bbZ)$. By \cite{Cas}, Thm. 6.3.10, the representation $W^i$
is of finite length, and the trace of $f^*_\pi \cdot f^*$ on $W^i$
is the same as the trace of that function on $V^i$. For the
following arguments we fix some isomorphism $\Qlb \simeq \bbC$,
and consider all occurring representations by base-change as
representations over $\bbC$. As $W^i$ is of finite length, its
character $\chi_{W^i}$ is a locally integrable function on
$G/\vpi^\bbZ \times \Bx/\vpi^\bbZ$ which is locally constant on
the subset of regular elements (cf. \cite{Le}, Thm. 5.2.4).
Further, by \cite{Ka}, Thm. A, the orbital integrals of $f_\pi$
over regular non-elliptic conjugacy classes all vanish (by
\cite{Le}, Kazhdan's results hold also in the equal-characteristic
case). We get therefore:

$$tr(f^*_\pi \cdot f^* | W^i) = \int_{(G^e/\vpi^\bbZ) \times
(\Bx /\vpi^\bbZ)} \chi_{W^i}(g,b)f^*_\pi(g)f^*(b)dgdb \,,$$

\medskip

where $G^e$ denotes the open set of regular elliptic elements in
$G$. Note that $\supp(f_\pi) \cap G^e/\vpi^\bbZ$ is an open, but
in general not a compact subset of $G/\vpi^\bbZ$. Let

$$C_1 \sub C_2 \sub \ldots \sub G^e/\vpi^\bbZ$$

\medskip

be an ascending sequence of compact-open subsets whose union is
$\supp(f_\pi) \cap G^e/\vpi^\bbZ$. Let $(f_j)$ be a sequence of
locally constant functions on $G$, such that $\supp(f_j) \sub C_j$
and $f_j|_{C_j} = f_\pi|_{C_j}$. Then we have

$$tr(f^*_\pi \cdot f^* | W^i) = \lim_{j \ra \infty}
\int_{(G/\vpi^\bbZ) \times (\Bx/\vpi^\bbZ)}
\chi_{W^i}(g,b)f^*_j(g)f^*(b)dgdb \,.$$

\medskip

(The only reason we switched to $\bbC$ is to make sense of this limit.)
The next step is to compute the integrals appearing in the formula
above. Fix $j$. For a given element $g \in C_j$ there is a compact
open subgroup $K_g$ of $G/\vpi^\bbZ$ such that:\\

\begin{itemize}
\item $g^{-1}K_gg = K_g$
\item $K_g g \sub C_j$, and
\item $f^*_j$
is constant on the coset $K_g g$ \,.
\end{itemize}

\medskip

$K_g$ can be taken to be arbitrarily small (i.e. contained in any
given open subgroup). As $C_j$ is compact, there are finitely many
$g_1,g_2,\ldots,g_{r}  \in C_j$ and corresponding compact-open
subgroups $K_{g_\nu}$ such that $C_j$ is covered by the cosets
$K_{g_\nu} g_\nu$. We may assume that none of these cosets is
contained in another coset. Then we can choose a normal
compact-open subgroup $K'_{g_1}$ of $K_{g_1}$ which is contained
in all the other subgroups $K_{g_\nu}$ and for which $g_1^{-1}
K'_{g_1} g_1 = K'_{g_1}$. Then there are finitely many cosets
$K'_{g_1} g_{1,\al} \sub K_{g_1} g_1$ whose intersection with any
of the cosets $K_{g_\nu} g_\nu$, $\nu > 1$, is empty, and whose
union is the complement of the union of the $K_{g_\nu} g_\nu$,
$\nu > 1$, in $C_j$. Proceeding inductively we then can actually
assume that the cosets $K_{g_\nu} g_\nu$ are mutually disjoint.
Remember that the subrepresentation $W^i$ depends on some compact
open subgroup $K$. Chose $K$ to be contained in all $K_{g_\nu}$.
Then we have

$$\begin{array}{ll}
& \sum_i (-1)^i \int_{(G/\vpi^\bbZ) \times (\Bx /\vpi^\bbZ)}
\chi_{W^i}(g,b)f^*_j(g)f^*(b)dgdb \\
 & \\
= &  \sum_i (-1)^i \sum_{1 \le \nu \le r} \int_{(K_{g_\nu} g_\nu)
\times (\Bx /\vpi^\bbZ)}
\chi_{W^i}(g,b^{-1})f^*_j(g)f(b)dgdb \\
 & \\
= & \sum_{1 \le \nu \le r} \int_{(K_{g_\nu} g_\nu) \times (\Bx
/\vpi^\bbZ)}f^*_j(g)tr((g,b^{-1}) |
H^*_c(M_{K_{g_\nu}}/\vpi^\bbZ))f(b)dgdb
\end{array}$$

\medskip

Now suppose that $\supp(f)$ is contained in the set of regular
elliptic elements of $\Bx$. Then we can use \ref{fixed point trace
formula} which states that

$$tr((g,b^{-1}) | H^*_c(M_{K_{g_\nu}}/\vpi^\bbZ))
= n \cdot \# \{h \in G / \varpi^{\bbZ}K_{g_\nu} \,|\, h^{-1} g_b h
= g^{-1} \} \,,$$

\medskip

where $g_b$ is any element in $G$ having the same characteristic
polynomial as $b$. Plugging this trace formula into the previously
derived expression gives:

$$\begin{array}{ll}
 & \sum_{1 \le \nu \le r}
\int_{(K_{g_\nu} g_\nu) \times \Bx}f^*_j(g)tr((g,b^{-1}) |
H^*_c(M_{K_{g_\nu}}/\vpi^\bbZ))f(b)dgdb \\
 & \\
= & \sum_{1 \le \nu \le r} n \cdot vol(K_{g_\nu}) \int_\Bx
f_j(g^{-1}_\nu) \# \{h \in G / \varpi^{\bbZ}K_{g_\nu} \,|\, h^{-1}
g_b h = g^{-1}_\nu \}
f(b)db \\
 & \\
= & \sum_{1 \le \nu \le r} n \cdot \int_\Bx \int_{G/\varpi^{\bbZ}}
f_{j,\nu,!}(h^{-1} g_b h) f(b)dhdb \\
& \\
= & n \cdot \int_\Bx \int_{G/\varpi^{\bbZ}} f_j(h^{-1} g_b h)
f(b)dhdb \,,
\end{array}$$

\medskip

where $f_{j,\nu,!}$ denotes the extension by zero of the function
$f_j|_{g^{-1}_\nu K_{g_\nu}}$. As $f$ has compact support, we can
pass to the limit as $j$ tends to infinity, and take the
alternating sum over all $i$ to get:

$$\begin{array}{ll}
 & tr(f \midc \Hom_G(H^*_c,\pi))
 = \sum_i (-1)^i tr(f^*_\pi \cdot f^* \midc V^i) \\
 & \\
= & \lim_{j \ra \infty} n \cdot \int_\Bx \left(
\int_{G/\varpi^{\bbZ}}
f_j(h^{-1} g_b h)dg \right) f(b)db \\
 & \\
= & n \cdot \int_\Bx \left( \int_{G/\varpi^{\bbZ}}
f_\pi(h^{-1} g_b h)dg \right) f(b)db \\
 & \\
= & n \cdot \int_\Bx \left( vol(K_\pi/\vpi^\bbZ)\int_{G/K_\pi}
f_\pi(h^{-1} g_b h)dg \right) f(b)db \\
 & \\
= & n \cdot \int_\Bx \left(\sum_{h \in G/K_\pi, h^{-1}g_bh \in
K_\pi}
\chi_\lambda(h^{-1} g_b h)\right) f(b)db \\
 & \\
= & n \cdot \int_\Bx \chi_\pi(g_b)f(b)db \,.
\end{array}$$\\

\medskip

If we now fix a regular elliptic element $b \in \Bx$, and if we
replace $f$ by a sequence of compactly supported functions on
$\Bx$ whose support converges to $\{ b \}$ and whose integral is
$1$, we get

$$tr(b \midc \Hom_G(H^*_c,\pi)) = n \cdot \chi_\pi(g_b)
= n \cdot (-1)^{n-1} \chi_{{\cJL}(\pi)}(b) \,,$$ \\

where the last equality is the character identity of the
Jacquet-Langlands correspondence. Because a virtual representation
of $\Bx$ is already determined by the restriction of its character
to the dense subset of regular elliptic elements, the theorem is
proved. \hfill $\Box$

\bigskip

\subsection{The $\vpi$-adic boundary}
\label{vpi-adic boundary}

\begin{para}\label{set-up middle-degree}
Let $\pi$ be an irreducible supercuspidal representation of $G$.
The aim of the next section is to show that $\Hom_G(H^i_c,\pi) = 0$
for $i \neq n-1$. As we have shown in the proof of \ref{main
theorem}, part (i), we can assume that $\vpi$ acts on $\pi$ as the
identity. Therefore, we consider the admissible representation

$$H^i_c(M_\infty/\vpi^\bbZ) = \lim_{\stackrel{\lra}{K}}
H^i_c(M_K / \vpi^\bbZ) \,,$$

\medskip

with the notation as in the proof of \ref{main theorem}, and show
that no subquotient of this representation is supercuspidal. To
this end we will introduce in this section yet another kind of
compactifications of the spaces $M^{(j)}_K$, to be denoted by
$\overline{M}^{(\vpi,j)}_K$. $\overline{M}^{(\vpi,j)}_K$ and the
boundary

$$\partial^{\vpi} M^{(j)}_K = \overline{M}^{(\vpi,j)}_K -
M^{(j)}_K \times_\hFnr \Fbh$$

\medskip

are analytic pseudo-adic spaces, with $\vpi$ being invertible on
the structure sheaf. We call $\partial^{\vpi} M^{(j)}_K$ the {\it
$\vpi$-adic boundary}.
\end{para}

\medskip

\begin{para}\label{vpi-adic comp}
Let $\frM^{(j)}_K = \Spec(\frR^{(j)}_K)$ be one of the affine
schemes of finite type over $\fronr$ from Thm. \ref{algebraicity}.
There is a closed point $\frx_K$ in the special fibre of
$\frM^{(j)}_K$ such that the completion of $\frM^{(j)}_K$ at
$\frx_K$ is isomorphic to $\cM^{(j)}_K$. Let
$\widehat{\frM}^{(j)}_K$ be the completion of

$$\frM^{(j)}_K \times_{\Spec(\fronr)} \Spec(\fro_\Fbh)$$

along the closed subscheme where $\vpi$ is zero. Denote by
$(\frM^{(j)}_K)^{ad}$ the analytic adic space associated to
$\widehat{\frM}^{(j)}_K$, cf. \cite{Hu3}, Prop. 1.9.1, and by

$$sp_{\widehat{\frM}^{(j)}_K}:
(\frM^{(j)}_K)^{ad} \lra \widehat{\frM}^{(j)}_K$$

\medskip

the specialization map. Put

$$\overline{M}^{(\vpi,j)}_K =
sp_{\widehat{\frM}^{(j)}_K}^{-1}(\frx_K) \,.$$

\medskip

This is a pseudo-adic subspace of $(\frM^{(j)}_K)^{ad}$, cf.
\cite{Hu3}, sec. 1.10, and Prop. \ref{top prop of vpi-adic
boundary} below. $\overline{M}^{(\vpi,j)}_K$ consists of all
$\vpi$-adically continuous valuations $v$ of the ring
$\frR^{(j)}_K \otimes_\fronr \fro_\Fbh$ such that

$$sp_{\widehat{\frM}^{(j)}_K}(v) =
\{ f \in \frR^{(j)}_K \otimes_\fronr \fro_\Fbh \midc |f|_v < 1
\}$$

\medskip

is equal to the maximal ideal of $\frR^{(j)}_K \otimes_\fronr
\fro_\Fbh$ corresponding to the closed point $\frx_K$. The
canonical map

$$\frR^{(j)}_K \otimes_\fronr \fro_\Fbh
\lra R^{(j)}_K \hat{\otimes}_\fronr \fro_\Fbh$$

\medskip

induces an open embedding

$$M^{(j)}_K \times_\hFnr \Fbh \hookrightarrow
\overline{M}^{(\vpi,j)}_K $$

\medskip

(cf. the proof of Lemma \ref{separation}). Most important for our
following reasoning is that the cohomology of these two spaces
coincides:
\end{para}

\medskip

\begin{prop}\label{coh of vpi-adic comp}
For any $i \ge 0$ and $r \ge 0$ the canonical map

$$H^i(\overline{M}^{(\vpi,j)}_K, \bbZ/\ell^r\bbZ) \lra H^i(M^{(j)}_K
\times_\hFnr \Fbh, \bbZ /\ell^r\bbZ)$$

\medskip

is an isomorphism.
\end{prop}

{\it Proof.} Let

$$(\frM^{(j)}_K \times_\fronr \fro_\Fbh)_s \stackrel{i}{\lra}
\frM^{(j)}_K \times_\fronr \fro_\Fbh \stackrel{j}{\longleftarrow}
(\frM^{(j)}_K \times_\fronr \fro_\Fbh)_\eta$$

\medskip

be the closed (resp. open) embedding of the closed (resp. generic)
fibre.\\

By \cite{Hu3}, Thm. 3.5.15, one has a canonical isomorphism

$$H^i(\overline{M}^{(\vpi,j)}_K, \bbZ/\ell^r\bbZ) \simeq
(R^i j_* (\bbZ/\ell^r\bbZ))_{\frx_K}\,,$$

\medskip

The maximal Hausdorff quotient $(M^{(j)}_K \times_\hFnr
\Fbh)^{an}$ of the adic space $M^{(j)}_K \times_\hFnr \Fbh$ is the
non-Archimedean analytic space that Berkovich associates to the
completion of $\frM^{(j)}_K \times_\fronr \fro_\Fbh$ at $\frx_K$,
cf. \cite{Be3}, sec. 1. By \cite{Be3}, Cor. 3.5, one has a
canonical isomorphism

$$H^i((M^{(j)}_K \times_\hFnr \Fbh)^{an}, \bbZ/\ell^r\bbZ) \simeq
(R^i j_* (\bbZ/\ell^r\bbZ))_{\frx_K} \,.$$

\medskip

Hence the cohomology of $\overline{M}^{(\vpi,j)}_K$ is canonically
isomorphic to the cohomology of $(M^{(j)}_K \times_\hFnr
\Fbh)^{an}$ . By \cite{Hu3}, Thm. 8.3.5, the cohomology of
$(M^{(j)}_K \times_\hFnr \Fbh)^{an}$ is canonically isomorphic to
the cohomology of $M^{(j)}_K \times_\hFnr \Fbh$. \hfill $\Box$

\medskip

{\it Remark.} One can show that $M^{(j)}_K \times_\hFnr \Fbh$ is
the interior of $\overline{M}^{(\vpi,j)}_K$ in
$(\frM^{(j)}_K)^{ad}$. If $F$ has characteristic zero then the
assertion of the preceding proposition is \cite{Hu4}, Cor.
3.9.\\

\begin{para}\label{vpi-adic boundary strata}
Recall from \ref{compactifications} the adic space

$$\overline{M}^{(j)}_K = t(R^{(j)}_K,R^{(j)}_K)_a \,.$$

\medskip

There is a canonical continuous map of topological spaces

$$sp: \overline{M}^{(\vpi,j)}_K \lra \overline{M}^{(j)}_K$$

\medskip

which is defined as follows, cf. \cite{Hu1}, Prop. 2.6. For a
point $v \in \overline{M}^{(\vpi,j)}_K$ with value group
$\Gamma_v$ let $c\Gamma_v(\frm)$ be the largest convex subgroup of
$\Gamma_v$ such that $|f|_v$ is cofinal for $c\Gamma_v(\frm)$ for
all $f$ in the maximal ideal of $\frR^{(j)}_K \otimes_\fronr
\fro_\Fbh$ corresponding to $\frx_K$. Then $v|c\Gamma_v(\frm)$
(cf. \cite{Hu1}, sec. 2) is a valuation of $\frR^{(j)}_K
\otimes_\fronr \fro_\Fbh$ which is continuous for the topology
defined by the maximal ideal of $\frR^{(j)}_K \otimes_\fronr
\fro_\Fbh$ corresponding to $\frx_K$. It extends therefore to a
valuation of $R^{(j)}_K$ which is continuous for the topology
defined by the maximal ideal of $R^{(j)}_K$, and we put

$$sp(v) = v|c\Gamma_v(\frm) \in \overline{M}^{(j)}_K$$

\medskip

If the valuation $v$ is in $M^{(j)}_K \times_\hFnr \Fbh$, then $v$
is continuous with respect to the topology defined by the maximal
ideal of $\frR^{(j)}_K \otimes_\fronr \fro_\Fbh$ corresponding to
$\frx_K$, and it extends then to $R^{(j)}_K$. By \cite{Hu1}, Prop.
2.6, the map $sp$ is continuous. ($sp$ is the restriction of the
map $r$ in \cite{Hu3}, Prop. 2.6.) We now consider the preimages
of the strata defined in sec. \ref{compactifications}. For an
element $A \in \cS_K$ put

$$\partial^{\vpi}_A M^{(j)}_K = sp^{-1}(\partial_A
M^{(j)}_K) \,,$$

\medskip

where $\partial^{\vpi}_A M^{(j)}_K$ is defined as in
\ref{compactifications}. Put for $h \in \{0,\ldots,n-1\}$

$$\partial^{\vpi}_h M^{(j)}_K
= \bigcup_{A \in \cS_{K,h}} \partial^{\vpi}_A M^{(j)}_K \,\,,
\,\,\,\,\, \partial^{\vpi}_{\ge h} M^{(j)}_K  = \bigcup_{h' \ge h,
A \in \cS_{K,h'}} \partial^{\vpi}_A M^{(j)}_K \,.$$

\medskip

$\partial^{\vpi}_{\ge h} M^{(j)}_K$ is a closed subspace of
$\overline{M}^{(\vpi,j)}_K$ because it is the preimage of all
strata $\partial_A M^{(j)}_K$ with $A \in \cS_{K,h'}$ and $h' \ge
h$, which is a closed subset by \ref{top prop of comp}. Hence we
have a descending sequence of closed subspaces

$$\overline{M}^{(\vpi,j)}_K = \partial^{\vpi}_{\ge 0} M^{(j)}_K
\supset \partial^{\vpi}_{\ge 1} M^{(j)}_K \supset \ldots \supset
\partial^{\vpi}_{\ge n-1} M^{(j)}_K $$

\medskip

with

$$\partial^{\vpi}_{\ge h} M^{(j)}_K -
\partial^{\vpi}_{\ge h+1} M^{(j)}_K = \partial^{\vpi}_{h} M^{(j)}_K \,.$$

\medskip

Denote by

$$j_A: \partial^{\vpi}_A M^{(j)}_K \lra \partial^{\vpi}_{\ge h} M^{(j)}_K \,,
\,\,\, j_h: \partial^{\vpi}_h M^{(j)}_K \lra \partial^{\vpi}_{\ge
h} M^{(j)}_K$$

\medskip

the inclusions, where $A \in \cS_{K,h}$.
\end{para}

\medskip

\begin{prop}\label{top prop of vpi-adic boundary}
(i) The subsets $\partial^{\vpi}_A M^{(j)}_K$, $\partial^{\vpi}_h
M^{(j)}_K$, and $\partial^{\vpi}_{\ge h} M^{(j)}_K$ are
pseudo-adic subspaces of $(\frM^{(j)}_K)^{ad}$.\\

(ii) $\partial^{\vpi}_{\ge h} M^{(j)}_K$ is proper over
$\Spa(\Fbh,\fro_\Fbh)$. $\partial^{\vpi}_A M^{(j)}_K$ and
$\partial^{\vpi}_h M^{(j)}_K$ are partially proper over
$\Spa(\Fbh,\fro_\Fbh)$. If $\cF$ is a sheaf of abelian groups on
$\partial^{\vpi}_A M^{(j)}_K$, $A \in \cS_{K,h}$, then

$$H^i_c(\partial^{\vpi}_A M^{(j)}_K,\cF)
= H^i(\partial^{\vpi}_{\ge h} M^{(j)}_K,(j_A)_!\cF) \,,$$

\medskip

and

$$H^i_c(\partial^{\vpi}_h M^{(j)}_K,\cF)
= H^i(\partial^{\vpi}_{\ge h} M^{(j)}_K,(j_h)_!\cF) \,.$$

\medskip
\end{prop}

{\it Proof.} (i) $\overline{M}^{(\vpi,j)}_K$ is a closed subspace
of $(\frM^{(j)}_K)^{ad}$ and hence convex and pro-constructible in
$(\frM^{(j)}_K)^{ad}$. So
$((\frM^{(j)}_K)^{ad},\overline{M}^{(\vpi,j)}_K)$ is a pseudo-adic
space in the sense of \cite{Hu3}, Def. 1.10.3. From now on we will
always take $(\frM^{(j)}_K)^{ad}$ as the ambient adic space and
call a subset $Z \sub (\frM^{(j)}_K)^{ad}$ a pseudo-adic space if
it is convex and locally pro-constructible in
$(\frM^{(j)}_K)^{ad}$. The map $sp: \overline{M}^{(\vpi,j)}_K \lra
\overline{M}^{(j)}_K$ is not only continuous but even a spectral
map, cf. \cite{Hu1}, Prop. 2.6. The preimage of a locally
pro-constructible set under a spectral map is again locally
pro-constructible, and the subsets $\partial_A M^{(j)}_K$ of
$\overline{M}^{(j)}_K$ are locally closed, hence locally
pro-constructible. $\partial^{\vpi}_{\ge h} M^{(j)}_K$ is closed
and hence pro-constructible and convex. So it is a pseudo-adic
space. Because in an analytic adic space the generalizations of a
point form a chain, cf. \cite{Hu2}, Lemma 1.1.10 (i), a subset $Z
\sub (\frM^{(j)}_K)^{ad}$ which contains all specializations of
$z$ in $(\frM^{(j)}_K)^{ad}$ for any point $z \in Z$ will be
convex. We call such a set closed under specialization. If $z$
lies in $\partial^{\vpi}_A M^{(j)}_K$, and if $z'$ is a
specialization of $z$, then $sp(z')$ is a specialization of
$sp(z)$. The adic space $\overline{M}^{(j)}_K$ is analytic, and
hence all specializations are secondary specializations, cf.
\cite{Hu3}, sec. 1.1.9. The support of a secondary specialization
of a point does not change, and so $sp(z')$ lies in the same
stratum as $sp(z)$. This means in particular that the set
$\partial^{\vpi}_A M^{(j)}_K$ is convex and hence a pseudo-adic
subspace of $(\frM^{(j)}_K)^{ad}$.
The same is then the case for the union $\partial^{\vpi}_h M^{(j)}_K$.\\

(ii) We show first that any valuation ring $(v,V)$ of
$(\frM^{(j)}_K)^{ad}$, in the sense of \cite{Hu3}, Def. 1.3.5,
with support $v \in \overline{M}^{(\vpi,j)}_K$ has a unique center
on $(\frM^{(j)}_K)^{ad}$. By \cite{Hu3}, Lemma 1.3.6, there is at
most one center on $(\frM^{(j)}_K)^{ad}$, and $(v,V)$ has a center
if and only if the image of $R = \frR^{(j)}_K \otimes_\fronr
\fro_\Fbh$ in the residue field $k(v)$ is contained in $V$. $V$ is
by definition contained in $k(v)^+$. Let $\frm$ be the maximal
ideal of $R$ corresponding to $\frx_K$. If $v$ is an element of
$\overline{M}^{(\vpi,j)}_K$ then $|f|_v < 1$ for all $f \in \frm$.
Hence $f(v)$, the image of $f$ in $k(v)$, is contained in the
maximal ideal of $k(v)^+$, and therefore in the maximal ideal of
$V$. If $f$ is any element of $R - \frm$, there is a unit $\alpha
\in (\fro_\Fbh)^\times$ and an element $f_1 \in \frm$ such that $f
= \alpha + f_1$. Then $f(v) = \alpha + f_1(v) \in V$. Hence, if
$v$ is in $\overline{M}^{(\vpi,j)}_K$, and if $(v,V)$ is a
valuation ring of $(\frM^{(j)}_K)^{ad}$, it has a center on
$(\frM^{(j)}_K)^{ad}$. This center is a specialization of $v$ and
has to lie in $\overline{M}^{(\vpi,j)}_K$. $(\frM^{(j)}_K)^{ad}$
is an affinoid adic space and in particular quasi-compact.
$\overline{M}^{(\vpi,j)}_K$ is closed in $(\frM^{(j)}_K)^{ad}$,
and hence quasi-compact too. By the valuative criterion for
properness \cite{Hu3}, Cor. 1.10.21, the space
$\overline{M}^{(\vpi,j)}_K$ is proper over $\Spa(\Fbh,\fro_\Fbh)$.
Any pseudo-adic subspace of $\overline{M}^{(\vpi,j)}_K$ which is
closed under specialization will then also verify the valuative
criterion for partially properness, and if such a subspace is
quasi-compact it will be proper over $\Spa(\Fbh,\fro_\Fbh)$. For
the partially proper spaces the cohomology with compact is defined
as the derived functor of the functor 'sections with compact
support'. The last assertion of (ii) follows immediately from
\cite{Hu3}, Lemma 5.4.2. \hfill $\Box$

\medskip

\begin{para}\label{action on algebraizations}
Next we want to show that the action of $G = GL_n(F)$ on the
spaces $M^{(j)}_K$ extends to the spaces
$\overline{M}^{(\vpi,j)}_K$, respects the strata and induces an
action on the cohomology of the strata. This is not automatic
because we defined $\overline{M}^{(\vpi,j)}_K$ as a subspace
of $(\frM^{(j)}_K)^{ad}$ which depends on the algebraic scheme
$\frM^{(j)}_K$ that we have chosen as in Thm. \ref{algebraicity}.
The choice of such an algebraization is not unique, and so, in
fact, the pseudo-adic spaces $\overline{M}^{(\vpi,j)}_K$ will not
be isomorphic for different choices of algebraizations. However,
by Prop. \ref{coh of vpi-adic comp}, the cohomology of the
$\overline{M}^{(\vpi,j)}_K$, for different algebraizations, is
isomorphic, and in the proof of Thm. \ref{coh of boundary} we will
show that the cohomology of the strata is also isomorphic, for
different algebraizations. But we will need to know that they are
canonically isomorphic. To this end we will use
a compatible system of algebraizations. \\

Let $G^{(0)} \sub G$ be the group of elements for which the
valuation of the determinant $v(det(g))$ is a multiple of $n$. It
is easily seen that with the notation of \ref{set-up
middle-degree}

$$H^i(M_\infty / \vpi^\bbZ) = Ind^G_{G^{(0)}} H^i_c(M^{(0)}_\infty)$$

\medskip

where

$$H^i_c(M^{(0)}_\infty) = \lim_{\stackrel{\lra}{K}}
H^i_c(M^{(0)}_K \times_\hFnr \Fbh, \Qlb) \,.$$

\medskip

Hence it suffices to extend the action of elements of $G^{(0)}$ to
the $\vpi$-adic boundary and its strata, and from now on we will
therefore consider only the 'height-$0$-components' of the various
spaces we defined. Moreover, we will consider only principal
congruence subgroups $K_m$, and instead of the subscript $K_m$ we
will write $m$. Let $\frS^{(0)}_m$ be a scheme of finite type over
$\fronr$ over which there is given a truncated Barsotti-Tate
$\fro$-module $\frX_m$ of level $m$, and a closed point $\fry_m$
such that the completion of $\frS^{(0)}_m$ at $\fry_m$ is
isomorphic to $\cM^{(0)}_{K_0}$ (the universal deformation space
of $\bbX$), and such that $\frX_m$ corresponds to
$X^{univ}[\vpi^m]$ under this isomorphism. Then, for $m' \le m$,
$\frX_m[\vpi^{m'}]$ will be isomorphic to $X^{univ}[\vpi^{m'}]$.
By the uniqueness theorem of algebraizations, cf. \cite{A2}, Thm.
1.7, there is a scheme $\frS'$, a closed point $\fry'$, and a
diagram of \'etale morphisms

$$\xymatrix{
& \frS' \ar[ld]_{\psi'} \ar[rd]^{\psi} &
 \\
\frS^{(0)}_{m'} & & \frS^{(0)}_{m}}$$

\medskip

such that $\psi'(\fry') = \fry_{m'}$ and $\psi(\fry') = \fry_{m}$,
and $(\psi')^*\frX_{m'} \simeq (\psi)^*\frX_{m}[\vpi^{m'}]$. Then,
when we take the pair $(\frS', (\psi)^*\frX_m)$ instead of
$\frS^{(0)}_m$, we find inductively a tower of schemes of finite
type over $\fronr$

$$\frS^{(0)}_0 \longleftarrow \frS^{(0)}_1 \longleftarrow \frS^{(0)}_2
\longleftarrow \ldots$$

\medskip

with closed points $\fry_m \in \frS^{(0)}_m$ in the special fibre.
The morphisms in this sequence are all \'etale and map $\fry_m$ to
$\fry_{m-1}$, and on $\frS^{(0)}_m$ there is a truncated
Barsotti-Tate $\fro$-module $\frX_m$ of level $m$ which has the
property that the pull-back of $\frX_{m'}$ for $m' \le m$ is
$\frX_m[\vpi^{m'}]$. We may clearly assume the schemes
$\frS^{(0)}_m$ to be affine. Then we define, as in the proof of
Thm. \ref{algebraicity}, the scheme $\frM^{(0)}_m$ as the scheme
over $\frS^{(0)}_m$ which parameterizes Drinfeld
level-$m$-structures on $\frX_m$. This gives a tower of affine
schemes of finite type over $\fronr$

$$\frM^{(0)}_0 \longleftarrow \frM^{(0)}_1 \longleftarrow \frM^{(0)}_2
\longleftarrow \ldots$$

\medskip

with closed points $\frx_m \in \frM^{(0)}_m$ such that the
completion of $\frM^{(0)}_m$ at $\frx_m$ is isomorphic to
$\cM^{(0)}_m$.\\

Now we can explain how the action of the group $G^{(0)}$ extends
to the spaces $\overline{M}^{(\vpi,0)}_K$ and that it respects the
subsets $\partial^{\vpi}_A M^{(0)}_K$. Given an element $g \in
G^{(0)}$, we can assume that the valuation of its determinant is
zero, because $\vpi$ acts anyway trivially on
$H^i_c(M_\infty/\vpi^\bbZ)$. Then we define, as in sec. \ref{group
actions}, for any sufficiently large $m$ a morphism

$$g: \frM^{(0)}_m \lra \frM^{(0)}_{m'} \,,$$

\medskip

where $m'$ depends on $m$ and $g$. Let $0 \le m' \le m$ and $r \in
\bbZ$ be integers such that

$$\xymatrixcolsep{5pc}\xymatrix{
\fro^n \sub \vpi^{-r}g\fro^n \sub \vpi^{-(m-m')}\fro^n \,.}$$

\medskip

Then we have an isomorphism

$$g': \vpi^{-m'}\fro^n/\fro^n \stackrel{\vpi^{-r}g}{\lra}
\vpi^{-m}\fro^n/\fro^n \stackrel{\vpi^{m-m'}}{\lra}
\vpi^{-m'}\fro^n/\fro^n$$

\medskip

as in sec. \ref{compactifications}. Let $f:\frT \ra \frS^{(0)}_m$
be a morphism of schemes over $\fronr$, and let $\phi:
\vpi^{-m}\fro^n/\fro^n \ra f^*\frX_m$ be a level-$m$-structure.
Let

$$\phi' = \phi|_{\vpi^{-m'}\fro^n/\fro^n} \circ g':
\vpi^{-m'}\fro^n/\fro^n \lra f^*\frX_m \,.$$

\medskip

$\phi'$ is then a level-$m'$-structure on
$(f^*\frX_m)[\vpi^{m'}]$. If we let $f': \frT \ra \frS^{(0)}_{m'}$
be the composition of $f$ with $\frS^{(0)}_m \ra \frS^{(0)}_{m'}$,
then $(f^*\frX_m)[\vpi^{m'}] = (f')^*\frX_{m'}$, and $\phi'$ is a
level-$m'$-structure on $(f')^*\frX_{m'}$. This defines the
morphism

$$g: \frM^{(0)}_m \lra \frM^{(0)}_{m'} \,,$$

\medskip

which maps the closed point $\frx_m$ to $\frx_{m'}$ and is, by its
very construction, compatible with the action on the completions.
Then there is a commutative diagram

$$\xymatrix{
\cM^{(0)}_m \ar[r] \ar[d]^{g} &  \frM^{(0)}_m \ar[d]^{g} &
(\frM^{(0)}_m)^{ad} \ar[l] \ar[r]^{sp} \ar[d]^{g}
& \frM^{(0)}_{m,s} \ar[d]^{g} \\
\cM^{(0)}_{m'} \ar[r] &  \frM^{(0)}_{m'} & (\frM^{(0)}_{m'})^{ad}
\ar[l] \ar[r]^{sp} & \frM^{(0)}_{m',s} }$$

\medskip

It follows that $\overline{M}^{(\vpi,0)}_m$ is mapped by $g$ to
$\overline{M}^{(\vpi,0)}_{m'}$, and this map is compatible with
specialization maps to the spaces $\overline{M}^{(0)}_m$ and
$\overline{M}^{(0)}_{m'}$, because these were constructed from the
formal schemes $\cM^{(0)}_m$ and $\cM^{(0)}_{m'}$. Hence there is
a commutative diagram

$$\xymatrix{
\overline{M}^{(\vpi,0)}_m \ar[r]^{sp} \ar[d]^g &
\overline{M}^{(0)}_m \ar[d]^{g} \\
\overline{M}^{(\vpi,0)}_{m'} \ar[r]^{sp} &
\overline{M}^{(0)}_{m'}}$$

\medskip

and $g$ maps therefore the subspace $\partial^{\vpi}_A M^{(0)}_m$
to $\partial^{\vpi}_{g^{-1}A} M^{(0)}_{m'}$, for any $A \in \cS_m
:= \cS_{K_m}$. $\partial^{\vpi}_{\ge h} M^{(0)}_m$ is hence mapped
to $\partial^{\vpi}_{\ge h} M^{(0)}_{m'}$, and we get induced maps
on the cohomology (or cohomology with compact support) of these
spaces.
\end{para}

\medskip

\begin{para}\label{finiteness hypothesis}
We define the $\ell$-adic cohomology groups (with compact support)
for the spaces $\partial^{\vpi}_A M^{(0)}_m$ and
$\partial^{\vpi}_{\ge h} M^{(0)}_m$ by firstly taking cohomology
with coefficients $\bbZ/\ell^r\bbZ$, then passing to the limit for $r
\ra \infty$ and finally tensoring with $\Qlb$. Unfortunately, we
are lacking here some basic finiteness results for the cohomology
of such spaces. We believe that these cohomology groups are
finite-dimensional, and the discussion in the next section will
show that otherwise there would occur some very strange phenomena.
Namely, it will be shown that the cohomology of the $\vpi$-adic
boundary $\partial^{\vpi} M^{(0)}_m = \partial^{\vpi}_{\ge 1}
M^{(0)}_m$ is indeed finite-dimensional, and the limit for $m \ra
\infty$ gives an admissible representation which sits in an exact
sequence

$$H^i_c(\partial^{\vpi}_1 M^{(0)}_\infty, \Qlb) \lra
H^i(\partial^{\vpi}_{\ge 1} M^{(0)}_\infty, \Qlb) \lra
H^i(\partial^{\vpi}_{\ge 2} M^{(0)}_\infty, \Qlb) \,,$$

\medskip

where the subscript $\infty$ indicates that the limit is taken for
$m \ra \infty$. The representation on the left can be shown to be
induced from a representation $\rho$ of a parabolic subgroup $P$.
However, if the cohomology of the strata is not
finite-dimensional, the representation $\rho$ would not be
admissible, and the unipotent radical of $P$ may not act
trivially. Similarly, the cohomology group
$H^i(\partial^{\vpi}_{\ge 2} M^{(0)}_\infty, \Qlb)$ is an
extension of subquotients of parabolically induced
representations, which would not be admissible if the cohomology
groups of the strata are not finite-dimensional. It would follow,
that the representation in the middle of the above exact sequence,
which is admissible, is an extension of subquotients of
representations which are themselves not admissible. Theoretically
this is of course possible, but in our situation it seems to us
very unlikely that something like that may happen. So we will be
working from now on under the hypothesis\\

(H) The cohomology groups $H^i_c(\partial^{\vpi}_A M^{(0)}_m,
\Qlb)$, $H^i(\partial^{\vpi}_h M^{(0)}_m, \Qlb)$ and
$H^i(\partial^{\vpi}_{\ge h} M^{(0)}_m, \Qlb)$ are
finite-dimensional $\Qlb$-vector spaces, for any $i$, $m$, $A$ and $h$.\\

We conclude with two remarks on how to prove this hypothesis,
which we conjecture to be true. Firstly, the subspaces
$\partial^{\vpi}_A M^{(0)}_m$ are of relatively simple type, and
we expect that a good theory of sub-constructible subsets and
their \'etale cohomology should imply the finiteness of the
cohomology. This is in fact a conjecture of R. Huber, cf. the
introduction of \cite{Hu4}. We understood that Huber is presently
working on such a theory, and our hypothesis may then be true by
the results of this theory. Secondly, one can also compute the
cohomology groups of the $\partial^{\vpi}_A M^{(0)}_m$ by means of
the morphism

$$sp: \partial^{\vpi}_A M^{(0)}_m \lra \partial_A M^{(0)}_m \,.$$

\medskip

The essential point to show would be that the sheaves
$R^isp_*(\bbZ/\ell^r\bbZ)$ are constructible on $\partial_A
M^{(0)}_m$. The stalk of this sheaf at a geometric point $x \in
\partial_A M^{(0)}_m$ can be identified with the cohomology of the
preimage $sp^{-1}(x)$. But this set has a kind of modular
interpretation. Namely, if the rank of $A$ is $h$, then the
connected part

$$(X^{univ}[\vpi^\infty] \otimes k(x))^\circ \sub
X^{univ}[\vpi^\infty] \otimes k(x)$$

\medskip

is of rank $h$. Suppose now that $h>0$ (otherwise we are in the
interior). The set $sp^{-1}(x)$ then looks like the generic fibre
of the deformation space of $(X^{univ}[\vpi^\infty] \otimes
k(x))^\circ$ together with an level-$m$-structure $\phi: A \ra
(X^{univ}[\vpi^\infty] \otimes k(x))^\circ$. $sp^{-1}(x)$ is a
space over $Frac(W_\fro(k(x)))$, but the field $k(x)$ carries a
non-trivial valuation, and so $Frac(W_\fro(k(x)))$ is a
higher-dimensional local field, equipped with a valuation of rank
at least two. A further analysis of $sp^{-1}(x)$ reveals that it
(or rather its interior) can be thought of as a rigid analytic
space over this higher-dimensional local field. But, as far as I
know, a theory of analytic spaces over higher-dimensional local
fields has not yet been developed.
\end{para}

\medskip

\subsection{Non-cuspidalness outside the middle degree}
\label{non-cuspidality}

\begin{para}
For $h \in \{1,\ldots,n-1\}$ let $F^h \sub F^n$ be the subspace
which is generated by the first $h$ standard basis vectors. For $m
> 0$ put

$$A_{m,h} = (F^h/\fro^n) \cap (\vpi^{-m}\fro^n/\fro^n)
\,\, \sub \, \vpi^{-m}\fro^n/\fro^n \,,$$

\medskip

i.e. $A_{m,h}$ is the submodule of $\vpi^{-m}\fro^n/\fro^n$
generated by the first $h$ standard generators. Let $P_h \sub G$
be the stabilizer of $F^h$ and $P_h^{(0)} = P_h \cap G^{(0)}$, cf.
sec. \ref{action on algebraizations}. Recall the action of $G$ on
the sets $\cS_{m,h} = \cS_{K_m,h}$ of labels of boundary
components. Firstly, $\cS_{m,h}$ can naturally be identified with
$K_m \bksl G /P_h$, and if $g \in G$ is such that $g^{-1}K_mg \sub
K_{m'}$ then the map from $\cS_{m,h}$ to $\cS_{m',h}$ can be
identified with the map

$$K_m \bksl G / P_h \lra K_{m'} \bksl G / P_h \,, \,\,
K_mxP_h \mapsto K_{m'}g^{-1}xP_h \,.$$

\medskip

The element $A_{m,h} \in \cS_{m,h}$ corresponds to the double
coset $K_m \cdot 1 \cdot P_h$, and is therefore mapped to $K_{m'}
\cdot 1 \cdot P_h$, if $g$ is in $P_h$. If $g$ is an element of
$P_h^{(0)}$, and if $g^{-1}K_mg \sub K_{m'}$ then $g$ maps
$\partial^{\vpi}_{A_{m,h}} M^{(0)}_m$ to
$\partial^{\vpi}_{A_{m',h}} M^{(0)}_{m'}$, cf. \ref{action on
algebraizations}. We have therefore an action of $P_h^{(0)}$ on

$$W^{(h,i,0)} = \lim_{\stackrel{\lra}{m}}
H^i_c(\partial^{\vpi}_{A_{m,h}} M^{(0)}_m, \Qlb) \,,$$

\medskip

and we put

$$W^{(h,i)} = Ind^{P_h}_{P_h^{(0)}} \,W^{(h,i,0)} \,.$$

\medskip

Similarly, for any $h \in \{1,\ldots,n-1\}$ there is an action of
$G^{(0)}$ on

$$V^{(h,i,0)} = \lim_{\stackrel{\lra}{m}}
H^i(\partial^{\vpi}_h M^{(0)}_m, \Qlb)  \,\,\, \mbox{ and } \,\,\,
V^{(\ge h,i,0)} = \lim_{\stackrel{\lra}{m}}
H^i(\partial^{\vpi}_{\ge h} M^{(0)}_m, \Qlb)
$$

\medskip

and we put

$$V^{(h,i)} = Ind^G_{G^{(0)}} \,V^{(h,i,0)} \,\,\, \mbox{ and }
\,\,\,  V^{(\ge h,i)} = Ind^G_{G^{(0)}} \,V^{(\ge h,i,0)} \,.$$
\end{para}

\medskip

\begin{thm}\label{coh of boundary}
We assume that the hypothesis (H) in \ref{finiteness hypothesis}
is fulfilled. Then the following assertions do hold.\\

(i) For $h \in \{1,\ldots,n-1\}$ and any $i$ the representation
$W^{(h,i)}$ of $P_h$ is admissible. The action of the unipotent radical
of $P_h$ on $W^{(h,i)}$ is trivial. \\

(ii) For $h \in \{1,\ldots,n-1\}$ and any $i$ the representation
of $G$ on $V^{(h,i)}$ is canonically isomorphic to the
representation $Ind^G_{P_h} W^{(h,i)}$.\\

(iii) For $h \in \{1,\ldots,n-1\}$ and any $i$ the representation
of $G$ on $V^{(\ge h,i)}$ does not have any supercuspidal
representation as a subquotient.\\

(iv) For any $i \neq n-1$ the representation of $G$ on

$$H^i_c(M_\infty/\vpi^\bbZ) := \lim_{\stackrel{\lra}{K}} H^i_c((M_K / \vpi^\bbZ)
\times_\hFnr \Fbh, \Qlb) \,,$$

\medskip

does not have any supercuspidal representation as a subquotient.
\end{thm}

{\it Proof.} (i) Fix $0 < m' \le m$. We want to compute the
invariants of $K_{m'} \cap P_h$ on
$H^i_c(\partial^{\vpi}_{A_{m,h}} M^{(0)}_m, \Qlb)$, and we want to
show that it is equal to the image of the map

$$H^i_c(\partial^{\vpi}_{A_{m',h}} M^{(0)}_{m'}, \Qlb) \lra
H^i_c(\partial^{\vpi}_{A_{m,h}} M^{(0)}_m, \Qlb) \,.$$

\medskip

From now on we will drop the superscript '$(0)$' everywhere in
this part of the proof. So we write $\frM_m$ instead of
$\frM^{(0)}_m$. Let $\frM_m = \Spec(\frR_m)$. Recall that there is
a action of $K_0$ on this scheme, via the level structure, which
is trivial on $K_m$. Put $\frR'_{m'} = (\frR_m)^{K_{m'}}$ and
$\frM'_{m'} = \Spec(\frR'_{m'})$. This scheme classifies Drinfeld
level-$m'$-structures on $\frX_m[\vpi^{m'}]$ over the scheme
$\frS_m$, cf. \ref{action on algebraizations}. Let $\frx'_{m'}$ be
the closed point in $\frM'_{m'}$ which is the image of $\frx_m$
under $\frM_m \ra \frM'_{m'}$. The morphism $pr: \frM_m \ra
\frM_{m'}$ induces a morphism

$$pr': \frM'_{m'} \ra \frM_{m'}$$

\medskip

because the action of $K_{m'}$ on $\frM_{m'}$ is trivial, and this
morphism induces an isomorphism on the completions, even on the
henselizations at the closed points $\frx'_{m'}$ resp.
$\frx_{m'}$. Let $\widehat{\frM}'_{m'}$ be the completion of
$\frM'_{m'} \times_\fronr \fro_\Fbh$ along the subscheme where
$\vpi$ is zero, let $(\frM'_{m'})^{ad}$ be the analytic adic space
associated to $\widehat{\frM}'_{m'}$, and denote by

$$sp_{\widehat{\frM}'_{m'}}: (\frM'_{m'})^{ad} \lra
\widehat{\frM}'_{m'}$$

\medskip

the specialization map. Let $\overline{M'}^{(\vpi)}_{m'}$ the
preimage of $\frx'_{m'}$ under the specialization map. Then we
have again another specialization map

$$sp': \overline{M'}^{(\vpi)}_{m'} \lra \overline{M}_{m'} \,,$$

\medskip

where $\overline{M}_{m'} = \overline{M}^{(0)}_{K_{m'}}$, in the
notation of \ref{compactifications}. For $B \in \cS_{m'}$ define
$\partial^{\vpi}_B M'_{m'} \sub (\frM'_{m'})^{ad}$ as
$(sp')^{-1}(\partial_B M_{m'})$. The morphism $pr': \frM'_{m'} \ra
\frM'_{m'}$ induces a morphism

$$pr': \overline{M'}^{(\vpi)}_{m'} \ra \overline{M}^{(\vpi)}_{m'}$$

\medskip

which maps $\partial^{\vpi}_B M'_{m'}$ to $\partial^{\vpi}_B
M_{m'}$. Our aim is to show that\\

($\ast$) $pr'$ induces an isomorphism $H^i_c(\partial^{\vpi}_B
M_{m'}, \bbZ/\ell^r\bbZ) \ra H^i_c(\partial^{\vpi}_B M'_{m'},
\bbZ/\ell^r\bbZ)$.\\

Define $\partial^{\vpi}_h M'_{m'}$ and $\partial^{\vpi}_{\ge h}
M'_{m'}$ as in sec. \ref{vpi-adic boundary strata}. Because the
cohomology with compact support of $\partial^{\vpi}_h M'_{m'}$ is
the direct sum of the cohomology groups with compact support of
the spaces $\partial^{\vpi}_B M'_{m'}$ for $B \in \cS_{m',h}$ it
suffices to show that the cohomologies of the $\partial^{\vpi}_h
M'_{m'}$ and $\partial^{\vpi}_h M_{m'}$ are isomorphic. Let $\cF$
be the constant sheaf associated to $\bbZ/\ell^r\bbZ$ on any of our
spaces. Consider the diagram with exact rows

$$\xymatrixcolsep{2pc}\xymatrix{
\ldots \ar[r] & H^i_c(\partial^{\vpi}_h M'_{m'},\cF) \ar[r] \ar[d]
& H^i(\partial^{\vpi}_{\ge h} M'_{m'},\cF) \ar[r] \ar[d]
& H^i(\partial^{\vpi}_{\ge h+1} M'_{m'},\cF) \ar[r] \ar[d] & \ldots \\
\ldots \ar[r] & H^i_c(\partial^{\vpi}_h M_{m'},\cF) \ar[r] &
H^i(\partial^{\vpi}_{\ge h} M_{m'},\cF) \ar[r] &
H^i(\partial^{\vpi}_{\ge h+1} M_{m'},\cF) \ar[r] & \ldots }
$$

\medskip

It follows from this diagram that it suffices to show that the
morphism

$$H^i(\partial^{\vpi}_{\ge h} M'_{m'},\cF) \lra
H^i(\partial^{\vpi}_{\ge h} M_{m'},\cF)$$

\medskip

is an isomorphism for any $h$ in $\{1,\ldots,n-1\}$. We will show
that $\partial^{\vpi}_{\ge h} M_{m'}$ has a finite covering by
closed subsets $Z_{h',\nu}$, $h \le h' \le n-1$, which are
pro-special in the sense of \cite{Hu3}, Def. 3.1.6. Then, by the
remark before Lemma 3.1.7 in \cite{Hu3}, the preimage
$Z'_{h',\nu}$ of $Z_{h',\nu}$ under the map

$$pr': \partial^{\vpi}_{\ge h} M'_{m'} \lra \partial^{\vpi}_{\ge h}
M_{m'}$$

\medskip

is pro-special too, and the union of all $Z'_{h',\nu}$ covers
$\partial^{\vpi}_{\ge h} M'_{m'}$. If $pr'$ induces an isomorphism
between the cohomology groups of the pseudo-adic spaces
$Z_{h',\nu}$ and $Z'_{h',\nu}$, then it maps the cohomology of
$\partial^{\vpi}_{\ge h} M_{m'}$ isomorphically to the cohomology
of $\partial^{\vpi}_{\ge h} M'_{m'}$, because the cohomology of
these spaces can be computed by the spectral sequence associated
to the coverings $(Z_{h',\nu})_{h',\nu}$ and
$(Z'_{h',\nu})_{h',\nu}$, cf. \cite{Hu3}, Cor. 2.6.10. \\

Now we will define the subsets $Z_{h',\nu}$. To this end, consider
the universal level-$m'$-structure

$$\phi_{m'}: \vpi^{-m'}\fro^n/\fro^n \lra \frX_{m'} \,,$$

\medskip

and let $\frZ \sub \cO(\frX_{m'})$ be the ideal defined by the
zero section. For $a \in \vpi^{-m'}\fro^n/\fro^n$ let
$\phi_{m'}(a): \frM_{m'} \ra \frX_{m'}$ be the corresponding
section and $\phi_{m'}(a)^\sharp: \cO(\frX_{m'}) \ra \frR_{m'}$
the corresponding ring homomorphism, where $\frM_{m'} \Spec(\frR_{m'})$.
For $B \in \cS_{m'}$ let $\frP_B$ be the ideal
of $\frR_{m'}$ generated by $\phi_{m'}(a)^\sharp(\frZ)$ for all $a
\in B$. Over the completion $R_{m'}$ of $\frR_{m'}$ at $\frx_{m'}$
the level structure $\phi_{m'}$ gives the universal level
structure on $X^{univ}[\vpi^{m'}]$. Therefore, $\frP_B$ generates
the ideal $\frp_B \sub R_{m'}$, cf. Prop. \ref{top prop of comp},
and it follows that $\frp_B \cap \frR_{m'} = \frP_B$. Put

$$\frQ_h = \prod_{B \in \cS_{m',h}} \frP_B \,, \,\,\,
\frq_h = \prod_{B \in \cS_{m',h}} \frp_B \,.$$

\medskip

Then we have an ascending sequence of ideals in $\frR_{m'}$:

$$\frQ_1 \sub \frQ_2 \sub \ldots \sub \frQ_{n-1} \sub
\frm_{\frR_{m'}} \,,$$

\medskip

where the maximal ideal $\frm_{\frR_{m'}}$ corresponds to the
closed point $\frx_{m'}$, and the ideal generated by $\frQ_h$ in
$R_{m'}$ is equal to $\frq_h$. Choose inductively generators
$f_{1,1},\ldots,f_{1,t_1},\ldots, f_{h,1},\ldots,f_{h,t_h}$ of
$\frQ_h$ such that $f_{1,1},\ldots,f_{1,t_1},\ldots,
f_{h-1,1},\ldots,f_{h-1,t_{h-1}}$ generate $\frQ_{h-1}$. Define
for $\nu  = 1,\ldots,t_{h+1}$:

$$Z_{h,\nu} = \overline{M}_{m'} \cap
\bigcap_{r > 0, 1 \le h' \le h, 1 \le \mu \le t_{h'}} \{v \in
(\frM_{m'})^{ad} \midc |\vpi|_v < |f_{h+1,\nu}|^r_v \,, \,\,
|f_{h',\mu}|_v < |f_{h+1,\nu}|^r_v \} \,.$$

\medskip

The subset $\overline{M}_{m'}$ is pro-special, because it is equal
to

$$\{v \in (\frM_{m'})^{ad} \midc \mbox{ for all }
f \in \frm_{\frR_{m'}}: \, |f|_v < 1 \} \,,$$

\medskip

and $\frm_{\frR_{m'}}$ contains $\vpi$, hence is open.
($\overline{M}_{m'}$ is even special because $\frm_{\frR_{m'}}$ is
finitely generated.) The set

$$\{v \in (\frM_{m'})^{ad} \midc |\vpi|_v < |f_{h,\nu}|^r_v \,,
\,\, |f_{h',\mu}|_v < |f_{h,\nu}|^r_v \}$$

\medskip

is also special (because of the condition $|\vpi|_v <
|f_{h,\nu}|^r_v$), and it is closed because its complement is

$$\{v \in (\frM_{m'})^{ad} \midc |f_{h,\nu}|^r_v \le |\vpi|_v \neq 0 \}
\, \cup \, \{v \in (\frM_{m'})^{ad} \midc |f_{h,\nu}|^r_v \le
|f_{h',\mu}|_v \neq 0 \}$$

\medskip

which is an open subset by the very definition of the topology,
cf. \cite{Hu1}, sec. 2. We claim that

$$\partial^{\vpi}_h M_{m'} \, \sub \,
\bigcup_{1 \le \nu \le t_{h+1}} Z_{h,\nu} \, \sub \,
\partial^{\vpi}_{\ge h} M_{m'} \,.$$

\medskip

If $v$ is an element of $Z_{h,\nu}$ then, for any $f\in \frQ_h$,
$|f|_v$ is not in $c\Gamma_v(\frm)$ ($\frm := \frm_{\frR_{m'}}$),
cf. \cite{Hu1}, sec. 2, hence $\frQ_h \sub \supp(sp(v))$, and so
$\frq_h \sub \supp(sp(v))$, which in turn implies that $\frp_B
\sub \supp(sp(v))$ for some $B \in \cS_{m',h}$. Consequently, $v$
lies in $\partial^{\vpi}_{\ge h} M_{m'}$. If $v$ is an element of
$\partial^{\vpi}_h M_{m'}$, then $\supp(sp(v))$ contains $\frq_h$
but does not contain $\frp_B$ for any $B \in \cS_{m',h+1}$, so
$\supp(sp(v))$ does not contain $\frq_{h+1}$, and there is hence
some $\nu$ with $|f_{h+1,\nu}|_{sp(v)} \neq 0$. Then
$|f_{h+1,\nu}|_v$ is in $c\Gamma_v(\frm)$, hence
$|f_{h+1,\nu}|^r_v > |f|_v$ for all $f \in \frQ_h$. So $v$ is in
$Z_{h,\nu}$. Therefore

$$\bigcup_{h' \ge h, 1 \le \nu \le t_{h'+1}} Z_{h',\nu} \, = \,
\partial^{\vpi}_{\ge h} M_{m'} \,.$$

\medskip

The morphism of schemes $\frM'_{m'} \ra \frM_{m'}$ corresponds to
a ring homomorphism $\frR_{m'} \ra \frR'_{m'}$, and the preimage
$Z'_{h,\nu}$ of $Z_{h,\nu}$ is defined by the images of the
functions $f_{h,\nu}$ in $\frR'_{m'}$. By Thm. 3.2.1 in
\cite{Hu3}, the cohomology of $Z_{h,\nu}$ is the same as the
cohomology of the scheme $\Spec(A(Z_{h,\nu}))$, where
$A(Z_{h,\nu})$ is the henselization of the affinoid ring

$$(\frR_{m'} \otimes_\fronr \Fbh, \frR_{m'} \otimes_\fronr \fro_\Fbh) \,,$$

\medskip

along the pro-special subset $Z_{h,\nu}$, cf. \cite{Hu3}, 3.1.12.
From the very definition of the ring $A(Z_{h,\nu})$ it is easily
seen that it depends only on the henzelization of $\frR_{m'}$ at
the maximal ideal $\frm_{\frR_{m'}}$. The ring homomorphism
$\frR_{m'} \ra \frR'_{m'}$ induces an isomorphism on the
henselizations, and therefore it induces an isomorphism
$A(Z_{h,\nu}) \ra A(Z'_{h,\nu})$. So Thm. 3.2.1 in \cite{Hu3}
shows that the map between the cohomology groups of $Z'_{h,\nu}$
and $Z_{h,\nu}$ is an isomorphism. This completes the proof
($\ast$). \\

Let $P_{m,h}$ be the stabilizer of $A_{m,h}$ in $K_0/K_m$. Because
the set $\cS_{K_m,h}$ can be identified with $K_0 / P_{m,h}$, the
group $P_{m,h}$ is the subgroup of $K_0/K_m$ which stabilizes
$\partial^{\vpi}_{A_{m,h}} M_m$. And because $K_{m'}/K_m$ acts
transitively on the geometric fibres of the canonical map

$$\frM_m^{ad} \lra (\frM'_{m'})^{ad}$$

\medskip

the subgroup $(K_{m'}/K_m) \cap P_{m,h}$ acts transitively on the
geometric fibres of

$$pr_{m,m'}: \partial^{\vpi}_{A_{m,h}} M_m \lra
\partial^{\vpi}_{A_{m',h}} M'_{m'} \,,$$

\medskip

This is a quasi-finite morphism, and the functor $(pr_{m,m'})_*$
is therefore exact (\cite{Hu3}, Prop. 2.6.4). $K_{m'}/K_m$ is a
finite $p$-group and has therefore no cohomology on $\ell^r$-torsion
groups. Hence we can conclude

$$\begin{array}{rl}
H^i_c(\partial^{\vpi}_{A_{m,h}} M_m \,, \,\cF)^{(K_{m'}/K_m) \cap
P_{m,h}}
 & = \, H^i_c \left( \partial^{\vpi}_{A_{m',h}} M'_{m'} \,, \, \big(
(pr_{m,m'})_*(\cF) \big)^{(K_{m'}/K_m) \cap P_{m,h}} \right) \\
 & \\
 & = \, H^i_c(\partial^{\vpi}_{A_{m',h}} M'_{m'} \,, \,\cF) =
H^i_c(\partial^{\vpi}_{A_{m',h}} M_{m'} \,, \,\cF)
\end{array}$$

\medskip

Where the last equality holds by ($\ast$). It follows from this
that

$$H^i_c(\partial^{\vpi}_{A_{m,h}} M_m, \Qlb)^{(K_{m'}/K_m) \cap P_{m,h}}
= H^i_c(\partial^{\vpi}_{A_{m',h}} M_{m'}, \Qlb) \,,$$

\medskip

and when we assume (H), then these spaces are finite-dimensional.
This shows that $W^{(h,i)}$ is an admissible representation of
$P_h$. That the unipotent radical acts trivially is a general fact
that was proved and used in this context by P. Boyer, \cite{Bo},
Lemme 13.2.3. \\

(ii) $\partial^{\vpi}_h M^{(0)}_m$ is the disjoint union of the
spaces $\partial^{\vpi}_A M^{(0)}_m$, for $A \in \cS_{m,h}$. The
action of $K_0/K_m$ on the set of labels $\cS_{m,h}$ is transitive
and the stabilizer of $A_{m,h}$ is by definition $P_{m,h}$. Hence
we see that

$$H^i_c(\partial^{\vpi}_h M^{(0)}_m, \Qlb) Ind^{K_0/K_m}_{P_{m,h}}
= H^i_c(\partial^{\vpi}_{A_{m,h}} M^{(0)}_m,
\Qlb) \,.$$

\medskip

When passing to the limit for $m \ra \infty$ we get

$$V^{(h,i,0)} = Ind^{G^{(0)}}_{P^{(0)}} W^{(h,i,0)}$$

\medskip

and then

$$V^{(h,i)} = Ind^G_{G^{(0)}} V^{(h,i,0)} Ind^G_{P^{(0)}} W^{(h,i,0)}
= Ind^G_P Ind^P_{P^{(0)}} W^{(h,i,0)}
= Ind^G_P W^{(h,i)} \,.$$

\medskip

(iii) From the long exact sequences

$$ \ldots \ra H^i_c(\partial^{\vpi}_h M^{(0)}_m, \Qlb) \ra
H^i(\partial^{\vpi}_{\ge h} M^{(0)}_m, \Qlb) \ra
H^i(\partial^{\vpi}_{\ge h+1} M^{(0)}_m, \Qlb) \ra \ldots $$

\medskip

we deduce an exact sequence

$$\ldots \ra V^{(h,i)} \ra
V^{(\ge h,i)} \ra V^{(\ge h+1,i)} \ra \ldots \,,$$

\medskip

where we have put

$$V^{(\ge h,i)} = Ind^G_{G^{(0)}} \left(
\lim_{\stackrel{\lra}{m}} H^i(\partial^{\vpi}_{\ge h} M^{(0)}_m,
\Qlb) \right) \,.$$

\medskip

By definition we have

$$\partial^{\vpi}_{\ge n-1} M^{(0)}_m = \partial^{\vpi}_{n-1} M^{(0)}_m$$

\medskip

and hence $V^{(\ge n-1,i)} = V^{(n-1,i)}$. By descending induction
on $h$, starting with $h = n-1$ and using (i) and (ii), we
conclude that $V^{(\ge h,i)}$ is a successive extension of
parabolically induced representations (with the unipotent radical
acting trivially), hence does
not have a supercuspidal representation as a subquotient.\\

(iv) By Prop. \ref{coh of vpi-adic comp} and Prop. \ref{top prop
of vpi-adic boundary} (ii), there is a long exact sequence

$$ \ldots \ra H^i_c(M^{(0)}_m \times_\hFnr \Fbh, \Qlb) \ra
H^i((M^{(0)}_m \times_\hFnr \Fbh, \Qlb) \ra
H^i(\partial^{\vpi}_{\ge 1} M^{(0)}_m, \Qlb) \ra \ldots \,.$$

\medskip

Put

$$H^i(M_\infty / \vpi^\bbZ) = \lim_{\stackrel{\lra}{m}} H^i((M_m /
\vpi^\bbZ) \times_\hFnr \Fbh, \Qlb) \,.$$

\medskip

Then we have an exact sequence

$$\ldots \ra H^i_c(M_\infty / \vpi^\bbZ) \ra
H^i(M_\infty / \vpi^\bbZ) \ra V^{(\ge 1,i)}\ra \ldots \,.$$

\medskip

If $i < n-1$ the cohomology with compact support $H^i_c(M_\infty /
\vpi^\bbZ)$ vanishes, and if $i > n-1$ the cohomology
$H^i_c(M_\infty / \vpi^\bbZ)$ vanishes, cf. Lemma \ref{finiteness
of cohomology}. Hence, if $i > n-1$, there is in the preceding
sequence a surjection

$$V^{(\ge 1,i-1)} \twoheadrightarrow H^i_c(M_\infty /
\vpi^\bbZ) \,.$$

\medskip

By part (iii) we can conclude that $H^i_c(M_\infty / \vpi^\bbZ)$
does not have a supercuspidal representation as a subquotient if
$i \neq n-1$. \hfill $\Box$

\bigskip

\begin{para} We want to conclude with some remarks concerning the
spaces $\overline{M}^{(\vpi,j)}_K$. In the definition of the space
$\overline{M}^{(\vpi,j)}_K$ we made use of the scheme of finite
type $\frM^{(j)}_K$. However it is possible to consider also
closely related spaces which are defined only in terms of the
formal schemes $\cM^{(j)}_K = \Spf(R^{(j)}_K)$. For the rest of
this paragraph we will suppress the superscript '(j)'. We equip
$R_K$ the $\vpi$-adic topology, and denote this topological ring
by $R^\diamond_K$. To the formal scheme $\Spf(R^\diamond_K)$ there
is an associated analytic adic space
$M^\diamond_K = t(\Spf(R^\diamond_K))_a$ and a specialization map

$$sp: M^\diamond_K \lra \Spf(R^\diamond_K) \,.$$

\medskip

Let $\overline{M}^\diamond_K$ be the preimage of the closed point
of $\Spf(R^\diamond_K)$ under the specialization map. It is in
this way that the author had considered the $\vpi$-adic
'compactifications' for the first time. By taking preimages one
gets subspaces $\partial^\diamond_A M_K \sub
\overline{M}^\diamond_K$, and the action of $GL_n(F)$ naturally
extends to the spaces $\overline{M}^\diamond_K$ and respects the
subspaces $\partial^\diamond_A M_K \sub \overline{M}^\diamond_K$.
Note that we did not use algebraizations to define these spaces,
and we think that it would be desirable to work only with these
spaces, and do without algebraizations. However, there are some
foundational problems one is facing when working with them. First,
one has to show that the fibre product $\overline{M}^\diamond
\times_\hFnr \Fbh$ does exist as a pseudo-adic space. (This is a
case that has not been dealt with so far, and is not covered by
\cite{Hu2} or \cite{Hu3}.) Let us suppose that it is well defined.
Then one has to show that Prop. \ref{coh of vpi-adic comp} is true
for $\overline{M}^\diamond_K \times_\hFnr \Fbh$. We think that
this is true but could not prove it yet. If one can show that
Prop. \ref{coh of vpi-adic comp} continues to hold for
$\overline{M}^\diamond_K \times_\hFnr \Fbh$, then one could work
entirely with the complete local rings $R_K$ and it would not be
necessary to work with algebraizations in section \ref{vpi-adic
boundary}. In the same vein, we think that it would be nice and
also possible to prove a Lefschetz type trace formula for the
space $\overline{M}_K$ studied in sec. \ref{The boundary of the
deformation spaces}, under the assumption that there are no fixed
points on the boundary. Then one could dispense with
algebraizations completely.
\end{para}

\bigskip

\section{Appendix}

\subsection{Algebraicity of the deformation rings in equal characteristic}

\hfill{\space}\newline

If $F$ has positive characteristic, the fact that the deformation
rings $R_m$ are completions of finitely generated
$\fronr$-algebras follows from the following

\begin{prop}\label{pi-action by polynomial in equal char}
If $char(F)>0$ there is a regular system $u_0 = \vpi, u_1,...,
u_{n-1}$ of parameters of $R_0$, a formal $\fro$-module
$\tilde{X}$ over $\frR_0 = \fronr[u_1,...,u_{n-1}] \sub R_0$ and
an isomorphism of formal $\fro$-modules over $\bbF$

$$\tilde{\iota}: \bbX \ra \tilde{X} \times_{\frR_0} \bbF$$

\medskip

such that

\medskip

(i) $(\tilde{X} \times_{\frR_0} R_0,\tilde{\iota})$ is a universal
deformation of $\bbF$ over $R_0$,

(ii) there is a coordinate $T$ on $\tilde{X}$ such that the
multiplication by $\vpi$ on $\tilde{X}$ is given by the polynomial

$$\vpi T + u_1T^q + ... + u_{n-1}T^{q^{n-1}} + T^{q^n} \,.$$

\end{prop}

{\it Proof.} Let $F({\underline u}) = F(u_1,...,u_{n-1})$ be the
so-called universal $\fro$-typical formal $\fro$-module described
in \cite{HG}, sec. 12. It is defined over $\fro[u_1,...,u_{n-1}]$,
but we consider it over $\fronr[[u_1,...,u_{n-1}]]$, and we
identify the universal deformation ring $R_0$ with
$\fronr[[u_1,...,u_{n-1}]]$. We take the universal deformation
$X^{univ}$ to be $F({\underline u})$. $F({\underline u})$ has the
property that

$$F({\underline u})(T_1,T_2) = T_1 + T_2$$

\medskip

and for $j = 1,...,n$:

$$[\vpi]_{F({\underline u})}(T) \equiv
u_jT^{q^j} \,\, \mod \,\, (\vpi,u_1,...,u_{j-1}), \, \deg (q^j+1) \,,$$

\medskip

where $u_n = 1$. The reduction of $F({\underline u})$ modulo the
maximal ideal of $R_0$ is the formal $\fro$-module $\bbX$ with:

$$\bbX(T_1,T_2) = T_1 + T_2 \,, \,\,\,
[\vpi]_{\bbX}(T) = T^{q^n} \,,$$

\medskip

cf. \cite{HG}, (12.5). However, it is not the case that
$[\vpi]_{F({\underline u})}(T)$ is a polynomial in $T$. We define
a deformation $\tilde{X}$ of $\bbX$ over
$\frR_0 = \fronr[u_1,...,u_{n-1}]$ by

$$\tilde{X}(T_1,T_2) = T_1 + T_2 \,, \,\,\,
[\vpi]_{\tilde{X}}(T) = \vpi T + u_1T^q + ... + u_{n-1}T^{q^{n-1}}
+ T^{q^n} \,.$$

\medskip

The reduction of $\tilde{X}$ modulo the ideal
$(\vpi,u_1,...,u_{n-1})$ is $\bbX$, and we take $\tilde{\iota}$ to
be the identity map. By the universal property of $F({\underline
u})$, there are elements $v_1,...,v_{n-1}$ in the maximal ideal
$\frm_{R_0} = (\vpi,u_1,...,u_{n-1})$ of $R_0 = \fronr[[u_1,...,u_{n-1}]]$,
such that there exists an isomorphism

$$\psi: F(v_1,...,v_{n-1}) \stackrel{\simeq}{\lra}
\tilde{X} \times_{\frR_0} R_0 $$

\medskip

of formal $\fro$-modules over $R_0$, and the reduction of $\psi$
modulo $(\vpi,u_1,...,u_{n-1})$ is the identity. We will show that
the elements $\vpi,v_1,...,v_{n-1}$ form a regular system of
parameters. This in turn implies that
$\tilde{X} \times_{\frR_0} R_0$ is a universal deformation too. \\

Write

$$\psi(T) = a_1T  + a_2T^2 + a_3T^3 + ... $$

\medskip

with $a_1 \equiv 1 \mod \frm_{R_0}$ and $a_i \equiv 0 \mod
\frm_{R_0}$ if $i>1$. Consider the identity of power series in
$T$:

\begin{numequation}\label{isomorphism}
\psi([\vpi]_{F({\underline v})}(T)) = [\vpi]_{\tilde{X}}(\psi(T))
\,.
\end{numequation}

Computing modulo $\vpi$ and degree $q+1$ we find that

$$a_1v_1 \equiv u_1a_1^q \mod (\vpi) \,,$$

hence $v_1 = a_1^{q-1}u_1 + \vpi\xi_{1,0}$ with $\xi_{1,0} \in
R_0$. Fix $1 < i < n$ and assume we had already proven that for $j
< i$ there exist $\xi_{j,0},...,\xi_{j,j-1} \in R_0$ such that

$$v_j = a_1^{q^j-1}u_j + \vpi\xi_{j,0} + ... + u_{j-1}\xi_{j,j-1} \,.$$

\medskip

Then we consider equation \ref{isomorphism} modulo
$(\vpi,u_1,...,u_{i-1})$
and degree $q^j + 1$ and find:

$$a_1v_j \equiv u_ja_1^{q^j} \mod (\vpi,u_1,...,u_{i-1}) \,.$$

\medskip

This shows that the map defined by $u_i \mapsto v_i$ induces on
the tangent space $\frm_{R_0}/(\frm_{R_0})^2$ a map which is given
by an upper triangular matrix with respect to the basis
$(\vpi,u_1,...,u_{n-1})$ with units on the diagonal. Hence
$(\vpi,v_1,...,v_{n-1})$
is a regular system of parameters too. $\Box$ \\

From the description of the rings $R_m$ as given in he proof of
\cite{Dr}, Prop. 4.3, which we recalled in the proof of
\ref{representability of def functors}, and the fact that
$X^{univ}$ may be chosen to be defined over
$\fronr[u_1,...,u_{n-1}]$ and such that $[\vpi]_{X^{univ}}(T)$ is
a polynomial it follows immediately that for any there is a
regular system of parameters

$$(\vpi,u_1, \ldots, u_{n-1})$$

\medskip

of $R_0$ with the following properties:\\

(i) If we put $\frR_0 = \fronr[u_1, \ldots, u_{n-1}]$, and
consider it as a subring of $R_0$, the image of universal
level-$m$-structure

$$\phi^{univ}_m: (\vpi^{-m}\fro/\fro)^n \lra \frm_{R_0}$$

\medskip

consists of elements which are integral over $\frR_0$. For $m \ge
0$ we let $\frR_m \sub R_m$ be the subring $\frR_m \sub R_m$ which
is generated by the image of $\phi^{univ}_m$ over $\frR_0$. It is
a free $\frR_0$-module of rank equal to $\# GL_n(\fro/(\vpi^m))$.
The ideal generated by the image of $\phi^{univ}_m$ is the unique
maximal ideal $\frm_m \sub \frR_m$ over $(\vpi,u_1, \ldots,
u_{n-1})$. The $\frm_m$-adic completion of $\frR_m$ is isomorphic
to $R_m$.

(ii) The maximal ideal $\frm_m$ is stable under the action of
$K_0$, and the isomorphism

$$\widehat{\frR_m}
\stackrel{\simeq}{\lra} R_m$$

\medskip

is $K_0$-equivariant. The completion on the left means
$\frm_m$-adic completion. \\

The reason why there is only one maximal ideal of $\frR_m$ over
$(\vpi,u_1, \ldots, u_{n-1})$ is the following: any
$\vpi^m$-torsion point $\xi$ of $X^{univ}$ satisfies an integral
equation of the form

$$\xi^r + a_1\xi^{r-1} + \ldots + a_r = 0$$

\medskip

with all coefficients $a_1,...,a_r$ in $(\vpi,u_1,...,u_{n-1})$.
Hence, any prime ideal over $(\vpi,u_1,\ldots$ $\ldots,u_{n-1})$
contains the ideal generated by all $\vpi^m$-torsion points, which
is clearly a maximal ideal.\\

\bigskip

\subsection{A continuity property of isolated fixed points}

\begin{para}\label{set-up continuity}
In this section $E$ denotes an algebraically closed
non-Archimedean field with a non-trivial valuation, and $\vpi$ is
a non-zero element of the maximal ideal of $\fro_E$. Let $A$ be an
affinoid $E$-algebra, and put $X = \Spa(A,A^\circ)$. Suppose $X$
is smooth over $\Spa(E,\fro_E)$. Let $\vphi:X \ra X$ be an
endomorphism of $X$ over $E$. We will assume that $\vphi$ has only
a finite number of fixed points and that these are all of
multiplicity one. Our aim in this section is to show that any
endomorphism $\psi$ of $X$ which is sufficiently closed to
$\vphi$, in the sense of \cite{Be4}, sec. 6, has also only
finitely many fixed points, that these are all of multiplicity one
and their total number
is equal to the number of fixed points of $\vphi$.\\

To make the relation of being close concrete, we choose an
affinoid generating system $f_1,\ldots,f_r$ of $A$ over $E$. Then,
by Cor. 6.3 of \cite{Be4}, for any $\vep \in \frE(X)$ (with the
notation of \cite{Be4}, sec. 6), there is a $t \in \bbZ_{>0}$ such
that, if

$$\vphi^\sharp(f_i) - \psi^\sharp(f_i) \in \vpi^t A^\circ$$

\medskip

then $d(\vphi,\psi) < \vep$. Here $\vphi^\sharp, \psi^\sharp: A
\ra A$ denote the corresponding ring homomorphisms, and $A^\circ$
is the subring of power bounded elements. By definition, we write
$d(\vphi,\psi) < \vep_t$ if the above relation holds for the fixed
set of affinoid generators and the element $\vpi$.
\end{para}

\medskip

\begin{lemma}
Let $x \in X$ be a fixed point of $\vphi$. Then there are open
affinoid neighborhoods
$U' = \Spa(B',(B')^\circ) \sub U = \Spa(B,B^\circ)$ of $x$,
with the following properties:\\

(i) $\vphi(U') \sub U$, and $x$ is the only fixed point of $\vphi$
in $U'$.\\

(ii) $U$ and $U'$ are isomorphic to polydiscs:

$$B' \simeq E \langle T_1,\ldots,T_n \rangle \,,$$

\medskip

and $B$ corresponds under the ring homomorphism $B \ra B'$ to

$$\{f = \sum_\nu \alpha_\nu T^\nu \in E \langle T_1,\ldots,T_n \rangle \midc
 \lim_{|\nu| \ra \infty} |\alpha_\nu||\vpi|^{-\frac{|\nu|}{\mu}}  0 \}$$

\medskip

for some $\mu \in \bbZ_{>0}$. \\

(iii) There is a $t > 0$ such that any morphism $\psi: X \ra X$
with $d(\vphi,\psi) < \vep_t$ maps $U'$ to $U$ and has a single
fixed point of multiplicity one on $U'$.
\end{lemma}

{\it Proof.} Because $X$ is smooth, there is for any point $x \in
X$ an \'etale map from an open neighborhood of $x$ to an affine
space, \cite{Hu3}, Cor. 1.6.10. By \cite{Be2}, Thm. 3.4.1, if $x$
is $E$-rational, there is even an isomorphism of a neighborhood of
$x$ with a neighborhood of $0$ in an affine space. So we may
assume $U$ is isomorphic to a polydisc. Then there is another
polydisc $U' \sub (U \cap \vphi^{-1}(U))$. This shows (i) and (ii).\\

Put $a_i = \vphi^\sharp(T_i) - T_i \in B'$, $i = 1,\ldots,n$. Let
$I_\vphi \sub B'$ be the ideal generated by $a_1,\ldots,a_n$. Let
$\diag: U \ra U \times_E U$ be the diagonal morphism and
$(1,\vphi): U' \ra U \times_E U$ be the graph of $\vphi$ on $U'$.
Then

$$U \times_{\scriptsize{\mbox{diag}}, U \times_E U, (1,\vphi)} U'
\simeq \Spa(B'/I_\vphi, (B'/I_\vphi)^\circ)$$

\medskip

is the fixed point locus of $\vphi$ on $U'$, hence equal to the
single point $x$, which is of multiplicity one, and so we have
$B'/I_\vphi \simeq E$. Denote by

$$\delta_\vphi = \det \left( \big(\frac{\partial a_i}{\partial T_j}
\big)_{1 \le i,j \le n} \right) \in B' \,,$$

\medskip

the Jacobi determinant. The image of $\delta_\vphi$ in
$B'/I_\vphi$ is invertible, because of our assumption that the
fixed point $x$ be of multiplicity one. By \cite{Be4}, Cor. 6.3,
there is a $t$ such that $\psi^{-1}(U) = \vphi^{-1}(U)$ if
$d(\vphi,\psi) < \vep_t$. Hence $\psi^\sharp$ maps $B$ to $B'$,
and we can consider the elements $b_i = \psi^\sharp(T_i) - T_i \in
B'$. Define $I_\psi$ to be the ideal generated by
$b_1,\ldots,b_n$. By the generalized Krasner lemma, cf. Cor. 1.7.2
in \cite{Hu3} or \cite{Be3}, Thm. 5.1, it is known that if the
element $b_i$ is sufficiently close to the element $a_i$ for all
$i = 1,\ldots,n$, then the image of $\delta_\psi$ (defined as
above but with $a_i$ being replaced by $b_i$) in $B'/I_\psi$ is
invertible and there is a (continuous) isomorphism of $E$-algebras
$B'/I_\psi \stackrel{\simeq}{\lra} B'/I_\vphi$. By increasing $t$
we may assume that this is the case. Then $\psi$ too has exactly
one fixed point on $U'$, which is of multiplicity one. \hfill
$\Box$

\medskip

\begin{prop}\label{continuity of fixed points}
Let $X$ be a smooth affinoid space, and $\vphi: X \ra X$ be as
above, i.e. with finitely many simple fixed points. Then there is
a $t > 0$ such that any morphism $\psi: X \ra X$ with
$d(\vphi,\psi) < \vep_t$ has finitely many fixed points on $X$,
each fixed point is of multiplicity one, and their number is equal
to the number of fixed points of $\vphi$.
\end{prop}

{\it Proof.} Let $f_1,\ldots,f_r$ be the affinoid generating
system of $A$ over $E$ from \ref{set-up continuity}. Then the
fixed points of $\vphi$ are exactly those $x \in X$ with
$|\vphi^\sharp(f_i) - f_i|_x = 0$, $i = 1,\ldots,r$. Let
$x_1,\ldots,x_s$ be the fixed points of $\vphi$. Choose for any
$j$ an open neighborhood $U'_j$ of $x_j$ with the properties of
the lemma above, and such that $U'_j \cap U'_{j'} = \emptyset$ for
$j \neq j'$. Then there is a $k \in \bbZ_{>0}$ such that

$$\max_{1 \le i \le r} \{|\vphi^\sharp(f_i) - f_i|_x\} >
|\vpi|^k_x$$

\medskip

for all $x \not \in \bigcup_{1 \le j \le s} U'_j$. Suppose
$d(\vphi,\psi) < \vep_{2k}$. Then it follows that for all $x \not
\in \bigcup_{1 \le j \le s} U'_j$ one has also

$$\max_{1 \le i \le r} \{|\psi^\sharp(f_i) - f_i|_x\} >
|\vpi|^k_x \,.$$

\medskip

Hence $\psi$ does not have any fixed point on the complement of
the $U'_j$. If now $t \ge 2k$ is sufficiently large, such that
part (iii) of the lemma holds for all fixed points of $\vphi$,
then $\psi$ has on each $U'_j$ a single fixed point of
multiplicity one. This proves the assertion. \hfill $\Box$

\bigskip

\subsection{A picture of the two boundaries}

\hfill{\space}\newline

The following picture is about the boundaries of the space
$M^{(0)}_{K_1}$ of deformations of a formal $\fro$-module of
height $n = 3$ with level-$1$-structure, and the residue field of
$\fro$ is $\bbF_2$. We denote this space here simply by $M_1$. The
characteristic-$\vpi$-boundary $\partial M_1$ has fourteen strata;
seven one-dimensional strata and seven one-point-strata.
$\partial M_1$ is drawn in the center of the picture.
 Let $\phi_i = \phi^{univ}(e_i)$, $i = 1,2,3$, be the images of the
standard basis vectors of $\vpi^{-1}\fro^3/\fro^3$ under the
universal level structure.
The seven one-dimensional strata are

$$\begin{array}{l}
\{|\phi_i|_v = 0\} \,, i = 1,2,3; \\
 \\
\{|\phi_i + \phi_j|_v = 0\} \,, 1 \le i < j \le 3;\\
 \\
\{|\phi_1 + \phi_2 + \phi_3|_v = 0\} \,.
\end{array}$$

\medskip

The last stratum is drawn with a triangle shape.
The strata of the $\vpi$-adic boundary $\partial^\vpi M_1$
are by definition the preimages under
the specialization map

$$sp: \partial^\vpi M_1 \lra \partial M_1 \,.$$

\medskip

If $A \sub \vpi^{-1}\fro^3/\fro^3$ is of rank two over $\fro/(\vpi)$, then
$\partial_A M_1$ is just a point and its preimage $\partial^\vpi_A M_1$
possesses itself a stratification. Namely, the closure of each stratum
$\partial^\vpi_B M_1$, with $B \sub A$ a rank-$1$-submodule, intersects
$\partial^\vpi_A M_1$ in exactly one point, which we have drawn as a
bold black dot.
So, in fact, $\partial^\vpi M_1$ has a natural stratification indexed by
the set of all {\it flags} in $\vpi^{-1}\fro^3/\fro^3$. This is of course
not something that is special to $M_1$. For any $n$ and any $m \ge 1$
the space $\partial^\vpi M_m$ has a natural stratification whose strata are
indexed by flags

$$\cA: \, \vpi^{-m}\fro^n/\fro^n = A_0 \supsetneq A_1 \supsetneq A_2
\ldots \supsetneq A_r = \{0\} \,,$$

\medskip

where each $A_i$ is free over $\fro/(\vpi^m)$ and a direct summand
of $A_{i-1}$. In order to define these refined strata
intrinsically, let

$$\phi^{univ}: \vpi^{-m}\fro^n/\fro^n \lra \frm_{R_m}$$

\bigskip

be the universal level-$m$-structure. For a point $v \in
\partial^\vpi M_m$, we consider the set of $\vpi^m$-torsion points

$$\{ \phi^{univ}(a)_v \midc a \in \vpi^{-m}\fro^n/\fro^n \} =
X^{univ}[\vpi^m](k(v)) \,,$$

\bigskip

where $\phi^{univ}(a)_v$ denotes the image of $\phi^{univ}(a)$ in
the residue field $k(v)$ at $v$. We introduce the following
notation:

$$|\phi^{univ}(a)|_v \ll |\phi^{univ}(a')|_v$$

\bigskip

if and only if for any $r \in \bbZ_{>0}$ one has
$|\phi^{univ}(a)|_v < |\phi^{univ}(a')|_v^r$, and

$$|\phi^{univ}(a)|_v \sim |\phi^{univ}(a')|_v$$

\bigskip

if and only if there are $r, r' \in \bbZ_{>0}$ such that

$$|\phi^{univ}(a)|_v^r < |\phi^{univ}(a')|_v \,\, \mbox{ and }
\,\, |\phi^{univ}(a')|_v^{r'} < |\phi^{univ}(a)|_v \,.$$

\bigskip

Then, for a given $v \in \partial^\vpi M_m$,
put $A_0(v) = \vpi^{-m}\fro^n/\fro^n$
and define inductively, for $i > 0$, the
subgroup $A_i(v) \sub A_{i-1}(v)$ to be the largest subgroup
satisfying

$$\mbox{ for all } a \in A_i(v) \,
\mbox{ and for all } a' \in A_{i-1}(v) - A_i(v):
|\phi^{univ}(a)|_v \ll |\phi^{univ}(a')|_v \,.$$

\bigskip

Then it is not difficult to see that $A_i(v)$ is free over
$\fro/(\vpi^m)$ and a direct summand of $A_{i-1}(v)$. (One reduces
modulo the ideal generated by $\phi^{univ}(a)_v$ for $a \in
A_i(v)$.) This defines a flag $\cA(v)$ in
$\vpi^{-m}\fro^n/\fro^n$, and $\partial^\vpi_\cA M_m$ consists of
all $v$ with $\cA(v) = \cA$. Moreover, the flag $\cA(v)$ defines a
sequence of canonical subgroups

$$X^{univ}[\vpi^m](k(v)) \supsetneq X^{univ}[\vpi^m]_{k(v),1}
\supsetneq \ldots \supsetneq X^{univ}[\vpi^m]_{k(v),r}$$

\medskip

where

$$X^{univ}[\vpi^m]_{k(v),i} = \{ \phi^{univ}(a)_v \midc a \in A_i(v) \}
\sub X^{univ}[\vpi^m](k(v)) \,.$$



\begin{figure}[htbp]
 \begin{center}
  \psfig{file=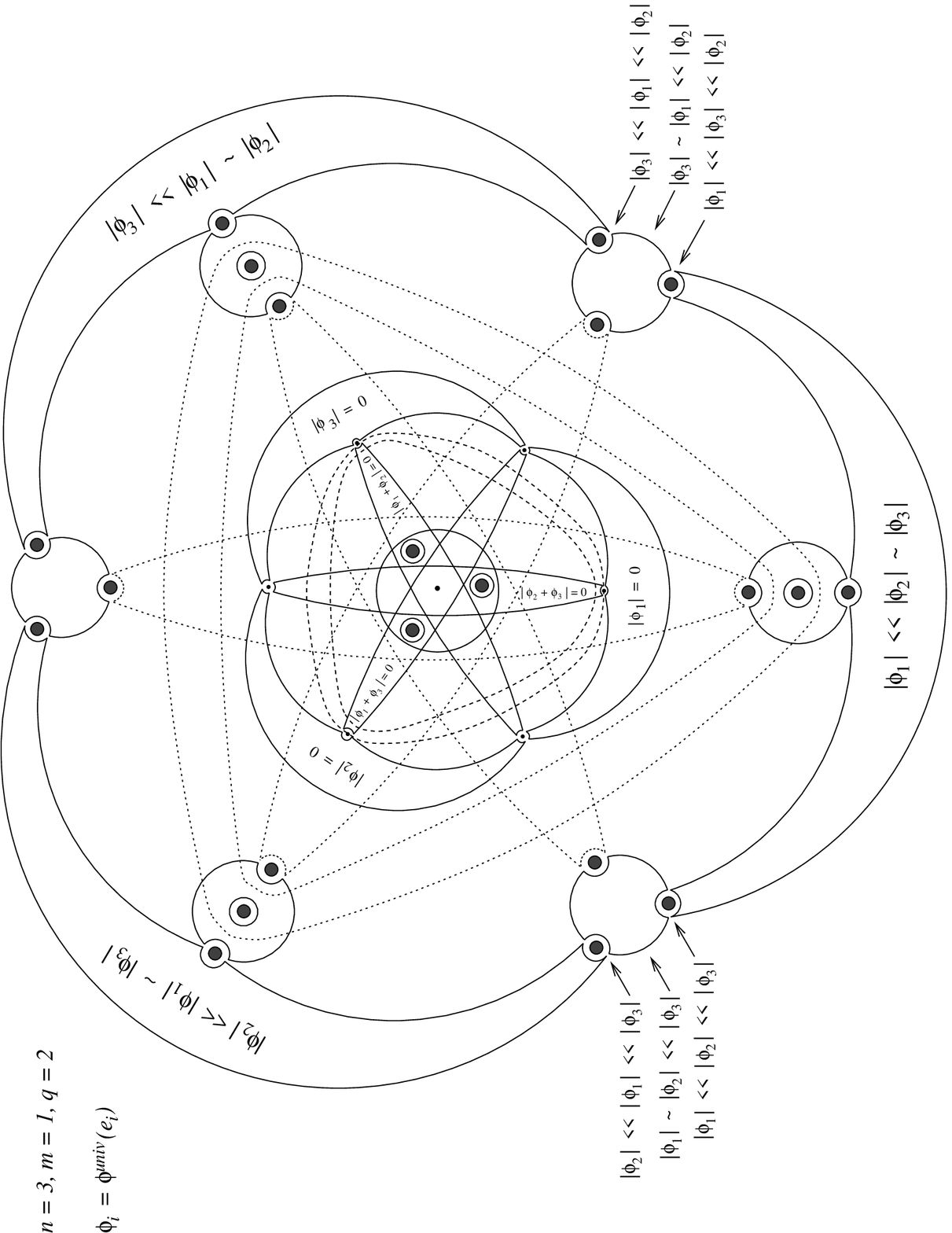,width=15.7cm}
 \end{center}
\end{figure}

\end{document}